\documentclass[12pt,a4paper,english,reqno 
]{amsart}
\usepackage{amsthm}
\usepackage{amssymb}
\usepackage{amscd}
\usepackage{amstext}

\DeclareFontFamily{U}{mathx}{}
\DeclareFontShape{U}{mathx}{m}{n}{<-> mathx10}{}
\DeclareSymbolFont{mathx}{U}{mathx}{m}{n}
\DeclareMathAccent{\widehat}{0}{mathx}{"70}
\DeclareMathAccent{\widecheck}{0}{mathx}{"71}

\usepackage[all]{xy}
\xyoption{2cell}
\UseTwocells

\usepackage{graphicx}
\usepackage{comment}
\usepackage{calrsfs}
\usepackage{framed}

\usepackage[colorlinks=true, allcolors=blue]{hyperref}

\usepackage{tikz}
\usetikzlibrary{positioning,intersections,cd,matrix, arrows, backgrounds}
\usetikzlibrary{shapes,shapes.geometric,shapes.misc}
\usetikzlibrary{decorations}
\usetikzlibrary{decorations.pathreplacing}
\def\Ar#1#2#3{\ar[from={#1}, to={#2}, #3]}
\newcommand\Nname[1]{|[alias=#1]|}

\usepackage[OT2,T1]{fontenc}

\usepackage[colorinlistoftodos]{todonotes}

\theoremstyle{plain}
\newtheorem{thm}{Theorem}[section]
\newtheorem{prp}[thm]{Proposition}
\newtheorem{lem}[thm]{Lemma}
\newtheorem{cor}[thm]{Corollary}
\newtheorem{clm}{Claim}

\newtheorem*{thm-nn}{Theorem}
\newtheorem*{prp-nn}{Proposition}
\newtheorem*{lem-nn}{Lemma}
\newtheorem*{cor-nn}{Corollary}
\newtheorem*{clm-nn}{Claim}
\newtheorem*{cnj-nn}{Conjecture}
\newtheorem*{prb-nn}{Problem}

\theoremstyle{definition}
\newtheorem{dfn}[thm]{Definition}
\newtheorem{exm}[thm]{Example}

\newtheorem*{dfn-nn}{Definition}

\newtheorem{rmk}[thm]{Remark}

\newcommand{\xyR}[1]{%
\xydef@\xymatrixrowsep@{#1}}

\newcommand{\xyC}[1]{%
\xydef@\xymatrixcolsep@{#1}}

\def\al{\alpha}
\def\be{\beta}

\def\de{\delta}
\def\ep{\varepsilon}
\def\ze{\zeta}
\def\et{\eta}
\def\th{\theta}

\def\ka{\kappa}
\def\la{\lambda}
\def\ro{\rho}
\def\si{\sigma}

\def\ph{\phi}
\def\ps{\psi}

\def\De{\Delta}

\def\om{\omega}
\def\Om{\Omega}
\def\Ph{\Phi}

\def\Ker{\operatorname{Ker}}

\def\Im{\operatorname{Im}}

\def\Hom{\operatorname{Hom}}

\def\Aut{\operatorname{Aut}}

\def\Mod{\operatorname{Mod}}

\newcommand{\ltMod}[1]{{}_{#1}\mathrm{Mod}}

\newcommand{\rtMod}[1]{\mathrm{Mod}_{#1}}

\newcommand{\biMod}[2]{{}_{#1}\mathrm{Mod}_{#2}}

\newcommand{\Ginv}{{G}\hyph\mathrm{inv}}
\newcommand{\Ggr}{{G}\hyph\mathrm{gr}}

\def\prj{\operatorname{prj}}

\def\calA{{\mathcal A}}
\def\calB{{\mathcal B}}
\def\calC{{\mathcal C}}
\def\calD{{\mathcal D}}
\def\calE{{\mathcal E}}

\def\bbN{{\mathbb N}}
\def\bbZ{{\mathbb Z}}

\def\op{^{\mathrm{op}}} 

\def\inv{^{-1}}

\def\implies{\text{$\Rightarrow$}}
\def\impliedby{\text{$\Leftarrow$}}

\def\incl{\hookrightarrow}
\def\iso{\cong}
\def\ds{\oplus}
\def\ox{\otimes}

\def\ovl{\overline}
\def\Ds{\bigoplus}
\def\DS{\bigoplus\limits}
\newcommand{\tDs}{\textstyle\Ds}

\def\dsm#1,#2..#3{\bigoplus_{{#1}={#2}}^{#3}}
\def\sm#1,#2..#3{\sum_{{#1}={#2}}^{#3}}

\def\id{1\kern-.25em{\text{{\rm l}}}} 

\def\isoto{\ \raise.3ex\hbox{$^{\sim}$}\kern-.8em\hbox{$\to$}\ } 
\newcommand\Cdot{\raisebox{1pt}{\scalebox{0.4}{$\bullet$}}}
\def\down{_{\Cdot}}

\def\ang#1{{\langle #1 \rangle}}
\def\ya#1{\xrightarrow{#1}}
\def\blank{\operatorname{-}}

\def\bg{%
\family{cmr}\size{20}{12pt}\selectfont}

\def\bigzerou{%
\smash{\lower1.7ex\hbox{\bg 0}}}

\def\repr[#1;#2;#3;#4;#5]{
\left(
\begin{matrix}#1\\#2\end{matrix}
#3
\begin{matrix}#4\\#5\end{matrix}
\right)}
\def\bmat#1{\begin{bmatrix} #1 \end{bmatrix}}
\def\smat#1{\begin{smallmatrix} #1 \end{smallmatrix}}

\def\sbmat#1{\left[\begin{smallmatrix} #1 \end{smallmatrix}\right]}

\def\GCat{G\text{-}\mathbf{Cat}}
\def\GGrCat{G\text{-}\mathbf{GrCat}}

\def\C{\mathcal{C}}

\def\oorbit{/\!_{_{o}}}

\def\GCat{G\text{-}\mathbf{Cat}}

\def\kMod{\ltMod\Bbbk}

\def\k{\Bbbk}

\def\To{\Rightarrow}

\def\hyph{\text{-}}
\newcommand{\bi}[3]{{}_{#2}{#1}_{#3}}
\def\Bi(#1,#2,#3){{}_{#1}{#2}_{#3}}
\def\Biup(#1,#2,#3){{}^{#1}{#2}^{#3}}
\def\Smor(#1,#2,#3){{}^{(#1)}{#2}^{(#3)}}

\newcommand{\sobj}[2]{#1^{(#2)}}  
\newcommand\smor[3]{{}^{(#1)}#2^{(#3)}}  

\newcommand{\Fgt}{\operatorname{Fgt}}

\newcommand{\rme}{\mathrm{e}}

\begin{document}
\title[Cohen-Montgomery duality for bimodules]%
{Cohen-Montgomery duality for bimodules and
singular equivalences of Morita type}
\author{Hideto Asashiba and Shengyong Pan}

\begin{abstract}
Let $G$ be a group and $\k$ a commutative ring.
All categories and functors are assumed to be $\Bbbk$-linear.
We define a $G$-invariant bimodule ${}_SM_R$ over $G$-categories $R, S$ and
a $G$-graded bimodule ${}_BN_A$ over $G$-graded categories $A, B$, and
introduce the orbit bimodule $M/G$ and the smash product bimodule $N\# G$.
We will show that these constructions are inverses to each other.
This will be applied to Morita equivalences,
stable equivalences of Morita type,
singular equivalences of Morita type,
and singular equivalences of Morita type with level
to show that
the orbit (resp.\ smash product) bimodule construction
transforms an equivalent pair of $G$-categories
(resp.\ $G$-graded categories) of each type
to an equivalent pair of $G$-graded categories
(resp.\ $G$-categories) of the same type.
\end{abstract}


\thanks{
HA is partially supported by Grant-in-Aid for Scientific Research 25287001,
25610003, 18K03207, and 25K06922 from JSPS, and by Osaka Central Advanced Mathematical Institute
(MEXT Promotion of Distinctive Joint Research Center Program
JPMXP0723833165).
SP is supported by Beijing Natural Science Foundation (1262017,1252011).
}

\address{Department of Mathematics, Faculty of Science, Shizuoka University,
836 Ohya, Suruga-ku, Shizuoka, 422-8529, Japan;}
\address{
Institute for Advanced Study, KUIAS, Kyoto University,
Yoshida Ushinomiya-cho, Sakyo-ku, Kyoto 606-8501, 
Japan; and}
\address{Osaka Central Advanced Mathematical Institute,
3-3-138 Sugimoto, Sumiyoshi-ku,
Osaka, 558-8585, Japan.}
\email{asashiba.hideto@shizuoka.ac.jp}


\address{School of Mathematics and Statistics, Beijing Jiaotong University,
Beijing, 100044, China.}
\email{shypan@bjtu.edu.cn}
\maketitle

\tableofcontents

\section*{Introduction}
We fix a commutative ring $\k$ and a group $G$.
In most cases, we considered $\k$-linear categories
($\k$-categories for short).
To include infinite coverings of $\k$-algebras into considerations,
we usually regard $\k$-algebras as locally bounded $\k$-categories
with finitely many objects,
and we will work with small $\k$-categories.
A $\k$-category $R$ with a $G$-action $X$ is called a $G$-{\em category},
and is sometimes denoted by $(R, X)$.

The orbit category\footnote{We denote the classical orbit category by $R\oorbit G$ to distinguish it from the orbit category $R/G$ defined in
Definition \ref{dfn:orbit-category}.
In the classical setting, $R\oorbit G$ is isomorphic to any skeleton of $R/G$.  See Remark \ref{rmk:classical-locbdd} for details.} $R\oorbit G$ for a locally finite-dimensional $\k$-category $R$ with a free and locally bounded $G$-action was introduced by Gabriel in \cite{Gab}, which was a central construction of a covering technique in representation theory of algebras, and played an important role to reduce problems on an algebra whose ordinary quiver has oriented cycles to an algebra without oriented cycles in its quiver.
In \cite{Asa97}, we generalized this to treat derived equivalences,
and in \cite{Asa99}, this was used to classify the representation-finite
selfinjective algebras over an algebraically closed field up to derived equivalences.
In \cite{Asa11}, this construction was further generalized to the orbit category $\calC/G$ for any $\k$-linear $G$-categories $\calC$
without any assumptions on $G$-actions
(see also \cite{CM06} and \cite{Ke} for this generalization),
and a notion of $G$-equivariant functors (as 1-morphisms of a 2-category of small $G$-categories) and $G$-precoverings were introduced by using a family of natural isomorphisms.
In \cite{CM06}, the inverse construction of the orbit category was formulated as a smash product $A\# G$ of a $G$-graded category $A$
and the group $G$.  In \cite{Asa11}, we generalized these inverse relations
to those between $G$-categories and $G$-graded categories.
This was formulated as 2-equivalences between 2-categories of
small $G$-categories and of small $G$-graded categories in \cite{Asa18}.

In this paper, we introduce similar constructions on bimodules.
More precisely,
for small $G$-categories $R$ and $S$,
we introduce the notion of $G$-{\em invariant} $S$-$R$-bimodules and their category denoted by
$\Ginv(\biMod{S}{R})$, and
for small $G$-graded $\k$-categories $A$ and $B$, we introduce
$G$-{\em graded} $B$-$A$-bimodules and their category denoted by $\Ggr(\biMod{B}{A})$.
Note here that the orbit category $R/G$ of $R$ by $G$
is a small $G$-graded $\k$-category, and the smash product $A\# G$ of $A$ and $G$ is a small $G$-category.
Then the Cohen-Montgomery duality theorem \cite{CoM84, Asa12} says that we have an equivalence $(R/G)\# G \simeq R$
(resp.\ $(A\# G)/G \simeq A$), by which we  identify the pairs of $G$-invariant (resp.\ of $G$-graded) categories (see also \cite{CM06}).
Here we introduce functors $?/G : \Ginv(\biMod{S}{R}) \to \Ggr(\biMod{S/G}{R/G})$ and
$?\# G : \Ggr(\biMod{B}{A}) \to \Ginv(\biMod{(B\# G)}{(A\# G)})$, and show that they are equivalences
and quasi-inverses to each other (by applying $A:=R/G$, $R:=A\# G$, etc.),
which have good relationships with tensor products and
preserve one-sided and two-sided projectivity of bimodules.
We apply this to stable equivalences of Morita type
(including Morita equivalences)
to have the following theorem (see Theorem \ref{thm:st-eq-M} for details)
in the $\k$-projective case, where a $\k$-category $\calC$ is said
to be $\k$-{\em projective}
if $\calC(x, y)$ are projective $\k$-modules for all objects $x, y$ in $\calC$:

\begin{thm-nn}
Assume that all of $R$, $S$, $A$ and $B$ are $\k$-projective. 
Then the following statements hold.

$(1)$
There exists a ``$G$-invariant stable equivalence of Morita type'' between $R$ and $S$ if and only if
there exists a ``$G$-graded stable equivalence of Morita type'' between $R/G$ and $S/G$.

$(2)$
There exists a ``$G$-graded stable equivalence of Morita type'' between $A$ and $B$ if and only if
there exists a ``$G$-invariant stable equivalence of Morita type'' between $A\# G$ and $B\# G$.
$($See {\rm Definition \ref{dfn:inv-graded-Morita}} for definitions of terminologies in the quotation marks$)$.
\end{thm-nn}

We note that a $G$-invariant (resp.\ $G$-graded) stable equivalence of Morita type
is defined to be a usual stable equivalence of Morita type with additional properties, and does not mean
an equivalence between stable categories of $G$-invariant (resp.\ $G$-graded) modules.
We also give the corresponding results for singular equivalences of Morita type (Theorem \ref{thm:sing-eq-M}),
and for singular equivalences of Morita type with level
(Theorem \ref{thm:sing-eq-M-l}).

The paper is organized as follows.
In Sect.\ 1, we review $G$-categories
(resp.\ $G$-graded categories) and their orbit categories by $G$
(resp.\ its smash product with $G$); and also review
our theorems stating that the orbit category construction and
the smash product construction are extended to
2-equivalences between the 2-categories of $G$-categories
and of $G$-graded categories.

In Sect. 2, we introduce $G$-invariant structure for
bimodules over $G$-categories, and $G$-gradings for
bimodules over $G$-graded categories.
To these two kinds of bimodules,
the two constructions in Sect.\ 1 are exported in sections 3
and 4, respectively.

In Sect.\ 5, we show that the exported constructions
give us an equivalence between the category of $G$-invariant
$S$-$R$-bimodules over $G$-categories $R, S$ 
(resp.\ of $G$-graded $B$-$A$-bimodules over $G$-graded categories $A, B$) and
the category of $G$-graded $S/G$-$R/G$-bimodules
(resp.\ of $G$-invariant $B\# G$-$A\# G$-bimodules).

In Sect.\ 6, we give explicit forms of finitely generated projective $G$-invariant bimodules (resp.\ finitely generated projective $G$-graded bimodules) over locally bounded categories with free $G$-actions
(resp.\ locally bounded $G$-graded categories).
Therefore, within this section, $\k$ is assumed to be a field.

In Sect.\ 7, we investigate further properties of smash products, and in Sect.\ 8, we apply these tools to show
our main theorems on Morita equivalences, stable equivalences
of Morita type, singular equivalences of Morita type,
and singular equivalences of Morita type with level.

In Sect.\ 9, we assume that $\k$ is an algebraically closed field and categories are given by finite bound quivers, and
we give examples of a $G$-invariant bimodule $M$
and a $G$-graded bimodule $M'$ which correspond to each other
by the orbit construction and the smash product construction.
In this example, both $M$ and $M'$ have their counter-parts $N$ and $N'$, respectively, such that the pair $(M, N)$ gives a $G$-invariant stable equivalence of Morita type and the pair $(M', N')$ gives a $G$-graded stable equivalence of Morita type.

\section*{Acknowledgments}
The basic part of this work was done during HA's visits
to Stuttgart in February, 2014 and February, 2015.
HA would like to thank Professor Steffen K{\"o}nig for his kind hospitality and support.
The results until stable equivalences of Morita type 
were announced at several meetings such as
China-Japan-Korea international conference on ring theory
held in China in 2015, 
Symposium on ring theory and representation theory held in Japan in 2015,
a seminar in Beijing Jiaotong University in 2016,
International Conference on Representations of Algebras (ICRA) at Syracuse U.S.A.\ in 2016, and so on.
The corresponding problem about singular equivalences of Morita type
was posed to HA by SP and solved by HA in 2016.
After a long break for about eight years,
precise investigation was done through discussions
using Zoom with SP that started in February, 2024.
In ICRA in 2024 at Shanghai, SP presented the result also on
singular equivalences of Morita type with level.
The example of the bimodule ${}_BM'_A$ in Sect.\ \ref{sec:exmples} is a reformulation of one half of an example of a pair of bimodules inducing stable equivalence of Morita type between Brauer tree algebras of Dynkin type $A_3$ with multiplicity 2 due to Professor Naoko~Kunugi.
She kindly gave HA a copy of the computation in March, 2015, for which the authors are very thankful.
The final computations of the example were done during HA's visit
to Beijing Jiaotong University in November, 2024.
HA is thankful for SP's kind hospitality and support.
The authors would like to thank the referee for valuable comments
especially for pointing out the missing assumption in Proposition 3.6.
\section{Preliminaries}

We denote the category of $\k$-modules by $\kMod$.
Let $\calC$ be a small $\k$-category.
Then we denote the class of objects (resp.\ morphisms) of $\calC$
by $\calC_0$ (resp.\ $\calC_1$).
A covariant (resp.\ contravariant) functor $\calC \to \kMod$ is called
a {\em left} (resp.\ {\em right}) $\calC$-{\em module}, and a natural transformation
between such functors is called a {\em morphism} between them regarded as modules.
The left (resp.\ right) $\calC$-modules and the morphisms between them form a $\k$-category,
which we denote by
$\ltMod\calC$ (resp.\ $\rtMod\calC$).
We set ${}_y\calC_x:= \calC(x, y)$
(resp.\ ${}_g\calC_f:= \calC(f, g)$ for all $f, g \in \calC_1$),
and we regard the category $\calC$
as a $\k$-bilinear functor
\[
\calC \times \calC\op \to \kMod, (y,x) \mapsto {}_y\calC_x.
\]
We also set $\calC_x:= \calC(x, \blank)$ and ${}_x\calC:= \calC(\blank, x)$
for all $x\in \calC_0$.
Then for any morphism $f \colon x \to y$ in $\calC$, note that
$\calC_f \colon \calC_y \to \calC_x$ (resp.\ ${}_f\calC \colon {}_x\calC \to {}_y\calC$)
is a morphism in $\ltMod\calC$ (resp.\ $\rtMod\calC$).

Let $\calD$ also be a small $\k$-category.
Then a $\k$-bilinear functor
\[
M \colon \calD \times \calC\op \to \kMod, (y,x) \mapsto {}_yM_x
\]
is called a $\calD$-$\calC$-{\em bimodule}, and a natural transformation between such functors
is called a {\em morphism} between them regarded as bimodules.
The $\calD$-$\calC$-bimodules and the morphisms between them form a $\k$-category,
which we denote by $\biMod\calD\calC$.
Let $M$ be a $\calD$-$\calC$-bimodule.
We say that ${}_\calD M$ is projective (resp.\ finitely generated projective) if
$M_x$ is a projective (resp.\ finitely generated projective) left $\calD$-module
for all $x \in \calC_0$; and that $M_\calC$ is projective (resp.\ finitely generated projective) if
${}_yM$ is a projective (resp.\ finitely generated projective) right $\calC$-module
for all $y \in \calD_0$.

Recall that the tensor product $\calD \ox_\k \calC$ is defined as follows.
\begin{itemize}
\item
$(\calD \ox_\k \calC)_0:= \calD_0 \times \calC_0$.
\item
For any $(y',x'),\, (y,x) \in (\calD \ox_\k \calC)_0$,
$\bi{(\calD \ox_\k \calC)}{(y',x')}{(y,x)}:= \bi{\calD}{y'}{y}\ox_\k \bi{\calC}{x'}{x}$.
\item
For any $g'\ox f' \in \bi{(\calD \ox_\k \calC)}{(y'',x'')}{(y',x')}$
and $g\ox f \in \bi{(\calD \ox_\k \calC)}{(y',x')}{(y,x)}$,
$(g'\ox f')(g\ox f):= g'g \ox f'f$.
\item
For any $(y,x) \in (\calD \ox_\k \calC)_0$,
$\id_{(y,x)}:= \id_y\ox \id_x$.
\end{itemize}
Then for each $\k$-linear category $\calE$,
we can identify each $\k$-bilinear functor
$F \colon \calD \times \calC \to \calE$ with a $\k$-linear functor
$F' \colon \calD \ox_\k \calC \to \calE$
defined by $F'(y, x):= F(y,x)$,\, $F'(g\ox f):= F(g, f)$
for all $(y,x)\in (\calD \ox_\k \calC)_0$ and
$(g,f) \in (\calD \ox_\k \calC)_1$.
In particular,
each bimodule $M$ above can be identified with a left
$\calD \ox_\k \calC\op$-module.

\subsection{\texorpdfstring{$G$}{G}-categories}

We first recall definitions of $G$-categories and their 2-category $\GCat$.

\begin{dfn}
\label{dfn:G-action}
(1) A $\k$-category with a $G$-\emph{action}, or a $G$-\emph{category} for short,
is a pair $(\mathcal{C}, X)$ of a $\k$-category $\mathcal{C}$
and a group homomorphism $X\colon G\to\Aut(\mathcal{C}), a\mapsto X_a$.
We often write $ax$ for $X_a(x)$ for all $a\in G$ and $x \in \calC_0 \cup \calC_1$
if there seems to be no confusion.

(2) Let $\calC = (\mathcal{C}, X)$ and $\calC' = (\mathcal{C}', X')$ be $G$-categories.
Then a $G$-\emph{equivariant} functor from $\calC$ to $\calC'$ is a pair
$(E,\ka)$ of a $\k$-functor
$E\colon\mathcal{C}\to\mathcal{C}'$ and
a family $\ka=(\ka_{a})_{a\in G}$ of natural isomorphisms
$\ka_{a}\colon X'_{a}E\To EX_{a}$ ($a\in G$)
 such that the diagrams
 \[
\xymatrix{
X'_{ba}E=X'_{b}X'_{a}E\ar@{=>}[r]^{X'_{b}\ka_{a}}\ar@{=>}[rd]_{\ka_{ba}} &
X'_{b}EX_{a}\ar@{=>}[d]^{\ka_{b}X_{a}}\\
 & EX_{ba}=EX_{b}X_{a}}
\]
commute for all $a,b\in G$.
We set
$$
\ka_a = (\ka_{a,x})_{x \in \calC_0}
$$
with $\ka_{a,x} \colon X'_a(E(x)) \to E(X_a(x))$ the $x$-component
of $\ka_a$.

(3) A $\k$-functor $E \colon \calC \to \calC'$ is called
a \emph{strictly} $G$-\emph{equivariant} functor
if  $(E, (\id_E)_{a\in G})$ is a $G$-equivariant functor,
i.e., if $X'_{a}E=EX_{a}$ for all $a\in G$.

(4) Let $(E,\ka),(E',\ka')\colon\mathcal{C}\to\mathcal{C}'$ be $G$-\emph{equivariant}
functors. Then a \emph{morphism} from $(E,\ka)$ to $(E',\ka')$
is a natural transformation $\eta\colon E\To E'$ such that
the diagrams
$$
\xymatrix{
X'_{a}E & EX_a\\
X'_aE' &E'X_a
\ar@{=>}^{\ka_a}"1,1";"1,2"
\ar@{=>}_{\ka'_a}"2,1";"2,2"
\ar@{=>}_{X'_a\et}"1,1";"2,1"
\ar@{=>}^{\et X_a}"1,2";"2,2"
}
$$
commute for all $a\in G$.
\end{dfn}

\begin{dfn}
\label{dfn:G-covering}
Let $\calB$ be a $\k$-category.
Then by $\De(\calB)$ we denote the $G$-category with the trivial
$G$-action.

(1)
Let $\calC = (\calC, X)$ be a $G$-category, and $\calB$ a $\k$-category.
Then an $G$-\emph{invariant} functor $\calC \to \calB$ is a $G$-equivariant functor $\calC \to \De(\calB)$.

(2)
A $G$-invariant functor $(E, \ka) \colon \calC \to \calB$ is called
a $G$-\emph{precovering} if both of the following are isomorphisms:
$$
\begin{aligned}
\bi{(E, \ka)}{y}{x}^{(1)} &\colon \Ds_{a \in G}\bi{\calC}{y}{ax}
\to \bi{\calB}{E(y)}{E(x)}, (f_a)_{a \in G} \mapsto \sum_{a \in G} E(f_a) \ka_{a,x},\\
\bi{(E, \ka)}{y}{x}^{(2)} &\colon \Ds_{b \in G}\bi{\calC}{by}{x}
\to \bi{\calB}{E(y)}{E(x)}, (f_b)_{b \in G} \mapsto \sum_{b \in G}  \ka_{b,x}\inv E(f_b).
\end{aligned}
$$
Note that both are isomorphisms if and only if so is one of them
by \cite[Proposition 1.6]{Asa11}.

(3)
A $G$-precovering is called a $G$-\emph{covering} if it is dense, namely
for each $y \in \calB_0$, there exists some $x \in \calC_0$
such that $E(x) \iso y$ in $\calB$.
\end{dfn}

\begin{dfn}\label{dfn:G-Cat}
A 2-category $G$-$\mathbf{Cat}$ of small $G$-categories is defined as follows.
\begin{itemize}
\item The objects are the small $G$-categories.
\item The 1-morphisms are the $G$-equivariant functors between objects.
\item The identity 1-morphism of an object $\calC$ is the 1-morphism
$(\id_\calC, (\id_{\id_\calC})_{a \in G})$.
\item The 2-morphisms are the morphisms of $G$-equivariant functors.
\item The identity 2-morphism of a 1-morphism $(E, \ka) \colon \calC \to \calC'$
is the identity natural transformation  $\id_E$ of $E$, which is clearly a 2-morphism.
\item The composite $(E',\ka')(E,\ka)$ of 1-morphisms
$\mathcal{C}\ya{(E,\ka)}  \mathcal{C}'\ya{(E',\ka')}  \calC''
$
is defined by
$$
(E',\ka')(E,\ka):= (E'E,((E'\ka_{a})(\ka_{a}'E))_{a\in G})\colon\mathcal{C}\to\mathcal{C}''.
$$
\item The vertical and the horizontal compositions of 2-morphisms are given
by the usual ones of natural transformations.
\end{itemize}
\end{dfn}

\subsection{\texorpdfstring{$G$}{G}-graded categories}
Next we recall definitions of $G$-graded categories and
their 2-category $\GGrCat$.

\begin{dfn}\label{dfn:graded-cat}
(1) A $G$-\emph{graded} $\k$-category is a category $\mathcal{B}$
together with a family of direct sum decompositions
$\bi{\calB} yx =\bigoplus_{a\in G} \bi{\calB} yx^{a}$
$(x,y \in \calB_0)$ of $\Bbbk$-modules such that
$\bi{\calB} zy^{b}\cdot \bi{\calB} yx^{a} \subseteq \bi{\calB} zx^{ba}$
for all $x,y\in\mathcal{B}$ and $a,b\in G$. \
It is easy to see that $\id_x \in \calB^1(x, x)$ for all $x \in \calB_0$.

$(2)$ A \emph{degree-preserving} functor is a pair $(H,r)$ of a
$\k$-functor $H\colon\mathcal{B}\to\mathcal{A}$ of $G$-graded categories
and a map $r\colon \mathcal{B}_0\to G$ such that
\[
H(\bi{\calB} yx^{r_{y}a})\subseteq \bi{\calA}{Hy}{Hx}^{ar_{x}}
\]
(or equivalently $H(\bi{\calB} yx^{a})\subseteq \bi{\calA}{Hy}{Hx}^{r_{y}\inv ar_{x}}$)
for all $x,y\in\mathcal{B}$ and $a\in G$.
This $r$ is called a \emph{degree adjuster} of $H$.

$(3)$ A $\k$-functor $H\colon\mathcal{B}\to\mathcal{A}$ of $G$-graded
categories is called a \emph{strictly} degree-preserving functor if
$(H,1)$ is a degree-preserving functor, where 1 denotes the constant
map $\mathcal{B}_0\to G$ with value $1\in G$, i.e., if
$H(\bi{\calB} yx^{a})\subseteq \bi{\calA}{Hy}{Hx}^{a}$
for all $x,y\in\mathcal{B}$ and $a\in G$.

$(4)$ Let $(H,r),(I,s)\colon\mathcal{B}\to\mathcal{A}$ be degree-preserving
functors. Then a natural transformation $\theta\colon H\To I$ is
called a \emph{morphism} of degree-preserving functors if
$\theta x\in \bi{\calA}{Ix}{Hx}^{s_{x}^{-1}r_{x}}$
for all $x\in\mathcal{B}$.
\end{dfn}

\begin{dfn}\label{dfn:G-GrCat}
A 2-category $G$-$\mathbf{GrCat}$ of small $G$-graded categories is defined as follows.
\begin{itemize}
\item The objects are the small $G$-graded categories.
\item The 1-morphisms are the degree-preserving functors between objects.
\item
The identity 1-morphism of an object $\calB$ is the 1-morphism $(\id_\calB, 1)$.
\item The 2-morphisms are the morphisms of degree-preserving functors.
\item
The identity 2-morphism of a 1-morphism $(H, r)\colon \calB \to \calA$
is the identity natural transformation $\id_H$ of $H$, which is a 2-morphism
(because $(\id_{H})x =\id_{Hx} \in \bi{\calA}{Hx}{Hx}^1 = \bi{\calA}{Hx}{Hx}^{r_x\inv r_x}$
 for all $x \in \calB$).
\item The composite $(H',r')(H,r)$ of 1-morphisms
$\mathcal{B}\ya{(H,r)} \mathcal{B}'\ya{(H',r')}  \mathcal{B}''$
 is defined by
$$
(H',r')(H,r):= (H'H,(r_{x}r'_{Hx})_{x\in\mathcal{B}})\colon\mathcal{B}\to\mathcal{B}''.
$$
\item The vertical and the horizontal compositions of 2-morphisms are given
by the usual ones of natural transformations.
\end{itemize}
\end{dfn}

\subsection{Orbit categories and smash products}
Finally we recall definitions of orbit categories and smash products, and their
relationships.

\begin{dfn}
\label{dfn:orbit-category}
Let $\calC$ be a $G$-category.
Then the {\em orbit category} $\mathcal{C}/G$ of
$\mathcal{C}$ by $G$ is a category defined as follows.
\begin{itemize}
\item $(\C/G)_0:=\C_0$;
\item For any $x,y\in G$, $\bi{(\C/G)} yx:=\Ds_{a\in G} \bi{\calC}{y}{ax}$;
and
\item For any $x\ya{f}y\ya{g}z$ in $\C/G$, we set
\[
gf:=\left(\sum_{\smat{a,b\in G\\ba=c}} g_{b}\cdot b(f_{a})\right)_{c\in G}.
\]
\item
For each $x\in(\mathcal{C}/G)_0$ its identity $\id_{x}:= \id_x^{\calC/G}$ in $\calC/G$
is given by $\id_{x}=(\delta_{a,1}\id_{x}^{\calC})_{a \in G}$,
where $\id_{x}^{\calC}$ is the identity of $x$ in $\calC$.
\end{itemize}
By setting $\bi{(\calC/G)}{y}{x}^a:= \bi{\calC}{y}{ax}$ for all $x, y \in \calC_0$
and $a \in G$, the decompositions
\[
\bi{(\calC/G)} yx = \Ds_{a\in G} \bi{(\calC/G)} yx^a
\]
makes $\calC/G$ a $G$-graded category.
\end{dfn}

\begin{dfn}
Let $\mathcal{B}$ be a $G$-graded category.
Then the \emph{smash product} $\mathcal{B}\#G$ is a category defined as follows.
\begin{itemize}
\item $(\mathcal{B}\#G)_0:=\mathcal{B}_0\times G$, we set $x^{(a)}:=(x,a)$
for all $x\in\mathcal{B}$ and $a\in G$.
\item $\bi{(\calB\# G)}{y^{(b)}}{x^{(a)}}
:= \{\Smor(b,f,a):= (b,f,a) \mid f \in \bi{\calB} yx^{b^{-1}a}\}
= \{b\} \times \bi{\calB} yx^{b^{-1}a} \times \{a\}$
for all $x^{(a)},y^{(b)} \in \mathcal{B}\#G$.
This definition make the union
\[
\bigcup_{x^{(a)},\, y^{(b)} \in (\calB \# G)_0}
\bi{(\calB\# G)}{y^{(b)}}{x^{(a)}}
\]
disjoint.
This is sometimes identified with
$\bi{\calB} yx^{b^{-1}a}$ by the correspondence
$f \leftrightarrow \Smor(b,f,a)$ if there seems to be
no confusion.
\item For any $x^{(a)},y^{(b)},z^{(c)}\in\mathcal{B}\#G$ the composition
is given by the following commutative diagram
\[
\xymatrix{
\bi{(\calB\# G)}{z^{(c)}}{y^{(b)}}\times\bi{(\calB\# G)}{y^{(b)}}{x^{(a)}} & \bi{(\calB\#G)}{z^{(c)}}{x^{(a)}}\\
\bi{\calB} zy^{c^{-1}b}\times \bi{\calB} yx^{b^{-1}a} & \bi{\calB} zx^{c^{-1}a},
\ar"1,1";"1,2"
\ar"2,1";"2,2"
\ar@{=}"1,1";"2,1"
\ar@{=}"1,2";"2,2"
}
\]
where the lower horizontal homomorphism is given by the composition of
$\mathcal{B}$.
\item
For each $x^{(a)} \in (\calB \# G)_0$ its identity $\id_{x^{(a)}}$ in $\calB \# G$
is given by $\id_x \in \bi{\calB} xx^1$.
\end{itemize}
$\mathcal{B}\#G$ has a free $G$-action defined as follows:
For each $c\in G$ and $x^{(a)}\in\mathcal{B}\#G$, ${c}x^{(a)}:=x^{(ca)}$; and
for each $\Smor(b,f,a) \in \bi{(\calB\# G)}{y^{(b)}}{x^{(a)}}$  with $f \in \bi{\calB} yx^{b^{-1}a}$,
noting that 
$\bi{(\calB\# G)}{y^{(b)}}{x^{(a)}}=\bi{(\calB\# G)}{y^{(cb)}}{x^{(ca)}}$, we set
${c}(\Smor(b,f,a)):=\Smor(cb,f,ca)$.
By the shorter notation, it becomes $cf = f$.
\end{dfn}

The following two propositions were proved in \cite{Asa11}.
\begin{prp}[{\cite[Proposition 5.6]{Asa11}}]
\label{smash-orb}
Let $\calB$ be a $G$-graded category.
Then there is a strictly degree-preserving equivalence
$\om_\calB\colon \calB \to (\calB \# G)/G$ of
$G$-graded categories.
\end{prp}

\begin{prp}[{\cite[Theorem 5.10]{Asa11}}]
\label{liberalization}
Let $\calC$ be a category with a $G$-action.
Then there is a $G$-equivariant equivalence
$\ze_\calC\colon \calC \to (\calC/G)\# G$.
\end{prp}

Note that we changed the notation $\ep_\calC$ used
in \cite{Asa11} to $\ze_\calC$ as used in
\cite{Asa-book}.

In fact, the orbit category construction and the smash product construction can be
extended to 2-functors $?/G \colon \GCat \to \GGrCat$
and $?\# G \colon \GGrCat \to \GCat$, respectively, and
they are inverses to each other as stated in the following theorem,
where $\om:= (\om_\calB)_\calB$ and
$\ze:= (\ze_\calC)_\calC$ are 2-natural isomorhisms.

\begin{thm}[{\cite[Theorem 7.5]{Asa12}}]
\label{thm-2-eq}
$?/G$ is strictly left $2$-adjoint to $?\# G$ and
they are mutual $2$-quasi-inverses.
\end{thm}

\begin{rmk}
\label{rmk:omega-zeta-identify}
$\om_\calB\colon \calB \to (\calB \# G)/G$ above is an equivalence in the 2-category
$\GGrCat$ and
$\ze_\calC\colon \calC \to (\calC/G)\# G$ above  is an equivalence in the 2-category
$\GCat$.
By these equivalences we identify  $(\calB \# G)/G$ with $\calB$, and
$(\calC/G)\# G$ with $\calC$ in the following sections.
Here we note that the quasi-inverse $\om'_\calB$ of $\om_\calB$ given in \cite{Asa12}
is not strictly degree preserving, which forced us to define degree adjusters.
\end{rmk}

\section{\texorpdfstring{$G$}{G}-invariant bimodules and \texorpdfstring{$G$}{G}-graded bimodues}

\subsection{\texorpdfstring{$G$}{G}-invariant bimodules}

\begin{dfn}
\label{dfn:G-inv}
Let $R = (R, X)$ and $S = (S, Y)$ be small $G$-categories.

(1) A {\em $G$-invariant} $S$-$R$-bimodule is a pair $(M, \ph)$ of
an $S$-$R$-bimodule $M$ and a family $\ph:= (\ph_a)_{a \in G}$ of
natural isomorphisms
$\ph_a\colon M \to \bi{M}{Y_a}{X_a}$,
where $\ph_a = (\bi{(\ph_a)} yx)_{(y,x)\in S_0 \times R_0}$, and
$\bi{(\ph_a)} yx \colon \bi{M} yx \to \bi{M}{ay}{ax}$ is in $\kMod$,
such that the following diagram commutes for all $a, b \in G$ and all $(y, x) \in S_0 \times R_0$:
$$
\xymatrix{
\bi{M} yx & \bi{M}{ay}{ax}\\
           & \bi{M}{bay}{bax}.
\ar^{\bi{(\ph_a)} yx}"1,1";"1,2"
\ar^{\bi{(\ph_b)}{ay}{ax}}"1,2";"2,2"
\ar_{\bi{(\ph_{ba})} yx}"1,1";"2,2"
}
$$

(2) Let $(M, \ph)$ and $(N, \ps)$ be $G$-invariant $S$-$R$-bimodules.
A {\em morphism}
$$
(M, \ph) \to (N, \ps)
$$ 
is an $S$-$R$-bimodule morphism $f\colon M \to N$
such that the following diagram commutes for all $a\in G$ and
all $(y, x) \in S_0 \times R_0$:
$$
\xymatrix{
\bi{M} yx& \bi{M}{ay}{ax}\\
\bi{N} yx & \bi{N}{ay}{ax}.
\ar^{\bi{(\ph_a)} yx}"1,1";"1,2"
\ar_{\bi{(\ps_a)} yx}"2,1";"2,2"
\ar_{\bi{f} yx}"1,1";"2,1"
\ar^{\bi{f}{ay}{ax}}"1,2";"2,2"
}
$$

(3)  Let $(L, \th), (M, \ph)$ and $(N, \ps)$ be $G$-invariant $S$-$R$-bimodules, and
$$
f \colon (L, \th)\to (M, \ph),\, g \colon (M, \ph) \to (N, \ps)
$$
morphisms of $G$-invariant $S$-$R$-bimodules.
Then as is easily seen, the composite 
$gf \colon L \to N$ in $\Bi(S, \Mod, R)$
turns out to be a morphism in $\Ginv(\Bi(S, \Mod, R))$, which
is defined to be the composite 
$gf \colon (L, \th) \to (N, \ps)$
in $\Ginv(\Bi(S, \Mod, R))$.

(4) Let $(M, \ph)$ be a $G$-invariant $S$-$R$-bimodule.
Then $\id_M$ of the $S$-$R$-bimodule $M$ turns out to be
the identity with respect to the composition defined above.

(5) The class of all $G$-invariant $S$-$R$-bimodules together with all morphisms between them
forms a $\k$-category, which we denote by $\Ginv(\biMod{S}{R})$.
\end{dfn}

\begin{rmk}
The commutativity of the diagram in (1) above for $a=b=1$ shows that
$\ph_1 = \id_M$ because $\ph_1^2 = \ph_1$ and $\ph_1$ is a natural
isomorphism.
This also shows that $(\bi{(\ph_a)} yx)\inv = \bi{(\ph_{a\inv})}{ay}{ax}$
for all $a\in G$ and all $(y, x)\in S_0 \times R_0$.
\end{rmk}

\begin{exm}
(1) $(R, (X_a)_{a\in G})$ is a $G$-invariant $R$-$R$-bimodule.

(2) Let $(y, x) \in S_0 \times R_0$.
Then
$\Ds_{a\in G}S_{ay}\ox_\k {}_{ax}R$ has the
{\em canonical} $G$-invariant structure $\ph = (\ph_b)_{b\in G}$ defined by the composite
\begin{equation}\label{eqn:can-G-inv}
{}_v(\ph_b)_u\colon \Ds_{a\in G}{}_v S_{ay}\ox_\k {}_{ax}R_u \ya{Y_b\ox_\k X_b} \Ds_{a\in G}{}_{bv} S_{bay}\ox_\k {}_{bax}R_{bu}
\isoto \Ds_{a\in G}{}_{bv} S_{ay}\ox_\k {}_{ax}R_{bu}
\end{equation}
for all $b \in G$ and $(v,u) \in B_0 \times A_0$,
thus, ${}_v(\ph_b)_u(s_a \ox r_a)_{a\in G}:= (bs_{b\inv a} \ox br_{b\inv a})_{a\in G}$
for all $s_a \in {}_vS_{ay}, r_a \in {}_{ax}R_u$.
\end{exm}

\subsection{\texorpdfstring{$G$}{G}-graded bimodules}

\begin{dfn}
Let $A$ and $B$ be $G$-graded small $\k$-categories.

(1) A $G$-graded $B$-$A$-bimodule is a $B$-$A$-bimodule $M$ together with
decompositions $\bi{M} yx = \Ds_{a\in G}\bi{M} yx^a$ in $\kMod$ for all $(y,x) \in B_0 \times A_0$
such that
$$
\bi{B}{y'}{y}^c \cdot \bi{M} yx^a \cdot \bi{A}{x}{x'}^b \subseteq \bi{M}{y'}{x'}^{cab}
$$
for all $a, b, c \in G$ and all $x, x' \in A_0, y, y' \in B_0$.

(2) Let $M$ and $N$ be $G$-graded $B$-$A$-bimodules.
Then a {\em morphism} $M \to N$ is a $B$-$A$-bimodule morphism $f \colon M \to N$
such that $f(\bi{M} yx^a) \subseteq \bi{N}{y}{x}^a$ for all $a \in G$ and all $(y,x) \in B_0\times A_0$.
Hence a morphism $f$ induces morphisms $\Bi(y,f,x)^a \colon \Bi(y,M,x)^a \to \Bi(y,N,x)^a$, and we may write
$\Bi(x,f,y) = \Ds_{a\in G}\Bi(y,f,x)^a \colon \Bi(y,M,x) \to \Bi(y,N,x)$ for all $(y,x) \in B_0\times A_0$.

(3) Let $L$, $M$ and $N$ be $G$-graded $B$-$A$-bimodules, and
$f \colon L \to M$, $g \colon M \to N$ morphisms of $G$-graded
$B$-$A$-bimodules.
Then as is easily seen,
the composite $gf \colon L \to N$ in $\Bi(B,\Mod, A)$ turns out
to be a morphism of $G$-graded $B$-$A$-bimodules, which is defined
to be the composite $gf$ of $f$ and $g$ as a morphism of $G$-graded $B$-$A$-bimodules.

(4) Let $M$ be a $G$-graded $B$-$A$-bimodule.
Then the identity $\id_M$ of $M$ as a morphism of $B$-$A$-bimodules
turns out to be the identity of $M$ as a morphism of
$G$-graded $B$-$A$-bimodules.

(5) The class of all $G$-graded $B$-$A$-bimodules together with all morphisms between them
forms a $\k$-category, which we denote by $\Ggr(\biMod{B}{A})$.
\end{dfn}

\begin{rmk}
Recall that the category $A\op$ turns out to be a $G\op$-graded category by setting
$$
\bi{(A\op)}{x'}{x}^a:= \bi{A}{x}{x'}^{a}
$$
to be the degree $a$ part for all $x, x' \in A_0$, $a \in G$.
Indeed, by denoting the composition of $A\op$ and of $G\op$ by $\circ$,
we have
$$
\begin{aligned}
&\, \bi{(A\op)}{x''}{x'}^{a'} \circ \bi{(A\op)}{x'}{x}^a
= \bi{A}{x'}{x''}^{a'} \circ \bi{A}{x}{x'}^{a}
= \bi{A}{x}{x'}^{a} \cdot \bi{A}{x'}{x''}^{a'}\\
\subseteq\, & \bi{A}{x}{x''}^{aa'}
= \bi{A}{x}{x''}^{a' \circ a}
= \bi{(A\op)}{x''}{x}^{a' \circ a}
\end{aligned}
$$
for all $x, x', x'' \in A_0$, $a \in G$.
\end{rmk}

\begin{exm}
Let $(y, x) \in B_0 \times A_0$.
Then $B_{y}\ox_\k {}_{x}A$ has a $G$-grading
defined by
\begin{equation}
\label{eq:can-G-gr}
B_{y}\ox_\k {}_{x}A
= \Ds_{c\in G} \left(\Ds_{ba = c}B_{y}^{b}\ox_\k {}_{x}A^{a}\right).
\end{equation}
Indeed, set $M^c:= \Ds_{ba = c}B_{y}^{b}\ox_\k {}_{x}A^{a}$
for all $c \in G$.
Then we have
$$
\begin{aligned}
&\, \bi{B}{y''}{y'}^{b'}\, \bi{M}{y'}{x'}^c\, \bi{A}{x'}{x''}^{a'}
= \textstyle\Ds_{ba = c}(\bi{B}{y''}{y'}^{b'}\, \bi{B}{y'}{y}^{b})\ox_\k 
(\bi{A}{x}{x'}^{a}\, \bi{A}{x'}{x''}^{a'})\\
\subseteq\, & 
\textstyle\Ds_{ba = c}\bi{B}{y''}{y}^{b'b}\ox_\k \bi{A}{x}{x''}^{aa'}
\subseteq
\textstyle\Ds_{b''a'' = b'ca'}\bi{B}{y''}{y}^{b''}\ox_\k \bi{A}{x}{x''}^{a''}
= \bi{M}{y''}{x''}^{b' c a'}
\end{aligned}
$$
for all $x',x'' \in A_0$, $y', y'' \in B_0$, $a',\, b',\, c \in G$.
This $G$-grading is called the {\em canonical}
$G$-grading of $B_{y}\ox_\k {}_{x}A$.
\end{exm}

\begin{rmk}
\label{rmk:G-gr-bimod-left}

If $G$ is an abelian group, then
we can introduce a $G$-grading on the category
$B \ox_k A\op$ using the formula
\eqref{eq:can-G-gr} as follows (called the total grading):
\begin{equation}
\label{eq:total-G-g}
\bi{(B\ox_\k A\op)}{(y',x')}{(y,x)} = \bi{B}{y'}{y} \ox_\k \bi{A}{x}{x'}
:= \Ds_{c \in G}\left(\Ds_{ba=c} \bi{B}{y'}{y}^b \ox_k \bi{A}{x}{x'}^a\right)
\end{equation}
for all $(y,x),\, (y',x') \in B_0 \times A_0$,
namely the degree $c$ part is given by
$$
\bi{(B\ox_\k A\op)}{(y',x')}{(y,x)}^c:=
\Ds_{ba=c} \bi{B}{y'}{y}^b \ox_k \bi{A}{x}{x'}^a
$$
for all $c \in G$.
Indeed, for any $c, c' \in G$ and
any $(y,x),\, (y',x'),\, (y'',x'') \in B_0 \times A_0$, we have
$$
\begin{aligned}
&\hspace{1.3em}\bi{(B\ox_\k A\op)}{(y'',x'')}{(y',x')}^{c'}\cdot \bi{(B\ox_\k A\op)}{(y',x')}{(y,x)}^c\\
&= \textstyle (\Ds_{b'a'=c'} \bi{B}{y''}{y'}^{b'} \ox_k \bi{A}{x'}{x''}^{a'})
(\Ds_{ba=c} \bi{B}{y'}{y}^b \ox_k \bi{A}{x}{x'}^a)\\
&\subseteq
\textstyle \Ds_{b'a'=c'}\Ds_{ba=c} \bi{B}{y''}{y}^{b'b}\ox_k \bi{A}{x}{x''}^{aa'}
\subseteq
\textstyle \Ds_{b''a''=c'c} \bi{B}{y''}{y}^{b''}\ox_k \bi{A}{x}{x''}^{a''}
\\
&= \bi{(B\ox_\k A\op)}{(y'',x'')}{(y,x)}^{c'c},
\end{aligned}
$$
where the last inclusion holds
because $b'baa' = b'a'ba$ for all $a,a',b,b' \in G$.

In this case, with this $G$-grading of $B \ox_k A\op$,
a direct sum decomposition $M = \Ds_{a \in G} M^a$ in $\ltMod{\k}$
of a $B$-$A$-bimodule $M$ gives a structure of a $G$-graded
left $B \ox_k A\op$-module if and only if it gives $M$
a structure of a $G$-graded $B$-$A$-bimodule because
for any $x,x' \in A_0$, $y,y' \in B_0$, and $c,c' \in G$, we have
$$
\begin{aligned}
\textstyle(\Ds_{ba=c} \bi{B}{y'}{y}^b \ox_k \bi{A}{x}{x'}^a) \bi{M}{y}{x}^{c'}
&= \textstyle\Ds_{ba=c} \bi{B}{y'}{y}^b\, \bi{M}{y}{x}^{c'}\, \bi{A}{x}{x'}^a, \text{ and}\\
\bi{M}{y'}{x'}^{cc'} &= \bi{M}{y'}{x'}^{bc'a}
\quad\text{for all $a, b \in G$ with $c = ba$}.
\end{aligned}
$$
Hence we can identify the category $\Ggr(\biMod{B}{A})$
of $G$-graded $B$-$A$-bimodules
with the category $\Ggr({}_{B \ox_\k A\op}\Mod)$ of
$G$-graded left $B \ox_\k A\op$-modules.

Note that if $G$ is not abelian, then the argument above does not work.
Nevertheless, the decomposition \eqref{eq:can-G-gr} gives us
a structure of a $G$-graded $B$-$A$-bimodule to the left $B\ox_\k A\op$-module
$(B\ox_\k A\op)_{(y,x)} = B_{y}\ox_\k {}_{x}A$.
\end{rmk}

\begin{rmk}
Let $M$ be a $G$-graded $B$-$A$-bimodule, $x \in A_0$ and $y \in B_0$.
Then
\begin{enumerate}
\item 
$M_x$ turns out to be a $G$-graded
left $S$-module by the decomposition
$M_x = \Ds_{a\in G} M_x^a$,
where for each $y \in B_0$, ${}_{y}(M_x^a):= {}_{y}M^a_x$, and
for each $f \in {}_{y'}B^b_y$, ${}_f(M^a_x) \colon {}_yM^a_x \to {}_{y'}M^{ba}_x$
is defined by $m \mapsto fm$
for all $m \in {}_yM^{a}_x$.
\item
Similarly, ${}_yM$ turns out to be a $G$-graded
right $R$-module by the decomposition
${}_yM = \Ds_{a\in G} {}_yM^a$.
\end{enumerate}
\end{rmk}

\begin{prp}
\label{prp:graded-prj}
Let $\Fgt \colon \Ggr(\ltMod{A}) \to \ltMod A$ be the forgetful functor,
and $P \in \Ggr(\ltMod A)_0$.
Then $P$ is projective in $\Ggr(\ltMod A)$ if and only if
$\Fgt(P)$ is projective in $\ltMod A$.
\end{prp}

\begin{proof}
This can be proved as in the case of $G$-graded algebras.
\end{proof}

The following is immediate by Remark \ref{rmk:G-gr-bimod-left} and
Proposition \ref{prp:graded-prj}
in the case where $G$ is abelian, but we mention that
it is even true for the non-abelian case, which will be slightly generalized and proved in  Proposition \ref{prp:graded-prj-bimod-eq}.

\begin{prp}
\label{prp:graded-prj-bimod}
Let $\Fgt \colon \Ggr(\biMod{B}{A}) \to \biMod{B}{A}$ be
the forgetful functor, and
$P \in \Ggr(\biMod{B}{A})_0$.
Then $P$ is projective in $\Ggr(\biMod{B}{A})$ if and only if
$\Fgt(P)$ is projective in $\biMod{B}{A}$.
\qed
\end{prp}

\subsection{Finitely generated projective \texorpdfstring{$G$}{G}-graded modules}

\begin{dfn}
Let $a \in G$.
Then we can define an automorphism $\la_a$ (resp.\ $\ro_a$) of $\Ggr(\ltMod{A})$
called the \emph{left shift} (resp.\ \emph{right shift}) by $a$ as follows:
Let $f \colon M \to N$ be in $\Ggr(\ltMod{A})$
with $M = \Ds_{b \in G} M^b$.
Then we set
$\la_a(M):= M$ (resp.\ $\ro_a(M):= M$) as a left $A$-module, and
$\la_a(M) = \Ds_{b \in G} (\la_a(M))^b$
(resp.\  $\ro_a(M) = \Ds_{b \in G} (\ro_a(M))^b$), where
$$
(\la_a(M))^b:= M^{a\inv b}\ (\text{resp.\ } (\ro_a(M))^b:= M^{b a\inv}).
$$
Moreover, we set $\la_a(f):= f$ (resp.\ $\ro_a(f):= f$).

If $G$ is abelian, then $\la_a$ and $\ro_a$ coincide, which will be denoted by
$\si_a$.
\end{dfn}

We recall the following definitions from \cite{Gab}.

\begin{dfn}
\label{dfn:loc-fd-cat}
Let $\calC$ be a skeletally small linear category over a field.
Then $\calC$ is called a \emph{locally finite-dimensional category}
if (i) $\calC$ is \emph{basic}, 
(i.e., distinct objects are not isomorphic in $\calC$),
(ii) $\bi{\calC}{x}{x}$ is a local algebra
for all $x \in \calC_0$,
and (iii) $\calC$ is \emph{Hom-finite}
(i.e., for any $x, y \in \calC_0$, $\bi{\calC}{y}{x}$ is finite-dimensional).
A locally finite-dimensional category $\calC$ is called
a \emph{locally bounded category}
if for each $x \in \calC_0$ the set
$\{y \in \calC_0 \mid \bi{\calC}{y}{x} \ne 0 \text{ or }
\bi{\calC}{x}{y} \ne 0\}$ is a finite set,
or equivalently, both $\bi{\calC}{x}{}$ and $\bi{\calC}{}{x}$
are finite-dimensional.
Note that a locally finite-dimensional category is necessarily a small category.
\end{dfn}

\begin{rmk}
\label{rmk:fg-proj-ind}
If $\calC$ is a locally finite-dimensional category, then
each finitely generated indecomposable projective right $\calC$-module $P$ has the form
${}_x\calC$ for some $x \in \calC_0$.
Indeed, since it is finitely generated there exists an epimorphism
$f \colon \Ds_{i=1}^n {}_{x_i}\calC \to P$ for some finitely many $x_i \in \calC_0$.
Since $P$ is projective, $f$ is a retraction.
Since each ${}_{x_i}\calC$ has a local endomorphism algebra,
$P \iso {}_{x_i}\calC$ for some $i$ by the Krull--Schmidt lemma.
\end{rmk}

We will give explicit forms of finitely generated projective objects
in the categories $\Ggr(\Mod_A)$ and $\Ggr(\ltMod{A})$
in the case where $A$ is a locally finite-dimensional
category.
For this sake, we need the following.

\begin{lem}
\label{lem:smashprod-almost-loc-findim}
Assume that $\k$ is a field, and that $A$ is a $G$-graded locally finite-dimensional category.
Then $A\# G$ is Hom-finite and the endomorphism algebra of each object is local.
Hence its skeleton is a locally finite-dimensional category.
In particular, each finitely generated indecomposable projective right $A\# G$-module has the form 
${}_{x^{(a)}}(A\# G)$ for some $x \in A_0$, $a \in G$
by Remark \ref{rmk:fg-proj-ind}.

In the following cases, $A \# G$ itself is basic,
and hence becomes a locally finite-dimensional category:
\begin{enumerate}
\item
$G$ is a torsion-free group.
\item
$A$ is a locally bounded category presented
as a path-category $\k(Q, I, W)$ of a bound quiver $(Q, I)$ whose $G$-grading is
given by a homogeneous $G$-weight on $(Q, I)$.
\end{enumerate}
\end{lem}

\begin{proof}
For any objects $x^{(a)}$ and $y^{(b)}$ in $A\# G$, we have an isomorphism
$$
\bi{(A\# G)}{y^{(b)}}{x^{(a)}} \iso \bi{A}{y}{x}^{b\inv a}
$$
of $\k$-vector spaces, the right hand side of which is finite-dimensional
as a subspace of $\bi{A}{y}{x}$.  Thus $A\# G$ is Hom-finite.

Consider the case where $x^{(a)} = y^{(b)}$.
Then we have an isomorphism
$$
\bi{(A\# G)}{x^{(a)}}{x^{(a)}} \iso \bi{A}{x}{x}^{1}
$$
of algebras, the right hand side of which turns out to be a local algebra
because so is $\bi{A}{x}{x}$.

We show that $A\# G$ becomes basic in the case (1).
Let $x^{(a)}$ and $y^{(b)}$ be objects in $A\# G$, and assume that
they are isomorphic, say there exist morphisms $f' \in \bi{(A\# G)}{y^{(b)}}{x^{(a)}}$
and $g' \in \bi{(A\# G)}{x^{(a)}}{y^{(b)}}$ such that
$g'f' = \id_{x^{(a)}}$ and $f'g' = \id_{y^{(b)}}$.
Then there exist $f \in \bi{A}{y}{x}^{b\inv a}$ and $g \in \bi{A}{x}{y}^{a\inv b}$
such that $f' = \smor{b}{f}{a}$ and $g' = \smor{a}{g}{b}$,
and that $gf = \id_x$ and $fg = \id_y$ in $A$.
Therefore $x$ and $y$ are isomorphic in the basic category $A$.
Hence $x = y$.

For each $c \in G$ with $c \ne 1$, we show that $\bi{A}{x}{x}^c$ is nilpotent.
Since $\bi{A}{x}{x} = \Ds_{d\in G} \bi{A}{x}{x}^d$ is finite-dimensional,
the set $H:= \{d \in G \mid \bi{A}{x}{x}^d \ne 0\}$ is finite.
But by (1), the set $\{c^n \mid 0 < n \in \bbZ\}$ is infinite, and hence
there exists some integer $n > 0$ such that $c^n \not\in H$.
Thus $\bi{A}{x}{x}^{c^n} = 0$.
Therefore $(\bi{A}{x}{x}^c)^n \subseteq \bi{A}{x}{x}^{c^n} = 0$, as desired.

Since $f \in \bi{A}{x}{x}^{b\inv a}$, if $a \ne b$, then $b\inv a \ne 1$
shows that $f$ would be nilpotent, a contradiction.  Hence $a = b$.
As a consequence, we have $x^{(a)} = y^{(b)}$.
Therefore, $A\# G$ is basic.

In the case (2), $A\# G$ is isomorphic to the path-category of a bound quiver
$(Q_{G,W}, I_{G,W})$ by \cite[Theorem 1.18]{Asa18}
(see also \cite[Theorem 6.2.18]{Asa-book}), and hence $A\# G$ is locally bounded.
\end{proof}

\begin{prp}
\label{prp:graded_proj_over_loc_fd}
Assume that $\k$ is a field,
$A$ is a $G$-graded locally finite-dimensional category.
Then each finitely generated projective object
in $\Ggr(\Mod_A)\linebreak[2] \ (\text{resp.\ }\Ggr(\ltMod A))$
has the following form:
$$
\Ds_{i=1}^n \la_{a_i}({}_{x_i}A)
\quad \left(\text{resp.\ } \Ds_{i=1}^n \ro_{a_i}(A_{x_i})\right)
$$
for some $a_1, \dots, a_n \in G$ and $x_1 \dots x_n \in A_0$,
$n \ge 0$.
\end{prp}

\begin{proof}
We have a $G$-invariant functor $(Q, \id) \colon A\# G \to A$
define by $Q(x^{(a)}):= x$ and $Q(\Smor(b,f,a)):= f$
for all $x^{(a)} \in (A\# G)_0$ and
$\Smor(b,f,a) \in \bi{(A\# G)}{y^{(b)}}{x^{(a)}}$,
which satisfies $Q = Q X'(a)$ for all $a \in G$,
where $X'$ is the $G$-action of $A\# G$.
Then it turns out to be a $G$-covering (see Definition \ref{dfn:G-covering}).
Indeed, it is surjective on objects, and
$$
(Q,\id)^{(1)}_{y,x} \colon \Ds_{c\in G} \bi{(A\# G)}{y^{(b)}}{x^{(ca)}}
\isoto \Ds_{c\in G} \bi{A}{y}{x}^{b\inv ca} = \bi{A}{y}{x}.
$$
is an isomorphism for all $x^{(a)}, y^{(b)} \in (A\# G)_0$.
By the same arguments as in \cite[Proposition 7.4.8]{Asa-book},
the push-down $Q\down \colon \Mod_{A\# G} \to \Mod_{A}$
of $Q$ has the following explicit form:
\begin{equation}
\label{eq:explicit-Qdown}
(Q\down M)(x) = \Ds_{a \in G}M(x^{(a)})
\end{equation}
for all $M \in \Mod_{A\# G}$ and $x \in A_0$.
This is done by regarding the natural isomorphism
$$
(Q\down M)(x) =  M\ox_{A\# G} \bi{A}{Q(\blank)}{x}
\iso
M \ox_{A\# G} (\textstyle\Ds_{a\in G}(A\# G)_{ax^{(1)}})\\
\iso
\textstyle\Ds_{c\in G} M(x^{(a)})
$$
as an identity.
Then $Q\down$ factors through the forgetful functor,
namely we have the following strictly commutative diagram:
$$
\begin{tikzcd}
\Mod_{A\# G} & \Mod_{A}\\
\Ggr(\Mod_{A}) & \Mod_{A}
\Ar{1-1}{1-2}{"Q\down"}
\Ar{2-1}{2-2}{"\Fgt" '}
\Ar{1-1}{2-1}{"Q'\down" '}
\Ar{1-2}{2-2}{equal}
\end{tikzcd},
$$
where $Q'\down$ is an equivalence induced by $Q\down$ as 
in \cite[Proposition 7.6.2]{Asa-book} (see also \cite[Theorem 6.6]{Asa11}).
By Lemma \ref{lem:smashprod-almost-loc-findim},
each indecomposable projective right $A\# G$-module has the form 
${}_{x^{(a)}}(A\# G)$ for some $x \in A_0$, $a \in G$.
Hence each indecomposable projective object $P$ in $\Ggr(\Mod_{A\# G})$
has the form $P \iso Q\down({}_{x^{(a)}}(A\# G))$.
Then by \eqref{eq:explicit-Qdown}, we have
$$
\begin{aligned}
P(y) &\iso  (Q\down({}_{x^{(a)}}(A\# G)))(y)
= \tDs_{b \in G}\bi{(A\# G)}{x^{(a)}}{y^{(b)}}
\iso \tDs_{b \in G}\bi{A}{x}{y}^{a\inv b}\\
&= \la_a(\tDs_{b \in G}\bi{A}{x}{y}^b)
= \la_a({}_xA)(y)
\end{aligned}
$$
for all $y \in A_0$.
Therefore, we have $P \iso \la_a({}_xA)$.
Since any finitely generated projective object 
in $\Ggr(\Mod_A)$ is expressed as a finite direct sum of
indecomposables,
$P$ has the form in the statement.
The similar proof works also for the remaining case.
\end{proof}

\section{Orbit bimodules}
Throughout this section $R = (R, X)$ and $S =( S, Y)$ are small $\k$-categories with $G$-actions,
and $E\colon R \to R/G$ and $F\colon S \to S/G$ are the canonical $G$-coverings.

\begin{dfn}
\label{dfn:orbit-bimod}
(1) Let $M = (M, \ph)$ be a $G$-invariant $S$-$R$-bimodule.
Then we form a $G$-graded $S/G$-$R/G$-bimodule $M/G$ as follows, which we call
the {\em orbit bimodule} of $M$ by $G$:
\begin{itemize}
\item
For each $(y, x) \in (S/G)_0\times (R/G)_0 = S_0 \times R_0$ we set
\begin{equation}
\label{eq:orbit-bimod}
\bi{(M/G)} yx:= \Ds_{a\in G}\bi{M}{y}{ax}.
\end{equation}
\item
For each $(y, x), (y', x') \in (S/G)_0\times (R/G)_0 = S_0 \times R_0$ and
each $(s, r)  \in \bi{(S/G)}{y'}{y} \times \bi{(R/G)}{x}{x'}$ we define a morphism
$$
\bi{(M/G)} sr \colon \bi{(M/G)} yx \to \bi{(M/G)}{y'}{x'}
$$
in $\kMod$ by
\begin{equation}
\label{eq:orbit-bimod-str}
\bi{(M/G)} sr(m):= s\cdot m\cdot r:= \left(\sum_{cba=d}s_c\cdot \ph_c(m_b)\cdot cb r_a\right)_{d \in G}
\end{equation}
for all $r = (r_a)_{a\in G} \in \Ds_{a\in G}\bi{R}{x}{ax'}, m = (m_b)_{b\in G} \in \Ds_{b\in G}\bi{M}{y}{bx}$, and
$s = (s_c)_{c\in G} \in \Ds_{c\in G}\bi{S}{y'}{cy}$.
By the naturality of $\ph_a$ ($a \in G$) we easily see that \eqref{eq:orbit-bimod-str} defines
an $(S/G)$-$(R/G)$-\linebreak[3]bimodule structure on $M/G$.
\item
We set $\bi{(M/G)} yx^a:= \bi{M}{y}{ax}$ for all $a\in G$ and all
$(y,x) \in S_0 \times R_0$.
We easily see that this defines a $G$-grading
on $M/G$ by \eqref{eq:orbit-bimod} and \eqref{eq:orbit-bimod-str}.
\end{itemize}

(2) Let $f \colon M \to N$ be in $\Ginv(\biMod{S}{R})$.
For each $(y, x) \in S_0 \times R_0$ we set
$$
\bi{(f/G)}yx:= \Ds_{a\in G}\bi{f}{y}{ax}.
$$
Then as is easily seen $f/G:= (\bi{(f/G)}yx)_{(y,x)\in S_0\times R_0}$ turns out to be
a morphism $M/G \to N/G$ in $\Ggr(\biMod{S/G}{R/G})$.
\end{dfn}

\begin{prp}
\label{prp:fun-orbit}
The correspondences in Definition \ref{dfn:orbit-bimod} define
a $\k$-functor
$$
?/G \colon \Ginv(\biMod{S}{R}) \to \Ggr(\biMod{S/G}{R/G}).
$$
\end{prp}

\begin{proof}
The $\k$-linearity of $?/G$ is clear.
Let $(M, \ph )\in \Ginv(\Bi(S,\Mod,R))_0$.
Then $\id_M/G = (\Ds_{a\in G} \Bi(y,(\id_M), ax))_{(y,x)\in S_0\times R_0}$, each $(y,x)$-entry of which is the identity of
$\Bi(y,(M/G),x) = \Ds_{a\in G}\Bi(y,M,ax)$.
Hence $\id_M/G$ is the identity of the $G$-graded $S/G$-$R/G$-bimodule $M/G$.

Let $(L,\th) \ya{f} (M, \ph) \ya{g} (N, \ps)$ be morphisms
in $\Ginv(\Bi(S,\Mod,R))$.
We have to show that $(gf)/G = (g/G)(f/G)$.
It is enough to show that for any $(y,x) \in S_0 \times R_0$,
we have $\Bi(y,((gf)/G),x) = \Bi(y,(g/G),x) \cdot \Bi(y,(f/G),x)$,
which is shown as follows:
$$
\begin{aligned}
\text{RHS} &= \left(\Ds_{a \in G} \Bi(y,g,ax) \right)\cdot
\left(\Ds_{a \in G} \Bi(y,f,ax) \right)\\
&= \Ds_{a \in G} (\Bi(y,g,ax)\cdot\Bi(y,f,ax))
= \Ds_{a \in G} \Bi(y,(gf),ax)
= \text{LHS}.
\end{aligned}
$$
\end{proof}

Throughout the rest of this section, to regard $R/G$ (reap.\ $S/G$) as an $R$-$R$-bimodule
(respect.\ $S$-$S$-bimodule) we use
the canonical $G$-covering functor $E \colon R \to R/G$
(resp.\ $F \colon S \to S/G$).
Similarly, any $(S/G)$-$(R/G)$-bimodule $M$ is regarded as $S$-$R$-bimodule via $E$ and $F$.

\begin{lem}
\label{lem:tensor-R-R/G}
We have isomorphisms
$$
R/G \ox_R R/G \iso R/G \ox_{R/G} R/G \iso R/G
$$
of $(R/G)$-$(R/G)$-bimodules.
\end{lem}

\begin{proof}
Let $y, x \in R_0$.
Then we have
\[
\begin{aligned}
\bi{(R/G \ox_R R/G)}yx &= \bi{(R/G)_E \ox_R {}_E(R/G)}yx\\
&= \left.\left(\Ds_{z\in R_0}\bi{R/G}y{Ez} \ox_\k \bi{R/G}{Ez}x\right)\right/\bi{I}yx, \text{ and}\\
\bi{(R/G \ox_{R/G} R/G)}yx &= {}_yR/G \ox_{R/G} R/G_x\\
&= \left.\left(\Ds_{z\in R_0}\bi{R/G}y{Ez} \ox_\k \bi{R/G}{Ez}x\right)\right/\bi{I'}yx,
\end{aligned}
\]
where
\[
\begin{aligned}
\bi{I}yx:&= \ang{h\ox E(r)f - hE(r) \ox f \mid (h,r,f) \in \bi{R/G}{y}{z'}\times \bi{R}{z'}{z}\times \bi{R/G}zx\,, \, z, z' \in R_0},\text{ and}\\
\bi{I'}yx&= \ang{h\ox gf - hg \ox f \mid (h,g,f) \in \bi{R/G}{y}{z'}\times \bi{R/G}{z'}{z}\times \bi{R/G}zx\,, \, z, z' \in R_0}.
\end{aligned}
\]
Therefore it is enough to show that $\bi{I}yx = \bi{I'}yx$.
Since $E(r) \in \bi{R/G}{z'}{z}$ for all $r \in \bi{R}{z'}{z}$,
it is obvious that $\bi{I}yx \subseteq \bi{I'}yx$.
To show the converse inclusion, it suffices to show that
$h\ox gf - hg \ox f \in \bi{I}yx$
for all $(h,g,f) \in \bi{R/G}{y}{z'}\times \bi{R/G}{z'}{z}\times \bi{R/G}zx$ and all $z, z' \in R_0$.
Write $f:= (f_a)_{a\in G} \in \Ds_{a\in G} \bi{R}{z}{ax}$,
$g := (g_b)_{b \in G} \in \Ds_{b \in G}\bi{R}{z'}{bz}$ and
$h = (h_c)_{c \in G} \in \Ds_{c \in G}\bi{R}{y}{cz'}$.
Then since
\[
\bi{E}{z'}{bz} \colon \bi{R}{z'}{bz} \to \bi{(R/G)}{z'}{bz} = \Ds_{d \in G}\bi{R}{z'}{dbz}
= \Ds_{d \in G}\bi{R}{z'}{dz} = \bi{(R/G)}{z'}{z}
\]
is defined by $\bi{E}{z'}{bz}(g_b):= 
(\de_{d,b}g_b)_{d\in G}$, we have
\[
g = (g_b)_{b \in G} = \sum_{b\in G} (\de_{d,b}g_b)_{d\in G} = \sum_{b \in G}\bi{E}{z'}{dz}(g_b).
\]
Therefore we have
\[
\begin{aligned}
h\ox gf - hg \ox f &= h\ox \left(\sum_{b \in G}\bi{E}{z'}{dz}(g_b)\right)f - h\left(\sum_{b \in G}\bi{E}{z'}{dz}(g_b)\right) \ox f\\
&= \sum_{b \in G} (h\ox \bi{E}{z'}{dz}(g_b)f - h\bi{E}{z'}{dz}(g_b) \ox f) \in \bi{I}yx.
\qedhere
\end{aligned}
\]

\end{proof}

\begin{prp}
\label{prp:orbit-bimod}
Let $M$ be a $G$-invariant $S$-$R$-bimodule.
Then
\begin{enumerate}
\item
$M \ox_R (R/G) \iso {}_FM/G$ as $S$-$(R/G)$-bimodules; and
\item
$(S/G) \ox_S M \iso M/G_E$ as $(S/G)$-$R$-bimodules.
\end{enumerate}
Hence in particular, we have  
\begin{enumerate}
\item[(3)]
$
M \ox_R (R/G)_E \iso {}_FM/G_E \iso {}_F(S/G) \ox_S M
$
as $S$-$R$-bimodules, and
\item[(4)]
$(S/G)\ox_S M \ox_R (R/G) \iso M/G$
as $G$-graded $(S/G)$-$(R/G)$-bimodules.
\end{enumerate}
\end{prp}

\begin{proof}
Let $(y, x) \in S_0 \times R_0$.

(1) We have the following isomorphisms natural in $x, y$:
$$
\begin{aligned}
\bi{(M \ox_R (R/G))}{y}{x} &= \bi{M}y{} \ox_R \bi{(R/G)}{E}{x}
= \bi{M}y{} \ox_R \left(\Ds_{a\in G} \bi{R}{}{ax}\right)\\
&\iso \Ds_{a\in G} \bi{M}y{} \ox_R \bi{R}{}{ax}
\iso \Ds_{a\in G} \bi{M}{y}{ax}\\
&= \bi{(M/G)}{Fy}{x}.
\end{aligned}
$$
Hence $M \ox_R (R/G) \iso 
{}_F(M/G)$ as $S$-$(R/G)$-bimodules.

(2) Similarly we have the following isomorphisms natural in $x, y$:
$$
\begin{aligned}
\bi{(S/G \ox_S M)}{y}{x} &= \bi{(S/G)}y{F} \ox_S \bi{M}{}{x}
=  \left(\Ds_{a\in G} \bi{S}{y}{Y_a}\right)\ox_S \bi{M}{}x\\
&\iso  \left(\Ds_{a\in G} \bi{S}{a\inv y}{}\right)\ox_S \bi{M}{}x
\iso  \Ds_{a\in G} \bi{M}{a\inv y}x\\
&\iso  \Ds_{a\in G} \bi{M}{y}{ax}
= \bi{(M/G)}{y}{Ex}.
\end{aligned}
$$
Hence $(S/G) \ox_S M \iso M/G_E$ as $(S/G)$-$R$-bimodules.

(3) This follows from (1) and (2).

(4) We have the following isomorphisms of $S$-$R$-bimodules:
\[
\begin{aligned}
(S/G)\ox_S M \ox_R (R/G) &\iso M/G \ox_R (R/G)\quad (\text{by } (2))\\
&\iso M \ox_R (R/G) \ox_R (R/G)\quad (\text{by } (1))\\
&\iso M \ox_R (R/G) \ox_{R/G} (R/G) \quad (\text{by } {\rm Lemma\ \ref{lem:tensor-R-R/G}})\\
&\iso M \ox_R (R/G)\\
&\iso M/G\quad (\text{by } (1)).
\end{aligned}
\]
Hence we have $(S/G)\ox_S M \ox_R (R/G) \iso M/G$ also as $(S/G)$-$(R/G)$-bimodules because each $S$-$R$-bimodule morphism between $(S/G)$-$(R/G)$-bimodules is a $S/G$-$R/G$-bimodule morphism.
\end{proof}

\begin{lem}
\label{lem:G-inv-generate}
For each $M \in \Ginv(\biMod{S}{R})_0$, there exists a small set $I$ and
a family $(x_i, y_i)_{i\in I} \in (R_0 \times S_0)^I$ such that
there exists an epimorphism
\begin{equation}
\label{eq:G-inv-generate}
\Ds_{i \in I} \left(\Ds_{a \in G} S_{ay_i}\ox_\k {}_{ax_i}R\right) \to M
\end{equation}
in the category $\Ginv(\biMod{S}{R})$.
\end{lem}

\begin{proof}
Since $M$ can be regarded as a left $S \ox_\k R\op$-module,
there exist a small set $I$ and
a family $(x_i, y_i)_{i\in I} \in (R_0 \times S_0)^I$ such that
there exists an epimorphism
\[
f \colon \Ds_{i \in I}  S_{y_i}\ox_\k {}_{x_i}R \to M
\]
in the category $\biMod{S}{R}$.
Thus for each pair $(x,y) \in R_0 \times S_0$, the linear map
\[
\bi{f}{y}{x} = (\bi{f}{y}{x}^i)_{i \in I} \colon \Ds_{i \in I} \bi{S}{y}{y_i}\ox_\k \bi{R}{x_i}{x} \to \bi{M}{y}{x}
\]
is an epimorphism.
For each $a \in G$, we define a linear map
\[
\bi{f}{y}{x}^{i,a} \colon \bi{S}{y}{ay_i}\ox_\k \bi{R}{ax_i}{x} \to \bi{M}{y}{x}
\]
by the following commutative diagram:
\[
\begin{tikzcd}[column sep=50pt]
\Nname{S1}\Bi(y,S,ay_i) \ox_\k \Bi(ax_i,R,x) & \Nname{M1}\Bi(y,M,x)\\
\Nname{S2}\Bi(a\inv y,S,y_i)\ox_\k \Bi(x_i,R,a\inv x) & \Nname{M2}\Bi(a\inv y, M, a\inv x)
\Ar{S2}{S1}{"Y_a \ox_\k X_a", "\iso" '}
\Ar{M2}{M1}{"\ph_a^M" ', "\iso"}
\Ar{S1}{M1}{"{\Bi(y,f,x)^{i,a}}" , dashed}
\Ar{S2}{M2}{"{\Bi(a\inv y,f,a\inv x)^{i}}" '}
\end{tikzcd}.
\]
Then it is easy to verify the commutativity of the following diagram:
\[
\begin{tikzcd}[column sep=50pt]
\Nname{S1}\Ds_{a\in G}\Bi(y,S,ay_i) \ox_\k \Bi(ax_i,R,x) & \Nname{M1}\Bi(y,M,x)\\
\Nname{S2}\Ds_{a \in G}\Bi(by,S,ay_i)\ox_\k \Bi(ax_i,R, bx) & \Nname{M2}\Bi(b y, M, b x)
\Ar{S1}{S2}{"\ph_b" '}
\Ar{M1}{M2}{"\ph_b^M"}
\Ar{S1}{M1}{"{(\Bi(y,f,x)^{i,a})_{a \in G}}"}
\Ar{S2}{M2}{"{(\Bi(by,f,bx)^{i,a})_{a\in G}}" '}
\end{tikzcd},
\]
where $\ph_b$ is defined in \eqref{eqn:can-G-inv}.
Set
$\Bi(y,\bar{f},x)^i:= (\Bi(y,f,x)^{i,a})_{a \in G}$, and define an $S$-$R$-bimodule morphism
\[
\bar{f}^i \colon \Ds_{a \in G} S_{ay_i}\ox_\k {}_{ax_i}R \to M
\]
by $\bar{f}^i:= (\Bi(y,\bar{f},x)^i)_{(x,y)\in R_0\times S_0}$.
The commutativity of the diagram above shows that $\bar{f}^i$ is a morphism
in $\Ginv(\biMod{S}{R})$ for all $i \in I$, and we obtain
a morphism
\[
\bar{f}:= (\bar{f}^i)_{i\in I} \colon \Ds_{i\in I}\left(\Ds_{a \in G} S_{ay_i}\ox_\k {}_{ax_i}R\right) \to M
\]
in $\Ginv(\biMod{S}{R})$,
which turns out to be an epimorphism in $\Ginv(\biMod{S}{R})$
because so is $f$ in $\biMod{S}{R}$.
\end{proof}

\begin{prp}
\label{prp:Ginv-prj}
Assume that $G$ is a finite group such that $|G|$ is not divided by the
characteristic $p$ of $\k$.
Let $\Fgt \colon \Ginv(\biMod{S}{R}) \to \biMod{S}{R}$ be the forgetful
functor $(M, \ph) \mapsto M$, and $P \in \Ginv(\biMod{S}{R})_0$.
Then the following are equivalent.
\begin{enumerate}
\item
$P$ is projective in $\Ginv(\biMod{S}{R})$.
\item
$\Fgt(P)$ is projective in $\biMod{S}{R}$.
\item
$P$ is a direct summand of an object of the form
\begin{equation}
\label{eq:free-bimod}
\Ds_{i\in I}\left(\Ds_{a \in G} S_{ay_i}\ox_\k {}_{ax_i}R\right)
\end{equation}
in $\Ginv(\biMod{S}{R})$ for some family
$(x_i, y_i)_{i \in I} \in (R_0 \times S_0)^I$
with $I$ a small set,
where for each $i \in I$, 
$\Ds_{a\in G}S_{ay}\ox_\k {}_{ax}R$ is assumed to have the canonical $G$-invariant structure defined by the formula \eqref{eqn:can-G-inv}.
\end{enumerate}
\end{prp}

\begin{proof}
(1)\implies (2).
By Lemma \ref{lem:G-inv-generate}, we have an epimorphism
\eqref{eq:G-inv-generate}, where $M$ is replaced by $P$.
Assume that $P$ is projective in $\Ginv(\biMod{S}{R})$.
Then this epimorphism becomes a retraction
and mapped to a retraction by $\Fgt$.
Hence $\Fgt(P)$ turns out to be a direct summand of
$\Fgt(\Ds_{i\in I}\Ds_{a\in G} S_{ay_i}\ox_\k {}_{ax_i}R)$,
which is projective in $\biMod{S}{R}$.
Thus $P$ is projective in $\biMod{S}{R}$.

(2)\implies (3).
For each $M \in \Ginv(\biMod{S}{R})$, we write $M = (M, \ph^M)$ to stress
that $M$ is a $G$-invariant bimodule.
Then $\Fgt(M, \ph^M) = M$.
Assume that $\Fgt(P, \ph^P) = P$ is projective in $\biMod{S}{R}$.
We have to show that $(P, \ph^P)$ is projective in $\Ginv(\biMod{S}{R})$.
As discussed in the proof of Lemma \ref{lem:G-inv-generate},
there exist a small set $I$ and
a family $(x_i, y_i)_{i\in I} \in (R_0 \times S_0)^I$ such that
there exists an epimorphism
\[
f \colon \Ds_{i \in I}  S_{y_i}\ox_\k {}_{x_i}R \to P
\]
in the category $\biMod{S}{R}$.
As in the proof of Lemma \ref{lem:G-inv-generate},
construct an epimorphism
\[
\bar{f}:= (\bar{f}^i)_{i\in I} \colon \Ds_{i\in I}\left(\Ds_{a \in G} S_{ay_i}\ox_\k {}_{ax_i}R\right) \to P
\]
in $\Ginv(\biMod{S}{R})$.
Since $P$ is projective in $\biMod{S}{R}$, $f$ has a section
$s \colon P \to \Ds_{i \in I}  S_{y_i}\ox_\k {}_{x_i}R$ such that
$fs = \id_P$ in $\biMod{S}{R}$.
Thus for each pair $(x,y) \in R_0 \times S_0$, the linear map
\[
\bi{s}{y}{x} = (\bi{s}{y}{x}^i)_{i \in I} \colon 
\bi{P}{y}{x} \to \Ds_{i \in I} \bi{S}{y}{y_i}\ox_\k \bi{R}{x_i}{x} 
\]
satisfies $\bi{f}{y}{x}\cdot \bi{s}{y}{x} = \id_{\bi{P}{y}{x}}$ in $\kMod$.
For each $a \in G$, we define a linear map
\[
\bi{s}{y}{x}^{i,a} \colon \bi{P}{y}{x} \to
\bi{S}{y}{ay_i}\ox_\k \bi{R}{ax_i}{x}
\]
by the following commutative diagram:
\[
\begin{tikzcd}[column sep=50pt]
\Nname{P1}\Bi(y,P,x) &
   \Nname{S1}\Bi(y,S,ay_i) \ox_\k \Bi(ax_i,R,x)\\
\Nname{P2}\Bi(a\inv y, P, a\inv x) &
   \Nname{S2}\Bi(a\inv y,S,y_i)\ox_\k \Bi(x_i,R,a\inv x)
\Ar{S2}{S1}{"Y_a \ox_\k X_a" ', "\iso"}
\Ar{P2}{P1}{"\ph_a^P", "\iso" '}
\Ar{P1}{S1}{"{\Bi(y,s,x)^{i,a}}" , dashed}
\Ar{P2}{S2}{"{\Bi(a\inv y,s,a\inv x)^{i}}" '}
\end{tikzcd}.
\]
Then as is easily seen, we have a commutative diagram
\[
\begin{tikzcd}[column sep=50pt]
\Nname{P1}\Bi(y,P,x)& 
\Nname{S1}\Ds_{a\in G}\Bi(y,S,ay_i) \ox_\k \Bi(ax_i,R,x)\\
\Nname{P2}\Bi(b y, P, b x)&
\Nname{S2}\Ds_{a \in G}\Bi(by,S,ay_i)\ox_\k \Bi(ax_i,R, bx)
\Ar{S1}{S2}{"\ph_b"}
\Ar{P1}{P2}{"\ph_b^P" '}
\Ar{P1}{S1}{"{(\Bi(y,s,x)^{i,a})_{a \in G}}"}
\Ar{P2}{S2}{"{(\Bi(by,s,bx)^{i,a})_{a\in G}}" '}
\end{tikzcd},
\]
where $\ph_b$ is defined in \eqref{eqn:can-G-inv}.
Set
$\Bi(y,\bar{s},x)^i:= {}^t(\Bi(y,s,x)^{i,a})_{a \in G}$, and define an $S$-$R$-bimodule morphism
\begin{equation}
\label{eq:s-bar}
\bar{s}^i \colon P \to \Ds_{a \in G} S_{ay_i}\ox_\k {}_{ax_i}R
\end{equation}
by $\bar{s}^i:= (\Bi(y,\bar{s},x)^i)_{(x,y)\in R_0\times S_0}$.
Here note that we used the assumption that
$G$ is a finite group to have the direct sum in \eqref{eq:s-bar}.
The commutativity of the diagram above shows that
$\bar{s}^i$ is a morphism
in $\Ginv(\biMod{S}{R})$ for all $i \in I$, and we obtain
a morphism
\[
\bar{s}:= (\bar{s}^i)_{i\in I} \colon P \to \Ds_{i\in I}\left(\Ds_{a \in G} S_{ay_i}\ox_\k {}_{ax_i}R\right)
\]
in $\Ginv(\biMod{S}{R})$.
It is not hard to see that $\bar{f}\bar{s} = |G|\id_{P}$.
Since $p\nmid |G|$, we have $\bar{f}\bar{s}$ is an isomorphism.
Thus $(P, \ph^P)$ is a direct summand of
$\Ds_{i\in I}\left(\Ds_{a \in G} S_{ay_i}\ox_\k {}_{ax_i}R\right)$
in $\Ginv(\biMod{S}{R})$.

(3)\implies (1).
It is standard to show that an object of the form
\eqref{eq:free-bimod} is projective in
$\Ginv(\biMod{S}{R})$.
Hence $(P, \ph^P)$ is projective in $\Ginv(\biMod{S}{R})$.
\end{proof}

\begin{dfn}
\label{dfn:fg-G-inv}
Let $M$ be an $S$-$R$-bimodule.
$M$ is called a \emph{finitely generated $G$-invariant} bimodule
(\emph{f.g. $G$-invariant} for short) if there exists
a finite set $I$ and an
epimorphism stated in Lemma \ref{lem:G-inv-generate}.
\end{dfn}

\begin{dfn}
\label{dfn:1sided-fgp}
Let $\calC$ and $\calD$ be $\k$-categories, and
$M$ a $\calD$-$\calC$-bimodule.
Then we say that $M$ is \emph{finitely generated projective} as a right $\calC$-module
(or $M_\calC$ is finitely generated projective for short) if ${}_yM$ is finitely generated projective right $\calC$-module for all $y \in \calD_0$.

Similarly, we say that $M$ is \emph{finitely generated projective} as a left $\calD$-module (or ${}_{\calD}M$ is finitely generated projective for short)
if $M_x$ is a finitely generated projective
left $\calD$-module for all $x \in \calC_0$.
\end{dfn}

By Proposition \ref{prp:orbit-bimod} (1) and (2) we obtain the following.

\begin{cor}\label{cor:proj-orbit}
Let $M$ be a $G$-invariant $S$-$R$-bimodule.
Then the following statements hold.
\begin{enumerate}
\item
Assume that $S$ is $\k$-projective.
If $M_R$ is finitely generated projective,
then so is $M/G_{R/G}$.

\item
Assume that $R$ is $\k$-projective.
If ${}_SM$ is finitely generated projective,
then so is ${}_{S/G}M/G$.
\end{enumerate}
\end{cor}

\begin{proof}
By Lemma \ref{lem:G-inv-generate},
we have an epimorphism
\[
f \colon \Ds_{i\in I}\Ds_{a \in G} S_{ay_i}\ox_\k {}_{ax_i}R \to M
\]
in $\Ginv({}_S\Mod_R)$ for some small set $I$ and a family $(x_i,y_i)_{i\in I} \in (R_0 \times S_0)^I$.
This yields an epimorphism
\[
{}_yf \colon \Ds_{i \in I}\Ds_{a \in G} {}_yS_{ay_i}\ox_\k {}_{ax_i}R \to {}_yM
\]
of right $R$-modules.
Since $S$ is $\k$-projective, each ${}_yS_{ay_i} \ox_\k {}_{ax_i}R$
is a projective $R$-module.
Since ${}_yM$ is a projective $R$-module, ${}_yf$ turns out to be a retraction.
By applying the functor $?\ox_R R/G$ $f$,
we obtain a retraction
\[
{}_yf \ox_R R/G \colon \Ds_{i \in I}\Ds_{a \in G} ({}_yS_{ay_i}\ox_k {}_{ax_i}R)\ox_R R/G \to {}_yM \ox_R R/G
\]
whose domain is a projective right $R/G$-module.
Hence by Proposition \ref{prp:orbit-bimod}(1),
${}_y(M/G) \iso {}_y(M \ox_R R/G) = {}_yM \ox_R R/G$ is a projective
$R/G$-module.
Since $M_R$ is finitely generated,
so is ${}_yM \ox_R R/G$ over $R/G$.
As a consequence, $M_R$ is finitely generated projective.

(2) This is proved similarly.
\end{proof}

\begin{prp}\label{prp:tensor-orbit}
Let $T = (T, Z)$ be a small $\k$-category with $G$-action, and
${}_TN_S$, ${}_SM_R$ be G-invariant bimodules.
Then
\begin{enumerate}
\item
${}_T(N \ox_S M)_R$ is a $G$-invariant bimodule.
\item
$(N \ox_S M)/G \iso (N/G)\ox_{S/G}(M/G)$ in $\Ggr(\biMod{T/G}{R/G})$.
\end{enumerate}
\end{prp}

\begin{proof}
(1) Let $(y,x) \in S_0 \times R_0$ and $a \in G$.
We set
\[\begin{aligned}
\bi{I}yx&:=\ang{v\ox su - vs \ox u \mid (v,s,u) \in \bi{N}{y}{z'}\times \bi{S}{z'}{z}\times \bi{M}zx, z, z' \in S_0}\\
&\subseteq \Ds_{z\in S_0}\bi{N}yz \ox_\k \bi{M}zx, \text{ and}\\
\bi{(\ph^{N\ox_{\k} M}_a)}yx&:= \Ds_{z \in S_0}\bi{(\ph^N_a)}yz \ox_\k \bi{(\ph^M_a)}zx \colon
\Ds_{z \in S_0}\bi{N_z\ox_{\k}{}_zM}yx \to \Ds_{z \in S_0}\bi{N_{az} \ox_{\k}{}_{az}M}{ay}{ax}.
\end{aligned}
\]
Here we write $\bi{(\ph_a)}yx:=\bi{(\ph^{N\ox_{\k} M}_a)}yx$ for short,
and it is an isomorphism in $\kMod$.
Note that for each $(v,s,u) \in \bi{N}{y}{z'}\times \bi{S}{z'}{z}\times \bi{M}zx\ (z, z' \in S_0)$,
we have
$$
\begin{aligned}
&\bi{(\ph_a)}yx(v\ox su - vs \ox u) = \bi{(\ph^N_a)}y{z'}(v)\ox \bi{(\ph^M_a)}{z'}{x}(su) - \bi{(\ph^N_a)}yz(vs) \ox \bi{(\ph^M_a)}zx(u)\\
&\hspace{6em}= \bi{(\ph^N_a)}{y}{z'}(v)\ox (as)\cdot\bi{(\ph^M_a)}{z}{x}(u) - \bi{(\ph^N_a)}{y}{z'}(v)\cdot (as)\ox \bi{(\ph^M_a)}{z}{x}(u),
\end{aligned}
$$
Hence we see that (we omit subscripts $y, z', z, x$ for $\ph_a^M, \ph_a^N$ below for simplicity)
$$
\begin{aligned}
&\bi{(\ph_a)}yx (\bi{I}yx) \\
&= \langle \ph^N_a(v)\ox (as)\cdot\ph^M_a(u) - \ph^N_a(v)\cdot (as)\ox \ph^M_a(u)\mid\\
&\hspace{5cm}(v,s,u) \in \bi{N}{y}{z'}\times \bi{S}{z'}{z}\times \bi{M}zx,\ z, z' \in S_0\rangle\\
&=
\ang{v'\ox s'u' - v's' \ox u' \mid (v',s',u') \in \bi{N}{ay}{az'}\times \bi{S}{az'}{az}\times \bi{M}{az}{ax},\ az, az' \in S_0}\\
&=
\ang{v'\ox s'u' - v's' \ox u' \mid (v',s',u') \in \bi{N}{ay}{t'}\times \bi{S}{t'}{t}\times \bi{M}{t}{ax},\ t, t' \in S_0}\\
&=\bi{I}{ay}{ax}.
\end{aligned}
$$
Therefore the isomorphism $\bi{(\ph_a)}yx$ induces the following
isomorphism $\ovl{\bi{(\ph_a)}yx}$:
\begin{equation}
\label{eq:tensor}
\begin{aligned}
\bi{(N\ox_S M)}yx =&\ {}_y N \ox_S M_x\\
=&\ \textstyle(\Ds_{z\in S_0}\bi{N}yz \ox_\k \bi{M}zx)/\bi{I}yx\\
\ya{\ovl{\bi{(\ph_a)}yx}}& \ \textstyle(\Ds_{z\in S_0}\bi{N}{ay}{az} \ox_\k \bi{M}{az}{ax})/\bi{I}{ay}{ax}\\
=&\  \textstyle(\Ds_{z\in S_0}\bi{N}{ay}{z} \ox_\k \bi{M}{z}{ax})/\bi{I}{ay}{ax}\\
=&\ {}_{ay}N\ox_S M_{ax}\\
=&\ \bi{(N\ox_S M)}{ay}{ax}.
\end{aligned}
\end{equation}
It is easy to verify that $\ovl{\ph_a}:= (\ovl{\bi{(\ph_a)}yx})_{(y,x)}$ is a natural transformation and that
$\ovl{\ph}:= (\ovl{\ph_a})_{a\in G}$ makes $N \ox_S M$ a $G$-invariant $T$-$R$-bimodule.

(2)
We have the following isomorphisms of $G$-graded $(T/G)$-$(R/G)$-bimodules:
$$
\begin{aligned}
(N \ox_S M)/G &\iso (T/G)\ox_T N \ox_S M\ox_R (R/G)\\
&\iso (T/G)\ox_T N \ox_S {}_F(M/G)\\
&\iso (T/G)\ox_T N \ox_S (S/G)\ox_S M\ox_R (R/G)\\
&\iso (T/G)\ox_T N \ox_S (S/G)\ox_{S/G} (S/G)\ox_S M\ox_R (R/G)\\
&\iso (N/G) \ox_{S/G} (M/G).
\end{aligned}
$$ 
\end{proof}

\begin{dfn}
\label{dfn:fgp-Ginv-bimod}
(1) An $S$-$R$-bimodule $P$ is called
{\em finitely generated projective $G$-invariant}
(\emph{f.g.\ projective $G$-invariant} for short)
if it is a direct summand of the form
$$
\Ds_{i\in I}\left(\Ds_{a \in G} S_{ay_i}\ox_\k {}_{ax_i}R\right)
$$
in $\Ginv(\biMod{S}{R})$ with $I$ a finite set
(see Proposition \ref{prp:Ginv-prj} in the case where
$G$ is a finite group such that $|G|$ is not divided
by the characteristic of $\k$).

(2) A $B$-$A$-bimodule $P$ is called {\em finitely generated projective $G$-graded} (\emph{f.g.\ projective $G$-graded} for short)
if it is a direct summand of an object of the form
$$
\Ds_{i=1}^n \ro_{b_i}(B_{y_i})\ox_\k \la_{a_i}({}_{x_i}A)
$$
in $\Ggr(\biMod{B}{A})$ for some integer $n \ge 0$,
$(y_i,x_i)_{i=1}^n \in (B_0 \times A_0)^n$,
and $(b_i,a_i)_{i=1}^n \in (G\times G)^n$,
$B_{y}\ox_\k {}_{x}A$ is assumed to have the canonical $G$-grading
defined by the formula \eqref{eq:can-G-gr}.
By Proposition \ref{prp:graded-prj-bimod},
$P$ is finitely generated projective $G$-graded in $\Ggr(\biMod{B}{A})$
if and only if $P$ is $G$-graded and
finitely generated projective as a $B$-$A$-bimodule.
\end{dfn}

\begin{prp}
\label{prp:orbitr-bi}
Let $(y,x) \in S_0 \times R_0$, and consider
the $G$-invariant module
$\Ds_{a\in G}S_{ay}\ox_\k {}_{ax}R$
having the canonical $G$-invariant structure \eqref{eqn:can-G-inv}.
Then
\begin{equation}
\label{eq:G-inv-slash-G}
(\textstyle\Ds_{a\in G}S_{ay}\ox_\k {}_{ax}R)/G \iso 
(S/G)_y \ox_\k {}_x(R/G),
\end{equation}
where the right hand side has the canonical $G$-grading \eqref{eq:can-G-gr}.
Therefore, if $P$ is a f.g.\  projective $G$-invariant $S$-$R$-bimodule,
then $P/G$ is a finitely generated projective $G$-graded $(S/G)$-$(R/G)$-bimodule.
\end{prp}

\begin{proof}
Set $P:= \Ds_{a\in G}S_{ay}\ox_\k {}_{ax}R$.
Then $P/G$ has the $G$-grading defined by ${}_v (P/G)_u = \Ds_{b\in G} {}_v (P/G)_u^b$
for all $(v,u) \in (S/G)_0 \times (R/G)_0$, where
\[
\begin{aligned}
{}_v (P/G)_u^b &= {}_vP_{bu} =\Ds_{a\in G}{}_v S_{ay}\ox_\k {}_{ax}R_{bu}\\
&\isoto \Ds_{a\in G}{}_v S_{ay}\ox_\k {}_{x}R_{a\inv b u}\\
&= \Ds_{a\in G}{}_v (S/G)_{y}^a\ox_\k {}_{x}(R/G)_{u}^{a\inv b}.
\end{aligned}
\]
Here, the isomorphism above is given by
$\Ds_{a \in G} \id_{\bi{S}{v}{ay}} \ox_\k \bi{(X_{a\inv})}{ax}{bu}$.
Therefore
\[
\begin{aligned}
P/G &= \Ds_{b \in G}(P/G)^b \\
&\iso \Ds_{b \in G}\left(\Ds_{a\in G}(S/G)_{y}^a\ox_\k {}_{x}(R/G)^{a\inv b}\right)
= (S/G)_y \ox_\k {}_x(R/G).
\end{aligned}
\]
Here, the latter has the canonical $G$-grading, and
the isomorphism above is given by
$\Ds_{b\in G}\Ds_{a \in G} \id_{\bi{S}{v}{ay}} \ox_\k \bi{(X_{a\inv})}{ax}{bu}$, which is degree preserving.

Since the functor $?/G$ is additive,
the remaining statement follows.
\end{proof}

\section{Smash products}
Throughout this section $A$ and $B$ are $G$-graded small $\k$-categories.

\begin{dfn}
(1) Let $M$ be a $G$-graded $B$-$A$-bimodule.
Then we define a $G$-invariant $(B\# G)$-$(A\# G)$-bimodule $M\# G$ as follows,
which we call the {\em smash product} of $M$ and $G$:
\begin{itemize}
\item
For each $(y^{(b)}, x^{(a)}) \in (B\# G)_0 \times (A\# G)_0$ we set
\begin{equation}
\label{eq:Hom-smash}
\bi{(M\# G)}{y^{(b)}}{x^{(a)}}:= \bi{M}{y}{x}^{b\inv a}.
\end{equation}
\item
For each $(y^{(b)}, x^{(a)}), ({y'}^{(b')}, {x'}^{(a')}) \in (B\# G)_0 \times (A\# G)_0$ and
each $(\Biup((b'),\be,{(b)}), \linebreak[3] \Biup((a),\al,{(a')})) \in \bi{(B\# G)}{{y'}^{(b')}}{y^{(b)}} \times \bi{(A\# G)}{x^{(a)}}{{x'}^{(a')}}$
with $(\be, \al) \in \bi{B}{y'}{y}^{{b'}\inv b} \times \bi{A}{x}{x'}^{a\inv a'}$,
we define a morphism $\bi{(M\# G)}{\be}{\al}$ in $\kMod$ by the following commutative diagram:
\begin{equation}
\label{eq:mor-smash}
\xymatrix@C=120pt{
\bi{(M\# G)}{y^{(b)}}{x^{(a)}} & \bi{(M\# G)}{{y'}^{(b')}}{{x'}^{(a')}}\\
\bi{M}{y}{x}^{b\inv a} & \bi{M}{y'}{x'}^{{b'}\inv a'}.
\ar^{\bi{(M\# G)}{\left(\Biup((b'),\be,{(b)})\right)}{\left(\Biup((a),\al,{(a')})\right)}}"1,1";"1,2"
\ar@{=}"1,1";"2,1"
\ar@{=}"1,2";"2,2"
\ar_{\bi{M}{\be}{\al}}"2,1";"2,2"
}
\end{equation}

Since $\deg(\be)\deg(m)\deg(\al) = ({b'}\inv b)(b\inv a)(a\inv a') =  {b'}\inv a'$
for all $m \in {}_yM_x^{b\inv a}= \bi{(M\# G)}{y^{(b)}}{x^{(a)}}$,
the bottom morphism is well-defined.
It is easy to verify that this makes $M\# G$ a $(B\# G)$-$(A\# G)$-bimodule.
\item
For each $(y^{(b)}, x^{(a)}) \in (B\# G)_0 \times (A\# G)_0$  and each $c \in G$
we define $\bi{(\ph_c)}{y^{(b)}}{x^{(a)}}$ by the following commutative diagram:
\begin{equation}
\label{eq:str-smash}
\xymatrix@C=10ex{
\bi{(M\# G)}{y^{(b)}}{x^{(a)}} & \bi{(M\# G)}{c\cdot y^{(b)}}{c\cdot x^{(a)}}\\
\bi{M}{y}{x}^{b\inv a} & \bi{(M\# G)}{y^{(cb)}}{x^{(ca)}}.
\ar^(0.45){\bi{(\ph_c)}{y^{(b)}}{x^{(a)}}}"1,1";"1,2"
\ar@{=}"1,1";"2,1"
\ar@{=}"1,2";"2,2"
\ar@{=}"2,1";"2,2"
}
\end{equation}
Then by letting
$\ph_c:=(\bi{(\ph_c)}{y^{(b)}}{x^{(a)}})_{(y^{(b)}, x^{(a)})}$, and $\ph:=(\ph_c)_{c\in G}$,
we have a $G$-invariant $(B\# G)$-$(A\# G)$-bimodule $(M\# G, \ph)$.
\end{itemize}

(2)
Let $f \colon M \to N$ be in $\Ggr(\biMod{B}{A})$.
For each $(y^{(b)}, x^{(a)}) \in (B\# G)_0 \times (A\# G)_0$,
we define
$\bi{(f\# G)}{y^{(b)}}{x^{(a)}}$ by the commutative diagram
\begin{equation}
\label{eq:smash-on-morph}
\xymatrix@C=15ex{
\bi{(M\# G)}{y^{(b)}}{x^{(a)}} & \bi{(N\# G)}{y^{(b)}}{x^{(a)}}\\
\bi{M}{y}{x}^{b\inv a} & \bi{N}{y}{x}^{b\inv a}.
\ar^{\bi{(f\# G)}{y^{(b)}}{x^{(a)}}}"1,1";"1,2"
\ar@{=}"1,1";"2,1"
\ar@{=}"1,2";"2,2"
\ar_{f|_{\bi{M}{y}{x}^{b\inv a}} =\, \Bi(y,f,x)^{b\inv a}}"2,1";"2,2"
}
\end{equation}
Then as is easily seen $f\# G:= (\bi{(f\# G)}{y^{(b)}}{x^{(a)}})_{(y^{(b)}, x^{(a)})}$
is a morphism $M\# G \to N\# G$ in the category $\Ginv(\biMod{(B\# G)}{(A\# G)})$.
\end{dfn}

\begin{prp}
\label{prp:fun-smash}
The smash product construction above is extended to a $\k$-functor
$$
?\# G \colon \Ggr(\biMod{B}{A}) \to \Ginv(\biMod{(B\# G)}{(A\# G)}).
$$
\end{prp}

\begin{proof}
By looking at \eqref{eq:smash-on-morph},
the $\k$-linearity of $?\# G$ is clear.
To verify the fact that $?\# G$ preserves identities,
let $M\in \Ggr(\biMod{B}{A})_0$.
Then for each $(y^{(b)}, x^{(a)}) \in (B\# G)_0 \times (A\# G)_0$,
the diagram \eqref{eq:smash-on-morph} for $f:= \id_M$
shows that
$\bi{(\id_M\# G)}{y^{(b)}}{x^{(a)}} = 
\id_M|_{\bi{M}{y}{x}^{b\inv a}} = 
\id_{\bi{M}{y}{x}^{b\inv a}}$, as required.

Finally to verify that $?\# G$ preserves compositions,
let $L \ya{f} M \ya{g} N$ be morphisms
in $\Ggr(\biMod{B}{A})$.
It is enough to show the equality
$\Bi(y^{(b)},((gf)\# G),x^{(a)}) = \Bi(y^{(b)},(g\# G),x^{(a)}) \cdot \Bi(y^{(b)},(f\# G),x^{(a)})$
for all $(y^{(b)}, x^{(a)}) \in (B\# G)_0 \times (A\# G)_0$.
By definition, this follows from the obvious
equality
$\left(g|_{\bi{M}{y}{x}^{b\inv a}} \right)\cdot
\left(f|_{\bi{L}{y}{x}^{b\inv a}} \right)
= (gf)|_{\bi{L}{y}{x}^{b\inv a}}$.
\end{proof}

\section{Cohen-Montgomery duality for bimodules}

Throughout this section,
we let $R, S$ be small $\k$-categories with $G$-actions, and
$A, B$ be $G$-graded small $\k$-categories.

\begin{dfn}
\label{dfn:ze-Mod-ze}
By Remark \ref{rmk:omega-zeta-identify},
we have $G$-equivariant equivalences $\ze_S \colon S \to (S/G)\# G$ and
$\ze_R \colon R \to (R/G)\# G$ (see \cite[Definition 5.8.1]{Asa-book} or \cite[Definition 8.1]{Asa11} for definition)
having quasi-inverses
$\ze'_S \colon (S/G)\# G \to S$ and
$\ze'_R \colon  (R/G)\# G\to R$, respectively
(see \cite[Definition 5.8.5]{Asa-book} for definition).
We define a functor
$$
\biMod{\ze_S}{\ze_R} \colon
\biMod{(S/G)\# G}{(R/G)\# G} \to
\biMod{S}{R}
$$
as follows.
For each $M \in (\biMod{(S/G)\# G}{(R/G)\# G})_0$,
$\biMod{\ze_S}{\ze_R}(M)$ is the $S$-$R$-bimodule
$M \circ (\ze_S \times \ze_R\op)$ as in the diagram
$$
S \times R\op \ya{\ze_S \times \ze_R\op} ((S/G)\# G)\times ((R/G)\# G)
\ya{M} \kMod.
$$
The functor $\biMod{\ze_S}{\ze_R}$ has a quasi-inverse $\biMod{\ze'_S}{\ze'_R}$, and hence it is an equivalence.
This induces the equivalence
\begin{equation}
\label{eq:zeta-bar}
\ovl{\ze}:=\Ginv(\biMod{\ze_S}{\ze_R}) \colon
\Ginv(\biMod{(S/G)\# G}{(R/G)\# G}) \to
\Ginv(\biMod{S}{R}).
\end{equation}
\end{dfn}

\begin{dfn}
\label{dfn:om-Mod-om}
By Remark \ref{rmk:omega-zeta-identify},
we have strictly $G$-degree preserving equivalences $\om_B \colon B \to (B\# G)/G$
and
$\om_A \colon A \to (A\# G)/G$ (see \cite[Definition 5.8.8]{Asa-book} or \cite[Definition 8.5]{Asa11} for definition)
having quasi-inverses
$\om'_B \colon (B\# G)/G \to B$ and
$\om'_A \colon  (A\# G)/G \to A$, respectively
(see \cite[Definition 5.8.11]{Asa-book} for definition).
We can define a functor
$$
\biMod{\om_B}{\om_A} \colon
\biMod{(B\# G)/G}{(A\# G)/G} \to
\biMod{B}{A}
$$
as follows.
For each $N \in (\biMod{(B\# G)/G}{(A\# G)/G})_0$,
$\biMod{\om_B}{\om_A}(N)$ is the $B$-$A$-bimodule
$N \circ (\om_B \times \om_A\op)$ as in the diagram
$$
B \times A\op \ya{\om_B \times \om_A\op} ((B\# G)/G)\times ((A\# G)/G)
\ya{N} \kMod.
$$
The functor $\biMod{\om_B}{\om_A}$ has a quasi-inverse $\biMod{\om'_B}{\om'_A}$, and hence it is an equivalence.
This induces the equivalence
$$
\ovl{\om}:=\Ggr(\biMod{\om_B}{\om_A}) \colon
\Ggr(\biMod{(B\# G)/G}{(A\# G)/G}) \to
\Ggr(\biMod{B}{A}).
$$
Since $\om'_A$ and $\om'_B$ are $G$-degree preserving
(in a general sense) but
not necessarily strictly degree preserving,
we do not use $\om'$ and the fact that $\om$ is an
equivalence as much as possible.
For example, in the proof of Theorem \ref{thm:CM-dual-bimod} below,
we do not make use of $\om'$, and hence the terminology
``degree-preserving'' is used only in the strict sense. 
\end{dfn}

\begin{thm}
\label{thm:CM-dual-bimod}
\begin{enumerate}
\item
The functor $?/G \colon \Ginv(\biMod{S}{R}) \to \Ggr(\biMod{S/G}{R/G})$ is an
equivalence with a quasi-inverse given by the composite
$$
(?\# G)'\colon
\Ggr(\biMod{S/G}{R/G}) \ya{?\# G} \Ginv(\biMod{(S/G)\# G}{(R/G)\# G})
\ya{\ovl{\ze}} \Ginv(\biMod{S}{R}).
$$
\item
The functor $?\# G \colon \Ggr(\biMod{B}{A}) \to \Ginv(\biMod{(B\# G)}{(A\# G)})$
is an equivalence with a quasi-inverse given by the composite
$$
(?/G)'\colon
\Ginv(\biMod{(B\# G)}{(A\# G)}) \ya{?/G} \Ggr(\biMod{(B\# G)/G}{(A\# G)/G})
\ya{\ovl{\om}} \Ggr(\biMod{B}{A}).
$$
\item
In particular, for each  $G$-invariant bimodule ${}_RM_S$
we have $((M/G)\# G)' \iso M$ as $S$-$R$-bimodules,
and for each $G$-graded bimodule ${}_BM_A$ we have
$((M\# G)/G)' \iso M$ as $B$-$A$-bimodules.
\item
As a consequence, both $?\# G$ and $?/G$ preserve small colimits and small limits, and in particular, \textcolor{blue}{are} right exact and left exact.
\end{enumerate}
\end{thm}

\begin{proof}
It is enough to prove the statements (1) and (2).

(1)
As before, we set $(E, \ph^E) \colon R \to R/G$ and
$(F, \ph^F) \colon S \to S/G$ to be the canonical $G$-coverings.
Recall that the $G$-equivariant equivalence
$(\ze_R, \ovl{\ph^E}) \colon R \to (R/G)\# G$ sends each morphism
$f \colon x \to x'$ in $R$ to
$\Biup((1),(Ef), {(1)}) \colon (Ex)^{(1)} \to (Ex')^{(1)}$,
where $\ovl{\ph^E}$ is a $G$-invariant structure given as a family of isomorphisms
$\ovl{\ph^E_{a}}x:= \Biup((1),(\ph^E_a x) ,{(a)}) \colon (Ex)^{(a)} = a((Ex)^{(1)}) \to (E(ax))^{(1)}$
for each $x \in R_0$ and $a \in G$ such that the diagram
\begin{equation}
\label{eq:nat-iso-pRh}
\begin{tikzcd}[column sep=60pt]
\Nname{x'}(Ex')^{(a)} & \Nname{x}(Ex)^{(a)}\\
\Nname{ax'}(Eax')^{(1)} & \Nname{ax}(Eax)^{(1)}
\Ar{x'}{ax'}{"\ovl{\ph^E_a}x'" '}
\Ar{x}{ax}{"\ovl{\ph^E_a}x"}
\Ar{x'}{x}{"{{}^{(a)}(Ef)^{(a)}}"}
\Ar{ax'}{ax}{"{{}^{(1)}(Eaf)^{(1)}}" '}
\end{tikzcd}
\end{equation}
is commutative for all $f \in \Bi(x,R,x')$, $x,x' \in R_0$
(see \cite[Lemma 5.8.2]{Asa-book}).

Consider the following diagram that illustrates our setting:
$$
\begin{tikzcd}
\Nname{R}\Ginv(\biMod{S}{R}) &  & \Nname{R'}\Ggr(\biMod{S/G}{R/G})\\
\Nname{R1} \Ginv(\biMod{S}{R}) & & \Nname{R1'}\Ggr(\biMod{S/G}{R/G})
\Ar{R}{R'}{"?/G"}
\Ar{R}{R1}{equal}
\Ar{R'}{R1}{"(?\# G)'" '}
\Ar{R'}{R1'}{equal}
\Ar{R1}{R1'}{"?/G"}
\ar["\et", Rightarrow, to path={([xshift=10, yshift=-10]R.south)--([xshift=30, yshift=-10]R.south)\tikztonodes}]
\ar["\ep", Rightarrow, to path={([xshift=-30, yshift=-10]R'.south)--([xshift=-10, yshift=-10]R'.south)\tikztonodes}]
\end{tikzcd}.
$$
To show that $?/G$ is an equivalence with a quasi-inverse
$(?\# G)'$ it is enough to construct natural isomorphisms:
$$
\begin{aligned}
\et &\colon  \id_{\Ginv(\biMod{S}{R})} \To (?\# G)' \circ  (?/G), \text{ and}\\
\ep &\colon (?/G) \circ (?\# G)' \To \id_{\Ggr(\biMod{S/G}{R/G})}.
\end{aligned}
$$

Let $M = (M, \ph) \in \Ginv(\biMod{S}{R})_0$.
Since $\Bi(Fy,(M/G),Ex) = \Ds_{a\in G} \Bi(y,M,ax)$,
we have
$\Bi((Fy)^{(b)},((M/G)\# G), {(Ex)^{(a)}}) = \Bi(Fy,(M/G),Ex)^{b\inv a}
= \Bi(y,M,{(b\inv a)x})$
for any $(Fy)^{(b)}\in ((S/G)\# G)_0, {(Ex)^{(a)}} \in ((R/G)\# G)_0$.
Thus
$$
\Bi(y,((M/G)\# G)',x) = \Bi((Fy)^{(1)},((M/G)\# G), {(Ex)^{(1)}}) = \Bi(y,M,x).
$$
Let $f\colon M = (M, \ph)\to  M = (M', \ph') \in \Ginv(\biMod{S}{R})_1$, and $y\in S_0, x \in R_0$. Then we have a commutative diagram
$$
\begin{tikzcd}[column sep=80pt]
\Nname{1} \Bi(y, ((M/G)\# G)', x) & \Nname{2} \Bi(y,((M'/G)\# G)', x)\\
\Nname{1a}\Bi((Fy)^{(1)}, (M/G)\# G), {(Ex)^{(1)}}) & \Nname{2a}\Bi((Fy)^{(1)},  ((M'/G)\# G), {(Ex)^{(1)}}) \\
\Nname{1b}\Bi(Fy, ( M/G)^1, Ex) &\Nname{2b}\Bi(Fy, (M'/G)^1, Ex)\\
\Nname{1c}\Bi(y, M, x) & \Nname{2c}\Bi(y, M', x) \\
\Ar{1}{2}{"{\Bi(y,((f/G)\# G)',x)}"}
\Ar{1a}{2a}{"{\Bi((Fy)^{(1)}, ((f/G)\# G), {(Ex)^{(1)}})}"}
\Ar{1b}{2b}{"{\Bi(Fy, ((f/G)^1), Ex)}"}
\Ar{1c}{2c}{"{\Bi(y, f, x)}"}
\Ar{1}{1a}{equal}
\Ar{2}{2a}{equal}
\Ar{1a}{1b}{equal}
\Ar{2a}{2b}{equal}
\Ar{1b}{1c}{equal}
\Ar{2b}{2c}{equal}
\end{tikzcd}.
$$
Thus we can define a natural isomorphism $\et_M \colon M \to ((M/G)\# G)'$ as the identity.

Let $f \colon N \to N'$ be a morphism in 
$\Ggr(\biMod{S/G}{R/G})$, and $Ex \in (R/G)_0,\, Fy \in (S/G)_0$.
Then we have a commutative diagram
$$
\begin{tikzcd}[column sep=130pt]
\Bi(Fy,((N\# G)'/G), Ex) & \Bi(Fy,((N'\# G)'/G), Ex)\\
\DS_{a \in G}\Bi(y, (N\# G)', ax) &
\DS_{a \in G}\Bi(y, (N'\# G)', ax)\\
\DS_{a \in G} \Bi((Fy)^{(1)}, (N \# G), {(Eax)^{(1)}}) &
\DS_{a \in G} \Bi({(Fy)}^{(1)}, (N' \# G), {(Eax)^{(1)}})\\[15pt]
\DS_{a \in G} \Bi((Fy)^{(1)}, (N \# G), {(Ex)^{(a)}}) &
\DS_{a \in G} \Bi({(Fy)}^{(1)}, (N' \# G), {(Ex)^{(a)}})\\
\DS_{a \in G} \Bi(Fy, N, Ex)^a & \DS_{a \in G} \Bi(Fy, N', Ex)^a\\
\Bi(Fy,N,Ex) & \Bi(Fy,N',Ex)
\Ar{1-1}{2-1}{equal}
\Ar{1-2}{2-2}{equal}
\Ar{2-1}{3-1}{equal}
\Ar{2-2}{3-2}{equal}
\Ar{3-1}{4-1}{"\wr"', "{\DS_{a \in G} \Bi((Fy)^{(1)}, (N \# G), \ovl{\ph^E_a x})}"}
\Ar{3-2}{4-2}{"{\DS_{a \in G} \Bi({(Fy')}^{(1)}, (N \# G), \ovl{\ph^E_a x'})}"', "\wr"}
\Ar{4-1}{5-1}{equal}
\Ar{4-2}{5-2}{equal}
\Ar{5-1}{6-1}{equal}
\Ar{5-2}{6-2}{equal}
\Ar{1-1}{1-2}{"{\Bi(Fy,((f\# G)'/G), Ex)}"}
\Ar{2-1}{2-2}{"{\DS_{a \in G}\Bi(y, (f\# G)', ax)}"}
\Ar{3-1}{3-2}{"{\DS_{a \in G} {}_{(Fy)^{(1)}}(f \# G)_{(Eax)^{(1)}}}" pos=0.45}
\Ar{4-1}{4-2}{"{\DS_{a \in G} \Bi(Fy^{(1)}, (f \# G), {(Ex)^{(a)}})}"}
\Ar{5-1}{5-2}{"{\DS_{a \in G} \Bi(Fy, f, Ex)^a }"}
\Ar{6-1}{6-2}{"{\Bi(Fy,f,Ex)}"}
\end{tikzcd}
$$
Therefore, we can define an isomorphism
$\ep_N \colon (N\# G)'/G \to N$ by
$$
\Bi(Fy,(\ep_N),Ex):=
\Ds_{a \in G} \Bi((Fy)^{(1)},(N\# G),\ovl{\ph^E_{a} x})
$$
for all $x \in R_0,\, y \in S_0$
that is natural in $N$ as shown above.
Thus $\ep:= (\ep_N)_{N \in \Ggr(\Bi(S/G, \Mod, R/G))}$
is a natural isomorphism
$$
(?/G) \circ (?\# G)' \To \id_{\Ggr(\biMod{S/G}{R/G})}.
$$
As a consequence, $?/G$ is an equivalence with a quasi-inverse
$(?\# G)'$.

By a general theory, $?/G$ becomes a left adjoint to $(?\# G)'$,
and we have an adjunction natural isomorphism $\th = (\th_{M,N})_{M,N}$, where
$$
\th_{M,N} \colon \Ggr(\biMod{S/G}{R/G})(M/G, N) \to
\Ginv(\biMod{S}{R})(M, (N\# G)')
$$
is defined by sending $f \colon M/G \to N$ to
$(f\# G)' \circ \et_M = (f\# G)'$
for all $M \in \Ggr(\biMod{S/G}{R/G})_0$ and
$N \in \Ginv(\biMod{S}{R})_0$.

We here note that the unit and the counit defined by
this $\th$ coincide with $\et$ and $\ep$ defined above.
The unit is given by $\th(\id_{M/G}) = \id_M = \et_M$, as
desired.
For the counit, it is enough to show that
$\th\inv(\id_{(N\# G)'}) = \ep_N$, or equivalently,
$$
(\ep_N \# G)' = \id_{(N\# G)'}
$$
because $\th(\ep_N) = (\ep_N \# G)'$.
This follows from the following commutative diagram:
$$
\begin{tikzcd}[column sep=80pt]
\Nname{1} \Bi(y, ((((N\# G)'/G) \# G)'), x) & \Nname{2} \Bi(y, ((N\# G)'), x)\\
\Nname{1a}\Bi((Fy)^{(1)}, (((N\# G)'/G) \# G), {(Ex)^{(1)}}) & \Nname{2a}\Bi((Fy)^{(1)}, (N\# G), {(Ex)^{(1)}}) \\
\Nname{1b}\Bi(Fy, (((N\# G)'/G))^1, Ex) &\Nname{2b}\Bi(Fy, N^1, Ex)\\
\Nname{1c}\Bi(y, ((N\# G)'), 1x) \\
\Nname{1d}\Bi((Fy)^{(1)}, (N\# G), {(Ex)^{(1)}}) & \Nname{2d}\Bi((Fy)^{(1)}, (N\# G), {(Ex)^{(1)}})\\
\Nname{1e}\Bi(Fy, N^1, Ex) & \Nname{2e}\Bi(Fy, N^1, Ex) \\
\Ar{1}{2}{"{\Bi(y,(\ep_N\# G)',x)}"}
\Ar{1a}{2a}{"{\Bi((Fy)^{(1)}, (\ep_N\# G), {(Ex)^{(1)}})}"}
\Ar{1b}{2b}{"{\Bi(Fy, (\ep_N^1), Ex)}"}
\Ar{1d}{2d}{"{\Bi((Fy)^{(1)}, (N\# G), \ovl{\ph^E_1 x})}"}
\Ar{1e}{2e}{"{\Bi(Fy, N, {(\id_{Ex})}) = \id_{(\bi{N^1}{Fy}{Ex})}}"}
\Ar{1}{1a}{equal}
\Ar{2}{2a}{equal}
\Ar{1a}{1b}{equal}
\Ar{2a}{2b}{equal}
\Ar{1b}{1c}{equal}
\Ar{2b}{2d}{equal}
\Ar{1c}{1d}{equal}
\Ar{1d}{1e}{equal}
\Ar{2d}{2e}{equal}
\end{tikzcd}.
$$

(2)
We set $(P, \ph^P) \colon A\# G \to (A\# G)/G$ and
$(Q, \ph^Q) \colon B\# G \to (B\# G)/G$ to be the canonical $G$-coverings.
Consider the following diagram that illustrates the setting
of the second assertion:
$$
\begin{tikzcd}
\Nname{A}\Ggr(\biMod{B}{A}) &  & \Nname{A'}\Ginv(\biMod{B\# G}{A\# G})\\
\Nname{A1}\Ggr(\biMod{B}{A}) &  & \Nname{A1'}\Ginv(\biMod{B\# G}{A\# G})
\Ar{A}{A'}{"?\# G"}
\Ar{A}{A1}{equal}
\Ar{A'}{A1}{"(?/G)'" '}
\Ar{A'}{A1'}{equal}
\Ar{A1}{A1'}{"?\# G"}
\ar["\et'", Rightarrow, to path={([xshift=10, yshift=-10]R.south)--([xshift=30, yshift=-10]R.south)\tikztonodes}]
\ar["\ep'", Rightarrow, to path={([xshift=-30, yshift=-10]R'.south)--([xshift=-10, yshift=-10]R'.south)\tikztonodes}]
\end{tikzcd}.
$$
To show that $?\# G$ is an equivalence with a quasi-inverse
$(?/ G)'$ it is enough to construct natural isomorphisms
$$
\begin{aligned}
\et' &\colon  \id_{\Ggr(\biMod{B}{A})} \To  (?/G)' \circ (?\# G), \text{ and}\\
\ep' &\colon (?\# G) \circ (?/G)' \To \id_{\Ginv(\biMod{B\# G}{A\# G})}.
\end{aligned}
$$
Let $N \in \Ggr(\biMod{B}{A})_0$, and $x \in A_0,\, y \in B_0$.
Then
$$
\begin{aligned}
\Bi(y,((N\# G)/G)',x) &= \Bi(Q(y^{(1)}),((N\# G)/G),{P(x^{(1)})})
= \Ds_{a \in G} \Bi(y^{(1)},N\# G, {a\cdot x^{(1)}})\\
&= \Ds_{a \in G} \Bi(y^{(1)},N\# G, {x^{(a)}})
= \Ds_{a \in G} \Bi(y,N,x)^a
= \Bi(y,N,x).
\end{aligned}
$$
Hence $\id_{\Ggr(\biMod{B}{A})} = (?/G)' \circ (?\# G)$, and
we can define $\et'$ as the identity natural transformation.

Let $M = (M, \ph) \in \Ginv(\biMod{B\# G}{A\# G})_0$, and
$x^{(a)} \in (A\# G)_0$, $y^{(b)} \in (B\# G)_0$.
Then we have an isomorphism
$$
\begin{aligned}
\Bi(y^{(b)}, ((M/G)'\# G), {x^{(a)}})
&= (\Bi(y, ((M/G)'), x))^{b\inv a}
= (\Bi(Q(y^{(1)}), (M/G), {P(x^{(1)})}))^{b\inv a}\\
&= \Bi(y^{(1)}, M,{(b\inv a)x^{(1)}})
\ya{\ph_b} \Bi(by^{(1)}, M,{ax^{(1)}})
= \Bi(y^{(b)}, M,{x^{(a)}}).
\end{aligned}
$$
Hence we can define a natural isomorphism $\ep'$ by setting
$$
\Bi(y^{(b)},(\ep'_M),{x^{(a)}}):=\ph_b \colon \Bi(y^{(b)}, ((M/G)'\# G), {x^{(a)}}) \to \Bi(y^{(b)}, M,{x^{(a)}}).
$$
To verify the naturality of $\ep'$,
let
$$
\begin{aligned}
\smor{a}{f}{a'} \in \Bi(\sobj{x}{a},(A\# G), \sobj{x'}{a'}) &= \{a\} \times \Bi(x,A,x')^{a\inv a'} \times \{a'\}, \\
\smor{b'}{g}{b} \in \Bi(\sobj{y'}{b'},(B\# G), \sobj{y}{b}) &= \{b'\} \times \Bi(y',B,y)^{{b'}\inv b}\times \{b\}.
\end{aligned}
$$ 
Then the naturality of $\ep'$ is equivalent to the commutativity of the diagram
\begin{equation}
\label{eq:nat-ep'}
\begin{tikzcd}[column sep=80pt]
\Nname{0}(\Bi(Q(\sobj y1),(M/G),{P(\sobj x1)}))^{b\inv a} &
\Nname{0'}(\Bi(Q(\sobj {y'}1),(M/G),{P(\sobj {x'}1)}))^{{b'}\inv a'}\\
\Nname{1}\Bi(\sobj{y}{1}, M, {(b\inv a) \sobj{x}{1}}) &
\Nname{1'}\Bi(\sobj{y'}{1}, M, {({b'}\inv a') \sobj{x'}{1}})\\
\Nname{a}\Bi(\sobj{y}{b}, M , \sobj{x}{a}) &
\Nname{a'}\Bi(\sobj{y'}{b'}, M , \sobj{x'}{a'})
\Ar{0}{0'}{"{\Bi(\om_B(g),(M/G),{\om_A(f)})}"}
\Ar{0}{1}{equal}
\Ar{0'}{1'}{equal}
\Ar{1}{1'}{"F", dashed}
\Ar{a}{a'}{"{\Bi(\smor{b'}{g}{b},M,\smor{a}{f}{a'})}" '}
\Ar{1}{a}{"{\Bi(\sobj{y}{1},(\ph_b),\sobj{x}{b\inv a})}" '}
\Ar{1'}{a'}{"{\Bi(\sobj{y'}{1},(\ph_{b'}),\sobj{x'}{{b'}\inv a'})}"}
\end{tikzcd}.
\end{equation}
We first show that the map $F$ above is given as follows:
$$
F = \Bi({b'}\inv (\smor{b'}{g}{b}),M,{{b'}\inv (\smor{a}{f}{a'})}) \circ \Bi(\sobj{y}{1}, (\ph_{{b'}\inv b}),{\sobj{x}{b\inv a}}).
$$
To show this, let $m \in \Bi(y^{(1)}, M, {x^{(b\inv a)}})$,
the precise form of which is
$(\de_{b\inv a, d}\, m)_{d \in G}$ as an element of
$\Bi(Q(\sobj y1),(M/G),{P(\sobj x1)})$.
Note that the precise forms of $\om_A(f)$ and $\om_B(g)$ is given
as follows (by a careful reading of \cite[Definition 5.8.8]{Asa-book}):
$$
\om_A(f) = (\de_{a\inv a', e} \smor{1}{f}{a\inv a'})_{e \in G},
\quad
\om_B(g) = (\de_{{b'}\inv b, c} \smor{1}{g}{{b'}\inv b})_{c \in G}.
$$
Then by the formula \eqref{eq:orbit-bimod-str}, we have
$$
\begin{aligned}
F(m) &= \om_B(g)\cdot m \cdot \om_A(f)\\
&= \left(\sum_{cde = u}\de_{{b'}\inv b, c}\smor{1}{g}{{b'}\inv b}
\cdot \ph_c(\de_{b\inv a, d}\, m) \cdot
cd\de_{a\inv a',e}\smor{1}{f}{a\inv a'}
\right)_{u \in G}\\
&= (\de_{{b'}\inv a', u} \smor{1}{g}{{b'}\inv b} \cdot 
\ph_{{b'}\inv b}(m) \cdot
({b'}\inv a) \,\smor{1}{f}{a\inv a'})_{u\in G}\\
&= {b'}\inv (\smor{b'}{g}{b})\cdot \ph_{{b'}\inv b}(m) \cdot {b'}\inv(\smor{a}{f}{a'})\\
&= (\Bi({b'}\inv (\smor{b'}{g}{b}),M,{{b'}\inv (\smor{a}{f}{a'})}) \circ \Bi(\sobj{y}{1}, (\ph_{{b'}\inv b}),{\sobj{x}{b\inv a}}))(m).
\end{aligned}
$$
For simplicity, we use the shorter forms.
Namely, $\Biup((d),f,{(c)})$ is denoted just by $f$ for all $c,d \in G$; the same for $g$; and $\Bi(v,(\ph_c),u)$ simply by $\ph_c$ for all $u \in (A\# G)_0,\, v \in (B\# G)_0,$ and $c \in G$.
Then the commutativity of the diagram \eqref{eq:nat-ep'}
follows from the naturality of $\ph_b$ 
(expressed by dashed arrows in the diagram below) for
all $b \in G$ and the fact that $\ph$ is the
$G$-invariant structure of $M$ (Definition \ref{dfn:G-inv} (1))
as the following commutative diagram shows:
$$
\begin{tikzcd}[column sep=60pt]
\Nname{1a}\Bi(y^{(1)}, M, {x^{(b\inv a)}}) &
\Nname{ba}\Bi(y^{({b'}\inv b)}, M, {x^{({b'}\inv a)}}) &
\Nname{ba'}\Bi(y^{({b'}\inv b)}, M, {{x'}^{({b'}\inv a')}})\\
&
\Nname{2-1a}\Bi({y'}^{(1)}, M, {x^{({b'}\inv a)}}) &
\Nname{1a'}   \Bi({y'}^{(1)}, M, {{x'}^{({b'}\inv a')}})\\
\Nname{3-ba}\Bi(y^{(b)}, M , x^{(a)}) &&
\Nname{b'a'}\Bi({y'}^{(b')}, M , {x'}^{(a')})
\Ar{1a}{ba}{"\ph_{{b'}\inv b}"}
\Ar{ba}{ba'}{"{\Bi(y^{({b'}\inv b)},M,{{b'}\inv f})}"}
\Ar{2-1a}{1a'}{"{\Bi({y'}^{(1)},M, {{b'}\inv f})}" '}
\Ar{3-ba}{b'a'}{"{\Bi(g,M,f)}" ', dashed}
\Ar{1a}{3-ba}{"\ph_b" '}
\Ar{ba}{2-1a}{"{\Bi({b'}\inv g,M,{x^{({b'}\inv a)}})}" '}
\Ar{ba'}{1a'}{"{\Bi({b'}\inv g, M, {{x'}^{({b'}\inv a')}})}"}
\Ar{1a'}{b'a'}{"\ph_{b'}", dashed}
\Ar{ba}{3-ba}{"\ph_{b'}", dashed}
\Ar{ba}{1a'}{"{\Bi({b'}\inv g,M,{b'}\inv f)}" ' pos=0.4, dashed}
\end{tikzcd}.
$$

As a consequence, $?\# G$ is an equivalence with a quasi-inverse $(?/ G)'$.

By a general theory, $?\# G$ becomes a left adjoint to $(?/ G)'$,
and we have an adjunction natural isomorphism $\th' = (\th'_{N,M})_{N,M}$, where
$$
\th'_{N,M} \colon \Ginv(\biMod{B\# G}{A\# G})(N\# G, M) \to
\Ggr(\biMod{B}{A})(N, (M/G)')
$$
is defined by 
$$
\th'_{N,M}(f):= (f/ G)' \circ \et'_N = (f/ G)'
$$
for all $f \colon N\# G \to M$,
all $N \in \Ggr(\biMod{B}{A})_0$ and all
$M \in \Ginv(\biMod{B\# G}{A\# G})_0$.

Note that the unit and the counit defined by
this $\th'$ coincide with $\et'$ and $\ep'$ defined above.
The unit is given by $\th'(\id_{N\# G}) = \et'_N (= \id_N)$, as
desired.
For the counit, it is enough to show that
$\ep'_M = (\th')\inv(\id_{(M/ G)'})$, or equivalently,
$$
(\ep'_M / G)' = \id_{(M/ G)'}
$$
because $\th'(\ep'_M) = (\ep'_M / G)'$.
This follows from the following commutative diagram:
$$
\begin{tikzcd}[column sep=80pt]
\Nname{1} \Bi(y, (((M/ G)'\# G) / G)', x) & \Nname{2} \Bi(y, (M/ G)', x)\\
\Nname{1a}\Bi(Q(\sobj{y}{1}), ((M/ G)'\# G) / G), {P(\sobj{x}{1})}) & \Nname{2a}\Bi(Q(\sobj{y}{1}), (M/ G), {P(\sobj{x}{1})}) \\
\Nname{1b}\Ds_{a\in G}\Bi(y^{(1)}, ((M/ G)'\# G)), {ax^{(1)}}) &\Nname{2b}\Ds_{a\in G}\Bi(y^{(1)}, M, {ax^{(1)}})\\
\Nname{1c}\Ds_{a\in G}\Bi(y^{(1)}, ((M/G)'\# G), x^{(a)})& \Nname{2c}\Ds_{a\in G}\Bi(y^{(1)}, M, x^{(a)}) \\
\Nname{1d}\Ds_{a\in G}(\Bi(y, (M/G)', x))^a  &
\Nname{2d}\Bi(y,(M/G)',x)\\
\Ar{1}{2}{"{\Bi(y,(\ep'_M/ G)',x)}"}
\Ar{1a}{2a}{"{\Bi(Q(\sobj{y}{1}), (\ep'_M/ G), {P(\sobj{x}{1})})}"}
\Ar{1b}{2b}{"{\Ds_{a\in G}\Bi(y^{(1)}, (\ep_M'), {ax^{(1)}})}"}
\Ar{1c}{2c}{"{\Ds_{a\in G}\Bi(y^{(1)}, (\ph_1), {x^{(a)}})}"}
\Ar{1d}{2d}{"{\id_{\Bi(y,(M/G)',x)}}"}
\Ar{1}{1a}{equal}
\Ar{2}{2a}{equal}
\Ar{1a}{1b}{equal}
\Ar{2a}{2b}{equal}
\Ar{1b}{1c}{equal}
\Ar{2b}{2c}{equal}
\Ar{2c}{2d}{equal}
\Ar{1c}{1d}{equal}
\end{tikzcd}. \qedhere
$$
\end{proof}

To prove Proposition \ref{prp:graded-prj-bimod}, we need the following.

\begin{lem}
\label{lem:G-gr-free}
Let $M \in \Ggr(\biMod{B}{A})$.
Then there exists a small set $I$,
$(x_i, y_i)_{i \in I} \in (A_0\times B_0)^I$ and 
$(a_i, b_i)_{i \in I} \in (G\times G)^I$
such that there exists an epimorphism
$$
\textstyle \Ds_{i\in I} \ro_{b_i}(B_{y_i}) \ox_\k \la_{a_i}({}_{x_i}A) \to M
$$
in $\Ggr(\biMod{B}{A})$.
\end{lem}

\begin{proof}
The functor $?\# G \colon \Ggr(\biMod{B}{A}) \to \Ginv(\biMod{(B\# G)}{(A\# G)})$
is an equivalence with a quasi-inverse given by the composite
$$
(?/G)'\colon
\Ginv(\biMod{(B\# G)}{(A\# G)}) \ya{?/G} \Ggr(\biMod{(B\# G)/G}{(A\# G)/G})
\ya{\ovl{\om}} \Ggr(\biMod{B}{A}).
$$
Apply Lemma \ref{lem:G-inv-generate} to the situation that
$R:= A\# G$, $S:= B\# G$, and $N:=M\#G \in \Ginv(\biMod{B\# G}{A\# G})_0$.
Here we set $(P, \ph^P) \colon A\# G \to (A\# G)/G$ and
$(Q, \ph^Q) \colon B\# G \to (B\# G)/G$ to be the canonical $G$-coverings.
Then there exist a small set $I$ and
a family $(x^{(a_i)}_i, y^{(b_i\inv)}_i)_{i\in I} \in ((A\# G)_0 \times (B\# G)_0)^I$ such that
there exists an epimorphism
\[
F\colon\textstyle\Ds_{i \in I} \left(\Ds_{c \in G} (B\# G)_{cy^{(b_i\inv)}_i}\ox_\k {}_{cx^{(a_i)}_i}(A\# G)\right) \to N
\]
in the category $\Ginv(\biMod{(B\# G)}{(A\# G)})$.
Applying the equivalence $(?/G)'$ to the above epimorphism $F$, we have an epimorphism
\[
(F/G)'\colon\textstyle
\left(\Ds_{i \in I} \left(\Ds_{c \in G} (B\# G)_{cy^{(b_i\inv)}_i}\ox_\k {}_{cx^{(a_i)}_i}(A\# G)\right)/G\right)' \to (N/G)'\iso M.
\]
Hence it is enough to show that
\begin{equation}
\label{eq:formula/G'}
\textstyle\left(\left(\Ds_{c \in G}(B\# G)_{cy^{(b_i\inv)}_i}\ox_\k {}_{cx^{(a_i)}_i}(A\# G)\right)/G\right)'\iso
\ro_{b_i}(B_{y_i}) \ox_\k \la_{a_i}({}_{x_i}A).
\end{equation}
By Proposition \ref{prp:orbitr-bi}, we have 
\[
\textstyle\left(\Ds_{c \in G} (B\# G)_{cy^{(b_i\inv)}_i}\ox_\k {}_{cx^{(a_i)}_i}(A\# G)\right)/G\iso
(B\# G)/G_{Q(y^{(b_i\inv)}_i)}\ox_\k{}_{P(x^{(a_i)}_i)}(A\# G)/G.
\]
Therefore, it is enough to show the following. 
\begin{clm}
$
\ovl{\om}((B\# G)/G_{Q(y^{(b_i\inv)}_i)}\ox_\k{}_{P(x^{(a_i)}_i)}(A\# G)/G)
\iso \ro_{b_i}(B_{y_i}) \ox_\k \la_{a_i}({}_{x_i}A).
$
\end{clm}
For any $y'\in B_0$ and $x'\in A_0$,
\[
\begin{aligned}
_{y'}(\mathrm{LHS})_{x'}&=_{Q(y'^{(1)})}\!\!(B\# G)/G_{Q(y_i^{(b_i\inv)})}\ox_\k {}_{P(x_i^{(a_i)})}(A\# G)/G_{P(x'^{(1)})}\\
&=\textstyle(\Ds_{d \in G} {}_{(y')^{(1)}}(B\# G)_{dy^{(b_i\inv)}_i})\ox_\k (\Ds_{c\in G}{}_{x^{(a_i)}_i}(A\# G)_{c{((x')^{(1)})}})\\
&= \textstyle(\Ds_{d \in G} {}_{y'}B^{db_i\inv}_{y_i})\ox_\k (\Ds_{c\in G}{}_{x}A^{a_i\inv c}_{x'}).
\end{aligned}
\]
On the other hand,
\[
\begin{aligned}
_{y'}(\mathrm{RHS})_{x'}
&=\textstyle\ro_{b_i}(_{y'}B_{y_i}) \ox_\k \la_{a_i}({}_{x_i}A_{x'})=\ro_{b_i}(\Ds_{d \in G}{}_{y'}B^d_{y_i}) \ox_\k \la_{a_i}(\Ds_{c\in G} {}_{x_i}A^c_{x'})\\
&=\textstyle(\Ds_{d \in G} {}_{y'}B^{db_i\inv}_{y_i})\ox_\k (\Ds_{c\in G}{}_{x}A^{a_i\inv c}_{x'}).
\end{aligned}
\]
Hence $_{y'}(\mathrm{RHS})_{x'}=_{y'}(\mathrm{LHS})_{x'}$.
As a consequence, there exists an epimorphism
$$
\textstyle \Ds_{i\in I} \ro_{b_i}(B_{y_i}) \ox_\k \la_{a_i}({}_{x_i}A) \to M
$$
in $\Ggr(\biMod{B}{A})$.
\end{proof}

Here we give a statement that is slightly finer than
Proposition \ref{prp:graded-prj-bimod}.

\begin{prp}
\label{prp:graded-prj-bimod-eq}
Let $\Fgt \colon \Ggr(\biMod{B}{A}) \to \biMod{B}{A}$ be
the forgetful functor, and
$P \in \Ggr(\biMod{B}{A})_0$.
Then the following are equivalent.
\begin{enumerate}
\item 
$P$ is projective in $\Ggr(\biMod{B}{A})$
\item
$\Fgt(P)$ is projective in $\biMod{B}{A}$.
\item 
$P$ is a direct summand of an object of the form
\begin{equation}
\label{eq:graded-prj-bimod}
\Ds_{i\in I} \ro_{b_i}(B_{y_i}) \ox_\k \la_{a_i}({}_{x_i}A)
\end{equation}
in $\Ggr(\biMod{B}{A})$ for some family
$(x_i, y_i)_{i \in I} \in (A_0 \times B_0)^I$
and $(a_i, b_i)_{i \in I} \in (G \times G)^I$
with $I$ a small set.
\end{enumerate}
\end{prp}

\begin{proof}
(1)\implies (3).
We apply Lemma \ref{lem:G-gr-free} to have an epimorphism
$$
\textstyle f \colon \Ds_{i\in I} \ro_{b_i}(B_{y_i}) \ox_\k \la_{a_i}({}_{x_i}A) \to P
$$
in $\Ggr(\biMod{B}{A})$
for some family
$(x_i, y_i)_{i \in I} \in (A_0 \times B_0)^I$
and $(a_i, b_i)_{i \in I} \in (G \times G)^I$
with $I$ a small set.
Hence if (1) holds, then $f$ splits and (3) holds.

(3)\implies (2).
Assume (3).
Then $\Fgt(P)$ turns out to be a direct summand
of $\Fgt(\Ds_{i\in I} \ro_{b_i}(B_{y_i}) \ox_\k \la_{a_i}({}_{x_i}A)) \iso
\Ds_{i\in I} B_{y_i} \ox_\k {}_{x_i}A$, which is projective in
$\biMod{B}{A}$.
Thus $\Fgt(P)$ is projective in $\biMod{B}{A}$.

(2)\implies (1).
This part is also proved by the same way as in the case of $G$-graded
algebras as follows.

For each $M \in \Ggr(\biMod{B}{A})$, we write $M = (M, d^M)$ to stress that
$M$ is a $G$-graded $B$-$A$-bimodule, where
$d^M$ is the $G$-grading of $M$.
Then $\Fgt(M, d^M) = M$.

Assume that $\Fgt(P, d^P) = P$ is projective in $\biMod{B}{A}$.
To show that $(P, d^P)$ is projective in $\Ggr(\biMod{B}{A})$, consider
the diagram
$$
\begin{tikzcd}
& (P, d^P)\\
(N, d^N) & (M, d^M)
\Ar{1-2}{2-2}{"f"}
\Ar{2-1}{2-2}{"g"}
\end{tikzcd}
$$
in $\Ggr(\biMod{B}{A})$, where $g$ is an epimorphism.
By the assumption, there exists a homomorphism $h \colon P \to N$
in $\biMod{B}{A}$ such that $f = gh$.
We define an $h' \colon (P, d^P) \to (N, d^N)$ in $\Ggr(\biMod{B}{A})$
such that $f = gh'$.

Let $x \in A_0,\, y \in B_0$ and
$p \in \bi{P}{y}{x} = \Ds_{a\in G}\bi{P}{y}{x}^a$.
Then $p = \sum_{a\in G} p^a$ for some $p^a \in \bi{P}{y}{x}^a\ (a \in G)$.
For each $a \in G$,
since $\bi{h}{y}{x}(p^a) \in \bi{N}{y}{x} = \Ds_{b \in G}\bi{N}{y}{x}^b$,
we have
$\bi{h}{y}{x}(p^a) = \sum_{b \in G} n^b$ for some unique
$(n^b)_{b\in G} \in \bigoplus_{b \in G}\bi{N}{y}{x}^b$.
Hence $\bi{f}{y}{x}(p^a) = \sum_{b \in G} \bi{g}{y}{x}(n^b)$.
Since $\bi{f}{y}{x}(p^a) \in \bi{M}{y}{x}^a$ and
$\bi{g}{y}{x}(n^b) \in \bi{M}{y}{x}^b$, we have
$$
\bi{g}{y}{x}(n^b) =
\begin{cases}
\bi{g}{y}{x}(n^a) & (b = a)\\
0 & (b \ne a)
\end{cases},
$$
and hence we have
\begin{equation}
\label{eq:fm=gn}
\bi{f}{y}{x}(p^a) = \bi{g}{y}{x}(n^a).
\end{equation}
Note that this $n^a \in \bi{N}{y}{x}^a$ is uniquely determined by $p^a$.
Therefore, we can define $\bi{h'}{y}{x} \colon \bi{P}{y}{x} \to \bi{N}{y}{x}$ by
$\bi{h'}{y}{x}(p):= \sum_{a\in G}n^a$.
Then we have
\begin{equation}
\label{eq:h'm=n}
\bi{h'}{y}{x}(p^a) = n^a,
\end{equation}
which means that $h'$ is degree-preserving.
Note by definition that the following holds
\begin{equation}
\label{eq:dfn-h'}
\bi{h'}{y}{x}(p):= \sum_{a\in G}n^a \quad \Longleftrightarrow \quad
\bi{f}{y}{x}(p^a) = \bi{g}{y}{x}(n^a) \ \text{for all $a \in G$.}
\end{equation}
The $\k$-linearity of $h'$ is obvious.
To show that $h'$ is a morphism in the category $\Ggr(\biMod{B}{A})$,
it remains to show the following:
$$
\bi{h'}{y'}{x'}(v \cdot p \cdot u) = v \cdot \bi{h'}{y}{x}(p) \cdot u
$$
for all $x, x' \in A_0$, $y, y' \in B_0$, and
$u \in \bi{A}{x}{x'},\, v \in \bi{B}{y'}{y},\, p \in \bi{P}{y}{x}$.

Since $f$ and $g$ are morphisms in $\Ggr(\biMod{B}{A})$, we have
$$
\bi{f}{y'}{x'}(v^c p^a u^b) = v^c \bi{f}{y}{x}(p^a) u^b,\quad
\bi{g}{y'}{x'}(v^c n^a u^b) = v^c \bi{g}{y}{x}(n^a) u^b
$$
for all
$n \in \bi{N}{y}{x}$, $a, b, c \in G$.
Now let $\bi{h'}{y}{x}(p) = \sum_{a \in G}n^a$, namely
$\bi{f}{y}{x}(p^a) = \bi{g}{y}{x}(n^a)$ for all $a \in G$.
Then $\bi{f}{y'}{x'}(v^c p^a u^b) = v^c \bi{f}{y}{x}(p^a) u^b =
v^c \bi{g}{y}{x}(n^a) u^b = \bi{g}{y'}{x'}(v^c n^a u^b)$
for all $a, b, c \in G$.
Hence
$$
\bi{f}{y'}{x'}\left(\sum_{cab = d} v^cp^a u^b\right) = \bi{g}{y'}{x'}\left(\sum_{cab = d} v^c n^a u^b\right),
$$
and therefore, since $vpu = \sum_{d\in G} \left(\sum_{cab = d} v^c p^a u^b \right)$,
we have the following by \eqref{eq:dfn-h'}:
$$
\bi{h'}{y'}{x'}(v p u) = \sum_{d\in G} \left(\sum_{cab = d} v^c n^a u^b \right)
= \left(\sum_{c \in G}v^c\right)\left(\sum_{a\in G}n^a\right)\left(\sum_{b \in G}u^b\right) = v \cdot \bi{h'}{y}{x}(p) \cdot u.
$$
By \eqref{eq:fm=gn} and \eqref{eq:h'm=n},
we have $f = gh'$ in $\Ggr(\biMod{B}{A})$.
Thus $(P, d^P)$ is projective in $\Ggr(\biMod{B}{A})$.
\end{proof}

\section{F.g.\ projective \texorpdfstring{$G$}{G}-invariant
(resp.\ \texorpdfstring{$G$}{G}-graded) bimodules}

In this section, we give explicit forms of
f.g.\ projective $G$-invariant (resp.\ $G$-graded) bimodules
over locally bounded categories.
Therefore, throughout this section,
{\em we assume that $\k$ is a field}.

\begin{rmk}
\label{rmk:skeleton-Ggr}
Let $A'$ be a skeleton of a $G$-graded category $A$.
Then $A'$ is also a $G$-graded category, and the inclusion
functor $\si_A \colon A' \to A$ becomes
a strictly degree-preserving equivalence.
Consider an $A'$-$A$-bimodule ${}_{\si_A} A$.
Then the tensor product functor
$\blank \ox_{A'}({}_{\si_A} A) \colon \Mod_{A'} \to \Mod_A$
turns out to be an equivalence.
Let $B'$ be a skeleton of a $G$-graded category $B$ with
the inclusion functor $\si_B \colon B' \to B$.
Consider a $B$-$B'$-bimodule $B_{\si_B}$.
Then the equivalence $\biMod{B'}{A'} \to \biMod{B}{A}$ is defined by
sending $M$ to $(B_{\si_B}) \ox_{B'} M \ox_{A'}({}_{\si_A}A)$.
This also yields an equivalence $\Ggr(\biMod{B'}{A'}) \to \Ggr(\biMod{B}{A})$.
\end{rmk}

\begin{rmk}
\label{rmk:classical-locbdd}
Let $(R, X)$ be a locally finite-dimensional category, where
the $G$-action $X$ is free and {\em locally bounded}, in the sense
that for each pair $(x,y)$ in $R$, the set $\{a \in G \mid \bi{R}{y}{ax} \iso \bi{R}{a\inv y}{x} \ne 0\}$ is finite. Then
\begin{enumerate}
\item 
By \cite[Proposition 3.1]{Gab},
the classical orbit category $R \oorbit G$ of $R$ by $G$ is defined,
and $R \oorbit G$ is shown to be a locally finite-dimensional category,
where
$(R\oorbit G)_0:= \{Gx \mid x \in R_0\},\ (Gx:= \{ax \mid a \in G\})$.
Note that $R \oorbit G$ is isomorphic to any skeleton $A$ of $R/G$ as a
$G$-graded category by \cite[Remark 2.2(1) and Proposition 2.11]{Asa11}
with a strictly degree-preserving isomorphism $A \to R\oorbit G$ that sends each $x \in A_0$ to $Gx$, by which we identify these $G$-graded categories.
Then we see that $A$ is a locally finite-dimensional category.
\item
If $R$ above is a locally bounded category,
then so is $A$.
Indeed, let $y \in A_0$.
Then for any $x \in A_0$, we have the following equivalences
$$
\begin{aligned}
\bi{A}{y}{x} \ne 0
&\iff \bi{R/G}{y}{x} \ne 0
\iff \Ds_{a\in G}\bi{R}{y}{ax} \ne 0\\
&\iff \exists a\in G, \bi{R}{y}{ax} \ne 0
\iff U(y) \cap Gx \ne \emptyset,
\end{aligned}
$$
where we set $U(y):= \{z \in R_0 \mid \bi{R}{y}{z} \ne 0\}$, which is a finite set because
$R$ is a locally bounded category.
Hence $\{x \in A_0 \mid \bi{A}{y}{x} \ne 0\}$ is a finite set.
Similarly, $\{x \in A_0 \mid \bi{A}{x}{y} \ne 0\}$ is
shown to be a finite set by using the fact that
$\bi{R}{x}{ay} \iso \bi{R}{a\inv x}{y}$ for all $a \in G$. Then by definition \ref{dfn:loc-fd-cat}, $A$ is a locally bounded category.
\end{enumerate}
\end{rmk}

\begin{lem}
\label{lem:ind-proj-G-inv}
Let $(R, X)$ and $(S, Y)$ be locally bounded categories
with free $G$-actions, and $(x,y) \in R_0 \times S_0$.
Then the f.g.\ projective $G$-invariant $S$-$R$-bimodule
$$
(M, \ph):= \Ds_{a\in G} S_{ay}\ox_\k {}_{ax}R
$$
with the canonical $G$-invariant structure has a local endomorphism algebra, and hence is an indecomposable projective object in $\Ginv(\biMod{S}{R})$.
\end{lem}

\begin{proof}
By \eqref{eq:G-inv-slash-G}, we have $M/G \iso (S/G)_y \ox_\k {}_x(R/G)$.
Choose skeletons $A$ and $B$ of $R/G$ and $S/G$, respectively
in such a way that $y \in B$ and $x \in A$.
As stated in Remark \ref{rmk:skeleton-Ggr},
both of the inclusion functors $A \to R/G$
and $B \to S/G$ become strictly degree-preserving equivalences.
By Remark \ref{rmk:classical-locbdd}(1)
both $A$ and $B$ are locally finite-dimensional categories,
and then as is easily seen, so is $B\ox_\k {A}\op$. 
Therefore, as an object in the category $\biMod{B}{A}$,
$B_y \ox_\k {}_x A = \Fgt(B_y \ox_\k {}_x A)$
has a local endomorphism algebra.
Indeed, since $B_y \ox_\k {}_x A = \bi{(B\ox_\k A\op)}{}{(y,x)}$,
its endomorphism algebra is given by $\bi{(B\ox_\k A\op)}{(y,x)}{(y,x)}
= \bi{B}{y}{y} \ox_\k \bi{A}{x}{x}$ by the Yoneda Lemma,
which turns out to be a local algebra
because so are $\bi{B}{y}{y}$ and $\bi{A}{x}{x}$.
In particular, $B_y \ox_\k {}_x A$ is indecomposable
in $\biMod{B}{A}$.  Since the forgetful functor is additive and sends any nonzero object to nonzero, $B_y \ox_\k {}_x A$ as
an object in $\Ggr(\biMod{B}{A})$ is indecomposable.
Moreover by Remark \ref{rmk:classical-locbdd}(2),
both $A$ and $B$ are locally bounded categories.
Thus the indecomposable object $B_y \ox_\k {}_xA$ is finite-dimensional,
and hence it has a local endomorphism algebra in $\Ggr(\biMod{B}{A})$.
By Remark \ref{rmk:skeleton-Ggr},
so does $(S/G)_y \ox_\k {}_x(R/G)$ as an object in
$\Ggr(\biMod{S/G}{R/G})$.
Therefore, $M \iso ((M/G)\# G)'$ has a local endomorphism algebra as well
by Theorem \ref{thm:CM-dual-bimod}.
\end{proof}

\begin{prp}
\label{prp:fgproj-Ginv-bimod}
Let $(R, X)$ and $(S, Y)$ be locally bounded catego\-ries
with free $G$-actions. 
Then each f.g.\ projective $G$-invariant $S$-$R$-bimodule $(P, \ps)$ has the following form\footnote{
$n=0$ means that $P = 0$.
}:
\begin{equation}
\label{eq:fgproj-Ginv}
\Ds_{i=1}^n \Ds_{a\in G}(S_{ay_i} \ox_\k {}_{ax_i}R)
\end{equation}
for some $(x_i, y_i)_{i=1}^n \in (R_0 \times S_0)^n$
with $n \ge 0$.
\end{prp}

\begin{proof}
By definition, $(P,\ps)$ is a direct summand of an object in
$\Ginv(\biMod{S}{R})$ of the form \eqref{eq:fgproj-Ginv}.
Then the assertion follows by Lemma \ref{lem:ind-proj-G-inv}
and the Krull--Schmidt--Azumaya theorem.
\end{proof}

\begin{rmk}
\label{rmk:G-stable-skeleton}
For a $G$-category $R = (R, X)$, we say that the action $X$
is {\em free on isomorphism classes} if for any $1 \ne a \in G$
and $x \in R_0$, we have $ax:= X(a)x \not\iso x$ in $R$.
If this is the case, we have the following remarks.

\begin{enumerate}
\item 
$R$ has a skeleton $R'$ such that
$R'_0$ is $G$-{\em stable}, i.e., $ax \in R'_0$ for all $a \in G$ and $x \in R'_0$; and $X$ restricts to a $G$-action $X'$ on $R'$, and hence
the inclusion functor $(R', X') \incl (R, X)$ turns out to be
a strictly $G$-equivariant equivalence.
Therefore, the tensor product functor
$\blank\ox_{R'} (\bi{R}{R'}{R}) \colon \Mod_{R'}$ $\to \Mod_R$
becomes an equivalence of $G$-categories.
\item 
Assume that $S$ is also a $G$-category whose $G$-action is free on isomorphism classes with a $G$-stable skeleton $S'$.
Then we have an equivalence
$$
\Ginv(\biMod{S'}{R'}) \to \Ginv(\biMod{S}{R})
$$
sending $(M, \ph)$ to
$(S \ox_{S'} M \ox_{R'}R,\, (Y(a) \ox_{S'} \ph_a \ox_{R'}X(a))_{a\in G})$
for all $(M, \ph) \in \Ginv(\biMod{S'}{R'})_0$.
\end{enumerate}
\end{rmk}

\begin{prp}
Assume that both $A$ and $B$ are $G$-graded
locally bounded categories.
If both $G$-actions on $A\# G$ and $B\# G$ are free on
isomorphism classes (e.g. this holds if 
both $A\# G$ and $B\# G$ are basic categories,
cf. Lemma \ref{lem:smashprod-almost-loc-findim}),
then each f.g.\ projective $G$-graded bimodule has the following form:
$$
\Ds_{i=1}^n \ro_{b_i}(B_{y_i}) \ox_\k \la_{a_i}({}_{x_i}A)
$$
for some $(a_i, b_i)_i \in (G\times G)^n$ and $(x_i, y_i)_i \in (A_0 \times B_0)^n$ with $n \ge 0$.
\end{prp}

\begin{proof}
It follows from the assumption by Remark \ref{rmk:G-stable-skeleton}(1) that each of $A\# G$ and $B \# G$ has a
$G$-stable skeleton, say $R$ and $S$, respectively.
We regard $R$ and $S$ as $G$-categories by the $G$-actions
induced from those of $A\# G$ and $S\# G$, respectively.
Then by Proposition \ref{lem:smashprod-almost-loc-findim},
both $R$ and $S$ are locally bounded categories
with free $G$-actions,
and hence by Proposition \ref{prp:fgproj-Ginv-bimod},
any f.g.\ projective object in 
$\Ginv(\biMod{S}{R})$ has the form
\begin{equation}
\label{eq:fgproj-S-R}
\Ds_{i=1}^n \Ds_{a\in G}(S_{ay_i^{(b_i\inv)}} \ox_\k {}_{ax_i^{(a_i)}}R)
\end{equation}
for some $(x_i^{(a_i)}, y_i^{(b_i\inv)}) \in R_0 \times S_0 \subseteq
(A\# G)_0 \times (B\# G)_0$ with $n \ge 0$.
By Remark \ref{rmk:G-stable-skeleton}(2),
the inclusion functors $R \incl A\# G$ 
and $S \incl B\# G$ induce an equivalence
$\Ginv(\biMod{S}{R}) \to \Ginv(\biMod{B\# G}{A\# G})$ that 
sends \eqref{eq:fgproj-S-R} to
$$
\Ds_{i=1}^n \Ds_{a\in G}((B\# G)_{ay_i^{(b_i\inv)}} \ox_\k {}_{ax_i^{(a_i)}}(A\# G)).
$$
Hence the statement follows by Theorem \ref{thm:CM-dual-bimod} and \eqref{eq:formula/G'}.
\end{proof}

\begin{rmk}
Assume that $G$ is abelian.
Then by Remark \ref{rmk:G-gr-bimod-left}, we have an identification
$\Ggr(\biMod{B}{A}) = \Ggr({}_{B \ox_\k A\op}\Mod)$.
Hence by Proposition \ref{prp:graded_proj_over_loc_fd}, each f.g.\ projective $G$-graded
$B$-$A$-bimodule has the form
$$
\Ds_{i=1}^n \si_{a_i}(B_{y_i} \ox_\k {}_{x_i}A)
$$
for some $(a_i)_{i=1}^n \in G^n$ and $(x_i, y_i)_{i=1}^n \in (A_0 \times B_0)^n$ with $n \ge 0$.

Note that we can easily verify the following by definitions:
$$
\si_a(B_y \ox_\k xA) = \si_a(B_y) \ox_\k xA = B_y \ox_\k \si_a(xA) 
$$
for all $a \in G$, $(x, y) \in A_0 \times B_0$.
Therefore in addition, for each $b \in G$, we have
$$
\si_{ab}(B_y \ox_\k xA) = \si_b(B_y) \ox_\k \si_a(xA). 
$$
\end{rmk}

\section{Properties of smash products}

In this section, we collect additional properties of smash products
giving relationships with tensor products and f.g.\ projectivity.

\begin{prp}
\label{prp:tensor-smash}
Let $C$ be a $G$-graded small $\k$-category, and ${}_CN_B, {}_BM_A$ $G$-graded bimodules.
Then
\begin{enumerate}
\item
$N\ox_B M$ is a $G$-graded $C$-$A$-bimodule.
\item
$(N\ox_B M)\# G \iso (N\# G)\ox_{B\# G}(M\# G)$
in $\Ginv(\biMod{C\# G}{A\# G})$.
\end{enumerate}
\end{prp}

\begin{proof}
(1)
Let $(z,x) \in C_0 \times A_0$.
Then
\[
\begin{aligned}
{}_z(N \ox_B M)_x &= ({}_zN) \ox_B (M_x)
= (\textstyle\Ds_{y\in B_0} {}_zN_y \ox_\k {}_yM_x)/{}_zI_x\\
&=(\textstyle\Ds_{y\in B_0}(\Ds_{b \in G}{}_zN_y^b) \ox_\k (\Ds_{a \in G}{}_yM^a_x))/{}_zI_x\\
&=(\textstyle\Ds_{y\in B_0}\Ds_{a,b \in G} {}_zN_y^b \ox_\k {}_yM^a_x)/{}_zI_x\\
&= (\textstyle\Ds_{y\in B_0}\Ds_{c\in G}\Ds_{\smat{a,b \in G\\ba=c}}{}_zN_y^b \ox_\k {}_yM^a_x)/{}_zI_x \\
&= \textstyle\Ds_{c\in G}(\Ds_{y\in B_0}\Ds_{\smat{a,b \in G\\ba=c}}{}_zN_y^b \ox_\k {}_yM^a_x)/\Ds_{c \in G}{}_zI^c_x \\
&= \textstyle\Ds_{c\in G}{}_z(N \ox_B M)^c_x,
\end{aligned}
\]
where we set
\[
\begin{aligned}
{}_zI_x^c &:= \langle v\ox su - vs \ox u \mid (v,s,u) \in \bi{N}{z}{y'}^{b'}\times \bi{B}{y'}{y}^d\times \bi{M}{y}{x}^{a'}, y, y' \in B_0,\\
&\hspace{21em} a',b', d \in G, b'da'=c \rangle\\
&\subseteq \textstyle\Ds_{y\in B_0}\Ds_{\smat{a,b \in G\\ba=c}}{}_zN_y^b \ox_\k {}_yM^a_x
\end{aligned}
\]
because for each $c \in G$, we have
$
\{(b'd, a'), (b', da') \mid a', d, b' \in G, b'da' = c\}
\subseteq \{(b,a) \in G \times G \mid ba = c\}$,
and hence
\[
{}_z(N \ox_B M)^c_x :=
\textstyle(\Ds_{y\in B_0}\Ds_{\smat{a,b \in G\\ba=c}}{}_zN_y^b \ox_\k {}_yM^a_x)/{}_zI^c_x.
\]

(2)
For any objects $(z^{(c)}, x^{(a)})$ in $(C\#G)_0 \times (A\#G)_0$,
we first show that both hand sides of (2) at $(z^{(c)}, x^{(a)})$ \textcolor{blue}{are} equal to
\[
\textstyle(\Ds_{y\in B_0} \Ds_{b\in G} \bi{N}{z}{y}^{c\inv b} \ox_\k \bi{M}{y}{x}^{b\inv a})/ {}_zI^{c\inv a}_x.
\]
The LHS at $(z^{(c)}, x^{(a)})$ is given by
$$
\bi{((N\ox_B M)\# G)}{z^{(c)}}{x^{(a)}}:= \bi{(N\ox_B M)}{z}{x}^{c\inv a}=
\textstyle(\Ds_{y\in B_0} \Ds_{b\in G} \bi{N}{z}{y}^{c\inv b} \ox_\k \bi{M}{y}{x}^{b\inv a})/ {}_zI^{c\inv a}_x.
$$
The RHS at $(z^{(c)}, x^{(a)})$ is given by
\[
\begin{aligned}
{}_{z^{(c)}}((N\# G)&\ox_{B\# G}(M\# G))_{x^{(a)}}
=\bi{(N\# G)\ox_{B\# G}(M\# G)}{z^{(c)}}{x^{(a)}}\\
&= \textstyle(\Ds_{y^{(b)}\in (B\# G)_0} {}_{z^{(c)}}(N\# G)_{y^{(b)}} \ox_\k {}_{y^{(b)}}(M\# G)_{x^{(a)}})/{}_{z^{(c)}}I_{x^{(a)}}\\
&=
\textstyle(\Ds_{y\in B_0}\Ds_{b \in G}({}_zN_y^{c\inv b} \ox_\k {}_yM^{b\inv a}_x))/{}_zI^{c\inv a}_x,
\end{aligned}
\]
where we used the fact that ${}_{z^{(c)}}I_{x^{(a)}}={}_zI_x^{c\inv a}$,
which is shown as follows:
\[
\begin{aligned}
&{}_{z^{(c)}}I_{x^{(a)}}\\
&:= \langle v\ox su - vs \ox u \mid (v,s,u) \in \bi{(N\# G)}{z^{(c)}}{y'^{(b)}}\times \bi{(B\# G)}{y'^{(b)}}{y^{(d)}} \times \bi{(M\# G)}{y^{(d)}}{x^{(a)}},\\
&\hspace{120pt} y, y' \in B_0,
b,d \in G \rangle\\
&:= \langle v\ox su - vs \ox u \mid (v,s,u) \in {}_zN_{y'}^{c\inv b}\times {}_{y'}B_y^{b\inv d} \times {}_yM_x^{d\inv a}, y, y' \in B_0,
b,d \in G \rangle\\
&= {}_zI_x^{c\inv a}
\subseteq \Ds_{y\in B_0}\Ds_{e \in G}{}_zN_y^{c\inv e} \ox_\k {}_yM^{e\inv a}_x, 
\end{aligned}
\]
where note that for any $a, c \in G$, we have
$
\{(c\inv b, b\inv a), (c\inv d, d\inv a)  \mid b, d \in G\}
\subseteq \{(c\inv e,e\inv a) \in G \times G \mid e \in G\}$.
Hence $(N\ox_B M)\# G = (N\# G)\ox_{B\# G}(M\# G)$
as $C\# G$-$A\# G$-bimodules.
\par
Next we show that they have the same $G$-invariant structures.
Namely, we show the commutativity of the following diagram
for all $d \in G$, $x^{(a)} \in A \# G$, $z^{(c)} \in C \# G$:
\[
\begin{tikzcd}
\Nname{NM}{}_{z^{(c)}}((N \ox_B M)\# G)_{x^{(a)}} &
\Nname{NMG} {}_{z^{(c)}}((N\# G)\ox_{B\# G}(M\# G))_{x^{(a)}}\\
\Nname{NMd}{}_{z^{(dc)}}((N \ox_B M)\# G)_{x^{(da)}}&
\Nname{NMGd} {}_{z^{(dc)}}((N\# G)\ox_{B\# G}(M\# G))_{x^{(da)}}
\Ar{NM}{NMd}{"\bi{(\ph_d)}{z^{(c)}}{x^{(a)}}" '}
\Ar{NMG}{NMGd}{"\Ds_{y^{(b)} \in (B\#G)_0}\bi{(\ph^N_d)}{z^{(c)}}{y^{(b)}} \ox_\k \bi{(\ph^M_d)}{y^{(b)}}{x^{(a)}}"}
\Ar{NMd}{NMGd}{equal}
\Ar{NM}{NMG}{equal}
\end{tikzcd}.
\]
Let $u \in {}_{z^{(c)}}((N \ox_B M)\# G)_{x^{(a)}}$.
Since the ${}_{z^{(c)}}((N \ox_B M)\# G)_{x^{(a)}}$ is equal to
\[
{}_z(N\ox_B M)_{x}^{c\inv a}= 
\textstyle(\Ds_{y\in B_0}\Ds_{b \in G}{}_zN_y^{c\inv b} \ox_\k {}_yM^{b\inv a}_x)/{}_zI^{c\inv a}_x,
\]
we may take $u = (v_y^b \ovl{\ox} w_y^b)_{y^{(b)}\in (B\# G)_0}$ with $v_y^b \in {}_zN_y^{c\inv b}$,
$w_y^b \in {}_yM^{b\inv a}_x$, where $v_y^b\ovl{\ox}w_y^b:=(v_y^b \ox w_y^b)_{y^{(b)}} + {}_zI^{c\inv a}_x$.
We have to show that
\[
\bi{(\ph_d)}{z^{(c)}}{x^{(a)}}(u)
= \ovl{\bi{(\ph^{(N\# G)\ox_{B\# G} (M\# G)}_a)}{z^{(c)}}{x^{(a)}}}(u).
\]
By \eqref{eq:tensor} and \eqref{eq:str-smash}, RHS is equal to
\[
\begin{aligned}
&\left(\textstyle\ovl{\Ds_{y^{(b)} \in (B\#G)_0}\bi{(\ph^N_d)}{z^{(c)}}{y^{(b)}} \ox_\k \bi{(\ph^M_d)}{y^{(b)}}{x^{(a)}}}\right)(u)\\
&=\ovl{\left(\bi{(\ph^N_d)}{z^{(c)}}{y^{(b)}}(v_y^b) \ox_\k \bi{(\ph^M_d)}{y^{(b)}}{x^{(a)}}(w_y^b)\right)}_{y^{(b)} \in (B\#G)_0}\\
&=\left(v_y^b \ovl{\ox}_\k w_y^b\right)_{y^{(b)} \in (B\#G)_0}\\
&=u\\
&= \bi{(\ph_d)}{z^{(c)}}{x^{(a)}}(u).
\end{aligned}
\]
Hence as $G$-invariant bimodules, we have
\[
(N\ox_B M)\# G = (N\# G)\ox_{B\# G}(M\# G)
\]
if we use our explicit definition of tensor products.
Since in a categorical sense, the tensor products are defined up to natural isomorphisms, we just say that they are (naturally) isomorphic.
\end{proof}

\begin{cor}
\label{cor:tensor-over-k-smash}
Let $(x, y) \in R_0 \times S_0$.
Then we have
$$
(B_y \ox_\k {}_xA)\# G \iso \Ds_{a \in G}(B \# G)_{y^{(a)}}\ox_\k {}_{x^{(a)}}(A\# G).
$$
\end{cor}

\begin{proof}
Regard
$\k$ as a category having only one object $*$,
and $B_{y_i}$ (resp.\ ${}_{x_i}A$) as a $B$-$\k$-bimodule
(resp.\ $\k$-$A$-bimodule).
Then by Proposition \ref{prp:tensor-smash}(2), we have
$(B_{y}\ox_\k {}_{x}A)\# G \iso (B_{y})\# G \ox_{\k\# G} ({}_{x}A) \# G$.
Note that the following hold:
\[
\begin{aligned}
{}_{*^{(a)}}(({}_{x}A) \# G) &= {}_{x^{(a)}}(A \# G), \text{ and}\\
{}((B_{y}) \# G)_{*^{(a)}} &= (B \# G)_{y^{(a)}}
\end{aligned}
\]
for all $a \in G$.

Indeed, for each ${x'}^{(b)} \in ({}_xA \# G)_0\ (x' \in A_0,\, b \in G)$, we have
\[
{}_{*^{(a)}}(({}_{x}A) \# G)_{{x'}^{(b)}} = \Bi(x,A,x')^{a\inv b}
= {}_{x^{(a)}}(A \# G)_{{x'}^{(b)}}.
\]
Moreover, for each morphism
${}^{(b)}f^{(c)} \in \bi{(A\# G)}{{x'}^{(b)}}{{x''}^{(c)}}
= \{b\} \times \bi{A}{x'}{x''}^{b\inv c} \times \{c\}$,
$$
{}_{*^{(a)}}(({}_{x}A) \# G)_{({}^{(b)}{f}^{(c)})}
= {}_{x^{(a)}}(A \# G)_{({}^{(b)}{f}^{(c)})}
$$
because
both are the map $\Bi(x,A,x')^{a\inv b} \to \Bi(x,A,x'')^{a\inv c}$
defined by the right multiplication by $f$.
The remaining equality is shown similarly.

Then by noting that $\k \# G \iso \k^G$, we have
\[
\begin{aligned}
(B_{y}\ox_\k {}_{x}A)\# G
&\iso (B_{y})\# G \ox_{\k\# G} ({}_{x}A) \# G\\
&=
\Ds_{a \in G}((B_{y})\# G)_{*^{(a)}} \ox_{\k} {}_{*^{(a)}}(({}_{x}A) \# G)\\
&=
\Ds_{b \in G}(B\# G)_{y^{(a)}} \ox_{\k} {}_{x^{(a)}}(A \# G).
\end{aligned}
\]
\end{proof}

\begin{prp}
\label{prp:one-sided-proj-smashprod}
Let ${}_BM_A$ be a finitely generated $G$-graded bimodule.
Then the following statements hold.
\begin{enumerate}
\item
Assume that $A$ is $\k$-projective, and
let $x \in A_0$.
If $M_x$ is  
finitely generated projective as a left $B$-module,
then so is $(M\# G)_{x^{(a)}}$ as a left $B\# G$-module
for all $a \in G$.
\item
Assume that $B$ is $\k$-projective, and
let $y \in B_0$. 
If ${}_yM$ is 
finitely generated projective as a right $A$-module,
then so is ${}_{y^{(b)}}(M\# G)$ as a right $A\# G$-module
for all $b \in G$.
\end{enumerate}
\end{prp}

\begin{proof}
Since ${}_BM_A$ is a finitely generated $G$-graded bimodule,
we have an epimorphism
\[
f \colon \Ds_{i= 1}^n B_{y_i}\ox_\k {}_{x_i}A \to M
\]
in the category $\Ggr(\biMod{B}{A})$ for some $x_i \in A_0, \, y_i \in B_0$, $n \in \bbN$.
By Theorem \ref{thm:CM-dual-bimod}, $?\# G$ is
right exact, and hence
$f$ yields an epimorphism
\[
f\# G \colon \Ds_{i= 1}^n (B_{y_i}\ox_\k {}_{x_i}A)\# G \to M\# G
\]
of $(B\# G)$-$(A\# G)$-bimodules.

(1) $f\# G$ yields an epimorphism
\[
(f\# G)_{x^{(a)}} \colon (\Ds_{i= 1}^n (B_{y_i}\ox_\k {}_{x_i}A)\# G)_{x^{(a)}} \to (M\# G)_{x^{(a)}}
\]
of left $B\# G$-modules for all $a \in G$.

We have to show that

(a) $((B_{y_i}\ox_\k {}_{x_i}A)\# G)_{x^{(a)}}$
is a projective $B\# G$-module; and

(b) $(f\# G)_{x^{(a)}}$ is a retraction.

\bigskip

(a)
By Corollary \ref{cor:tensor-over-k-smash},
for each $i = 1, \dots, n$, we have
\[
\begin{aligned}
((B_{y_i}\ox_\k {}_{x_i}A)\# G)_{x^{(a)}}
&=
\Ds_{b \in G}(B\# G)_{y_i^{(b)}} \ox_{\k} {}_{x_i^{(b)}}(A \# G)_{x^{(a)}}\\
&\iso \Ds_{b \in G}(B\# G)_{y_i^{(b)}} \ox_{\k} {}_{x_i}A_x^{b\inv a},
\end{aligned}
\]
which is a projective left $B\# G$-module because
$\Ds_{b\in G}{}_{x_i}A_x^{b\inv a} \iso {}_{x_i}A_x$ is $\k$-projective by assumption,
and hence for each $b \in G$,
${}_{x_i}A_x^{b\inv a}$ is also $\k$-projective,
which shows that $(B\# G)_{y_i^{(b)}} \ox_{\k} {}_{x_i}A_x^{b\inv a}$
is projective over $B\# G$.

(b)
Set $N:= \Ds_{i=1}^n B_{y_i} \ox_\k {}_{x_i}A$ for short.
Then we have an epimorphism $f \colon N \to M$
in the category $\Ggr(\biMod{B}{A})$, which  yields an epimorphism
$f_x \colon N_x \to M_x$ in the category $\Ggr(\ltMod{B})$
for all $x \in A_0$.
Since $M_x$ is projective in $\Ggr(\ltMod{B})$
by Proposition \ref{prp:graded-prj},
$f_x$ is a retraction in $\Ggr(\ltMod{B})$.
Thus there exists a morphism $s_x \colon M_x \to N_x$
in the category $\Ggr(\ltMod{B})$ satisfying $f_x s_x = \id_{M_x}$.
We define a morphism
$s_{x^{(a)}} \colon (M\# G)_{x^{(a)}} \to (N\# G)_{x^{(a)}}$
of left $B\# G$-modules
that is a section of $(f\# G)_{x^{(a)}}$,
which will show the assertion.

Take any $y^{(b)} \in (B\# G)_0$.
Then $\Bi(y^{(b)}, (M\# G), x^{(a)}) = \Bi(y,M,x)^{b\inv a}$.
Therefore, we can define a morphism
$\Bi(y^{(b)},s,x^{(a)}) \colon \Bi(y^{(b)},M\# G,x^{(a)}) \to
\Bi(y^{(b)},N\# G,x^{(a)})$ of $\k$-modules as the restriction of
$s_x$ because
$s_x(\Bi(y,M,x)^{b\inv a}) \subseteq \Bi(y,N,x)^{b\inv a}$.
We define $s_{x^{(a)}} \colon (M\# G)_{x^{(a)}} \to (N\# G)_{x^{(a)}}$
by setting
$s_{x^{(a)}}:= (\Bi(y^{(b)},s,x^{(a)}))_{y^{(b)} \in (B\# G)_0}$.
This becomes a morphism of left $B\# G$-modules.
Indeed, since
$s_x$ is a morphism of $G$-graded left $B$-modules for all $x \in A_0$, we have the following commutative diagram:
\[
\begin{tikzcd}
\Nname{yMx}\Bi(y,M,x)^{b\inv a} &&& \Nname{yNx}\Bi(y,N,x)^{b\inv a}\\
&\Nname{yMGx}\Bi(y^{(b)}, (M\# G), x^{(a)}) & \Nname{yNGx}\Bi(y^{(b)}, (N\# G), x^{(a)})\\
&\Nname{y'MGx}\Bi({y'}^{(b')}, (M\# G), x^{(a)}) & \Nname{y'NGx}\Bi({y'}^{(b')}, (N\# G), x^{(a)})\\
\Nname{y'Mx}\Bi(y',M,x)^{{b'}\inv a} &&& \Nname{y'Nx}\Bi(y',N,x)^{{b'}\inv a}
\Ar{yMGx}{yNGx}{"{\Bi(y^{(b)},s,x^{(a)})}"}
\Ar{y'MGx}{y'NGx}{"{\Bi({y'}^{(b')},s,x^{(a)})}"}
\Ar{yMx}{yMGx}{equal}
\Ar{yNx}{yNGx}{equal}
\Ar{y'Mx}{y'MGx}{equal}
\Ar{y'Nx}{y'NGx}{equal}
\Ar{yMGx}{y'MGx}{"{\Bi(g,M\# G, x^{(a)})}" '}
\Ar{yNGx}{y'NGx}{"{\Bi(g,N\# G, x^{(a)})}" }
\Ar{yMx}{yNx}{""}
\Ar{yMx}{yNx}{"{\Bi(y,s,x)}"}
\Ar{y'Mx}{y'Nx}{"{\Bi({y'},s,x)}"}
\Ar{yMx}{y'Mx}{"{\Bi(g,M, x)}" '}
\Ar{yNx}{y'Nx}{"{\Bi(g,N, x)}" }
\end{tikzcd}
\]
for all $g \in (B\# G)(y^{(b)}, {y'}^{(b')}) = \Bi(y',B,y)^{{b'}\inv b}$.
Here, for each $y^{(b)} \in (B\# G)_0$, we have
\[
{}_{y^{(b)}}(f\# G)_{x^{(a)}}\cdot {}_{y^{(b)}}s_{x^{(a)}} = 
(f_x \cdot s_x)|_{\Bi(y,M,x)^{b\inv a}} = \id_{\Bi(y,M,x)^{b\inv a}}.
\]
Hence $(f\# G)_{x^{(a)}}$ is a retraction of $B\# G$-modules.

The proof of statement (2) is similar to that of (1).
\end{proof}

\begin{prp}\label{prp:smsh-gr-inv}
Let $P$ be a f.g.\ projective  $G$-graded $B$-$A$-bimodule.
Then $P\# G$ is a f.g.\ projective $G$-invariant $B\# G$-$A\# G$-bimodule.
\end{prp}

\begin{proof}
Since ${}_BP_A$ is a finitely generated
$G$-graded bimodule,
we have a retraction
\[
\Ds_{i= 1}^n B_{y_i}\ox_\k {}_{x_i}A \to P
\]
in $\Ggr(\biMod{B}{A})$ for some $(x_i, y_i)_{i=1}^n \in (A_0 \times B_0)^n$ with
$n$ a non-negative integer.
The functor $?\# G$ sends it to a retraction
\[
\Ds_{i= 1}^n (B_{y_i}\ox_\k {}_{x_i}A)\# G \to P\# G
\]
in $\Ginv(\biMod{B\# G}{A\# G})$, where for each $i = 1, \dots, n$
we have
\[
\begin{aligned}
(B_{y_i}\ox_\k {}_{x_i}A)\# G
&=
\Ds_{b \in G}(B\# G)_{y_i^{(b)}} \ox_{\k} {}_{x_i^{(b)}}(A \# G)\\
&=\Ds_{b \in G}(B\# G)_{b(y_i^{(1)})} \ox_{\k} {}_{b(x_i^{(1)})}(A \# G)
\end{aligned}
\]
by Corollary \ref{cor:tensor-over-k-smash}.
Hence 
$P\# G$ is a f.g.\ projective $G$-invariant $(B\# G)$-$(A\# G)$-bimodule.
\end{proof}

\section{Applications}

Throughout this section we assume that
$R, S$ are small $\k$-categories with $G$-actions, and
$A, B$ are $G$-graded small $\k$-categories.

\subsection{Morita equivalences}

\begin{dfn}
\label{dfn:inv-graded-Morita}
A bimodule $\bi{L}{\calD}{\calC}$ is said to be \emph{finitely generated projective on each side}
if both ${}_\calD L$ and $L_\calC$ are finitely generated projective
(see Definition \ref{dfn:1sided-fgp}).
\begin{enumerate}
\item
A pair $({}_SM_R, {}_RN_S)$ of bimodules is said to induce a $G$-{\em invariant Morita equivalence} between $R$ and $S$ if  
the bimodules ${}_SM_R$ and ${}_RN_S$ are
$G$-invariant and finitely generated projective on each side
such that 
$$
\begin{aligned}
N \ox_S M &\iso R \text{ and}\\
M \ox_R N &\iso S
\end{aligned}
$$
as $G$-invariant bimodules.
\item
A pair $({}_BM_A, {}_AN_B)$ of bimodules is said to induce a $G$-{\em graded Morita equivalence} between $A$ and $B$ if 
the bimodules ${}_BM_A$ and ${}_AN_B$ are
$G$-graded and finitely generated projective on each side
such that 
$$
\begin{aligned}
N \ox_B M &\iso A\text{ and} \\
M \ox_A N &\iso B
\end{aligned}
$$
as $G$-graded bimodules.
\end{enumerate}
\end{dfn}

\begin{lem}
\label{lem:orbi-smash-fgp}
Let $M$ be an $S$-$R$-bimodule.
If $(M/G)\# G$ is finitely generated projective as
a right $(R/G)\# G$-module and as a left $(S/G)\# G$-module,
then
$M$ is finitely generated projective as a right $R$-module
and as a left $S$-module.
\end{lem}

\begin{proof}
We only show that $M$ is finitely generated as a right $R$-module
because the rest is proved similarly.
Take any $y \in S_0$.
It is enough to show that ${}_yM$ is a f.g.\ projective right $R$-module.
Note that ${}_yM \iso {}_y((M/G)\# G)'$ as right $R$-modules.
By the equivalences
$$
\biMod{(S/G)\#G}{(R/G)\#G} \ya{\biMod{\ze_S}{(R/G)\# G}} \biMod{S}{(R/G)\# G},
\ \text{ and }\ 
\biMod{}{(R/G)\# G} \ya{\biMod{}{\ze_R}} \biMod{}{R},
$$
the f.g.\ projective right $(R/G)\# G$-module ${}_{y^{(1)}}(M/G)\# G$
is sent first to the f.g.\ projective right $(R/G)\# G$-module ${}_{y^{(1)}}({}_{\ze_{S}}(M/G)\# G)$,
and then is sent to the f.g.\ projective right $R$-module ${}_y((M/G)\# G)'$,
which is isomorphic in $\Mod_R$ to ${}_yM$.
\end{proof}

Using statements in previous sections we obtain the following.

\begin{thm}
Assume that all of $R$, $S$, $A$ and $B$ are $\k$-projective.
Then the following statements hold.
\begin{enumerate}
\item[(1)]
Let ${}_SM_R$ and ${}_RN_S$ be $G$-invariant bimodules.
Then the pair $({}_SM_R, {}_RN_S)$ induces a $G$-invariant Morita equivalence between $R$ and $S$ if and only if the pair $(M/G, N/G)$ induces a $G$-graded
Morita equivalence between $R/G$ and $S/G$.
\item[(2)]
Let ${}_BM_A$ and ${}_AN_B$ be $G$-graded bimodules.
Then the pair $({}_BM_A, {}_AN_B)$ induces a $G$-graded Morita equivalence
between $A$ and $B$ if and only if
the pair $(M\# G, N\# G)$ induces a $G$-invariant Morita equivalence
between $A\# G$ and $B\# G$.
\end{enumerate}
\end{thm}

\begin{proof}
(1) (\implies).  Assume that a pair $({}_SM_R, {}_RN_S)$ of bimodules induces a $G$-invariant Morita equivalence between $R$ and $S$.
Then the bimodules
${}_SM_R$ and ${}_RN_S$ are
$G$-invariant and finitely generated projective on each side such that 
\begin{align}
N \ox_S M &\iso R \text{ and}\label{eqn:NM}\\
M \ox_R N &\iso S  \label{eqn:MN}
\end{align}
as $G$-invariant bimodules.
By Corollary \ref{cor:proj-orbit}, the four modules ${}_{S/G}M/G, M/G_{R/G}$,
${}_{R/G}N/G, N/G_{S/G}$ are finitely generated projective.
Apply the functor $?/G$ to \eqref{eqn:NM} to have $(N \ox_S M)/G \iso R/G$, which shows that 
$$
(N/G) \ox_{S/G}(M/G) \iso R/G
$$
as $G$-graded $R/G$-$R/G$-bimodules by Proposition \ref{prp:tensor-orbit}(2).

Similarly from \eqref{eqn:MN} we obtain the remaining isomorphism.

(2) (\implies). This is proved similarly by using
Propositions \ref{prp:one-sided-proj-smashprod}
and \ref{prp:tensor-smash}. 

(1) (\impliedby).
Assume that the pair $(M/G, N/G)$ induces a $G$-graded
Morita equivalence between $R/G$ and $S/G$.
Then by (2) (\implies), we see that
the pair $((M/G)\# G, (N/G)\# G)$ induces a $G$-invariant Morita equivalence.
Namely, it is a pair of $G$-invariant bimodules that are finitely generated projective on each side, 
and we have isomorphisms
\begin{equation}
\label{eq:Morita_eq_orbit_smash}
\begin{aligned}
((M/G)\# G)\ox_{(R/G)\#G} ((N/G)\# G) &\iso (S/G)\#G,\\
((N/G)\# G)\ox_{(S/G)\#G} ((M/G)\# G) &\iso (R/G)\#G.
\end{aligned}
\end{equation}
By Lemma \ref{lem:orbi-smash-fgp}, we see that
${}_SM,\, M_R,\, {}_RN$ and $N_S$ are finitely generated projective.
From \eqref{eq:Morita_eq_orbit_smash} it holds
by Proposition \ref{prp:tensor-smash} that
\begin{equation}
\label{eq:Morita_eq_smash}
\begin{aligned}
((M/G)\ox_{R/G} (N/G))\# G &\iso (S/G)\#G,\\
(N/G)\ox_{S/G} (M/G))\# G &\iso (R/G)\#G.
\end{aligned}
\end{equation}
Therefore, by Proposition \ref{prp:tensor-orbit},
we have
$$
\begin{aligned}
((M\ox_{R}N)/G)\# G &\iso (S/G)\#G,\\
((N\ox_{S} M)/G)\# G &\iso (R/G)\#G.
\end{aligned}
$$
Hence by Theorem \ref{thm:CM-dual-bimod},
$$
\begin{aligned}
M\ox_{R}N &\iso S,\\
N\ox_{S} M&\iso R
\end{aligned}
$$
as $S$-$S$-bimodules and as $R$-$R$-bimodules, respectively.
As a consequence, the pair $(M, N)$ induces a Morita equivalence.

(2) (\impliedby).
Similar proof as above works.
\end{proof}

\subsection{Stable equivalences of Morita type}

\begin{dfn}
\label{dfn:inv-graded-stable-Morita}
\begin{enumerate}
\item
A pair $({}_SM_R, {}_RN_S)$ of bimodules is said to induce a $G$-{\em invariant stable equivalence of Morita type} between $R$ and $S$ if the bimodules
${}_SM_R$ and ${}_RN_S$ are
$G$-invariant and finitely generated projective on each side such that 
$$
\begin{aligned}
N \ox_S M &\iso R \ds {}_RP_R \text{ and}\\
M \ox_R N &\iso S \ds {}_SQ_S
\end{aligned}
$$
 as $G$-invariant bimodules
for some f.g.\ projective $G$-invariant bimodules ${}_RP_R$ and ${}_SQ_S$
(see Definition \ref{dfn:fgp-Ginv-bimod}).
\item
A pair $({}_BM_A, {}_AN_B)$ of bimodules is said to induce a $G$-{\em graded stable equivalence of Morita type} between $A$ and $B$ if the bimodules 
${}_BM_A$ and ${}_AN_B$ are
$G$-graded and finitely generated projective on each side
such that 
$$
\begin{aligned}
N \ox_B M &\iso A \ds {}_AP_A\text{ and} \\
M \ox_A N &\iso B \ds {}_BQ_B
\end{aligned}
$$
as $G$-graded bimodules
for some f.g.\ projective $G$-graded  bimodules
${}_AP_A, {}_BQ_B$ (see Definition \ref{dfn:fgp-Ginv-bimod}).
\end{enumerate}
\end{dfn}

\begin{rmk}
\label{rmk:Ggr-prj}
Note that $\bi{P}{A}{A}$ is a $G$-graded bimodule, but it does not
need to be a $G$-graded left $A^{\rme}$-module because $G$ is not necessarily
abelian.
Nevertheless, the projectivity of $P$ in $\Ggr(\biMod{A}{A})$ and
that in $\biMod{A}{A}$ are equivalent by Lemma \ref{lem:G-gr-free}.
\end{rmk}
Using statements in previous sections we obtain the following.

\begin{thm}
\label{thm:st-eq-M}
Assume that all of $R$, $S$, $A$ and $B$ are $\k$-projective. 
Then the following statements hold.
\begin{enumerate}
\item
A pair $({}_SM_R, {}_RN_S)$ of bimodules induces a $G$-invariant stable equivalence of Morita type between $R$ and $S$ if and only if the pair $(M/G, N/G)$ induces a $G$-graded
stable equivalence of Morita type between $R/G$ and $S/G$.
\item
A pair $({}_BM_A, {}_AN_B)$ of bimodules induces a $G$-graded stable equivalence of Morita type between $A$ and $B$ if and only if
the pair $(M\# G, N\# G)$ induces a $G$-invariant stable equivalence of Morita type
between $A\# G$ and $B\# G$.
\end{enumerate}
\end{thm}

\begin{proof}
(1)(\implies). Assume that a pair $({}_SM_R, {}_RN_S)$ of bimodules induces a $G$-invariant stable equivalence of Morita type between $R$ and $S$.
Then the bimodules
${}_SM_R$ and ${}_RN_S$ are
$G$-invariant and finitely generated projective on each side such that 
\begin{align}
N \ox_S M &\iso R \ds {}_RP_R \text{ and}\label{eqn:NM-st}\\
M \ox_R N &\iso S \ds {}_SQ_S  \label{eqn:MN-st}
\end{align}
 as $G$-invariant bimodules
for some finitely generated projective $G$-invariant
bimodules ${}_RP_R, {}_SQ_S$.
Apply the functor $?/G$ to \eqref{eqn:NM-st} to have $(N \ox_S M)/G \iso (R \ds P)/G$, which shows that 
$$
(N/G) \ox_{S/G}(M/G) \iso R/G \ds P/G
$$
as $G$-graded $R/G$-$R/G$-bimodules by Proposition \ref{prp:tensor-orbit}(2).
Here the bimodules $\bi{M/G}{S/G}{R/G}$ and $\bi{N_G}{R/G}{S/G}$
are $G$-graded bimodules and finitely generated projective on each side
by Corollary \ref{cor:proj-orbit},
and ${}_{R/G}P/G_{R/G}$ is a f.g.\ projective $G$-graded bimodule by Proposition \ref{prp:orbitr-bi}.
Similarly from \eqref{eqn:MN-st} we obtain the remaining isomorphism.

(2)(\implies). This is shown similarly by using Propositions \ref{prp:one-sided-proj-smashprod} and \ref{prp:smsh-gr-inv}.

(1)(\impliedby).
Assume that the pair $(M/G, N/G)$ induces a $G$-graded
stable equivalence of Morita type between $R/G$ and $S/G$.
Then by definition, it is a pair of $G$-graded bimodules
that are finitely generated projective on each side, and
we have
\begin{align*}
N/G \ox_{S/G} M/G &\iso R/G \ds {}_{R/G}P'_{R/G} \text{ and}\\
M/G \ox_{R/G} N/G &\iso S/G \ds {}_{S/G}Q'_{S/G}  
\end{align*}
for some finitely generated projective $G$-graded bimodules
${}_{R/G}P'_{R/G}$ and ${}_{S/G}Q'_{S/G}$.
Then by (2)(\implies), we see that
the pair $((M/G)\# G, (N/G)\# G)$ induces a $G$-invariant stable equivalence of Morita type.
Namely, it is a pair of $G$-invariant bimoudles 
that are finitely generated projective on each side, and we have isomorphisms
\begin{equation}
\label{eq:st-Morita_eq_orbit_smash}
\begin{aligned}
((M/G)\# G)\ox_{(R/G)\#G} ((N/G)\# G) &\iso ((S/G)\#G) \ds (Q'\#G),\\
((N/G)\# G)\ox_{(S/G)\#G} ((M/G)\# G) &\iso ((R/G)\#G) \ds (P'\#G).
\end{aligned}
\end{equation}
By Lemma \ref{lem:orbi-smash-fgp}, we see that the bimodules
$\bi{M}{S}{R},\, \bi{N}RS$ are $G$-invariant
and finitely generated projective on each side.
From \eqref{eq:st-Morita_eq_orbit_smash} it holds
by Proposition \ref{prp:tensor-smash} that
\begin{equation}
\label{eq:st-Morita_eq_smash}
\begin{aligned}
((M/G)\ox_{R/G} (N/G))\# G &\iso (S/G\ds Q')\#G,\\
(N/G)\ox_{S/G} (M/G))\# G &\iso (R/G\ds P')\#G.
\end{aligned}
\end{equation}
Therefore, by Proposition \ref{prp:tensor-orbit},
we have
$$
\begin{aligned}
((M\ox_{R}N)/G)\# G &\iso (S/G\ds Q')\#G,\\
((N\ox_{S} M)/G)\# G &\iso (R/G\ds P')\#G.
\end{aligned}
$$
Hence by Theorem \ref{thm:CM-dual-bimod},
$$
\begin{aligned}
M\ox_{R}N &\iso S\ds (Q'\# G)'\\\
N\ox_{S} M&\iso R\ds (P'\# G)'
\end{aligned}
$$
as $S$-$S$-bimodules and as $R$-$R$-bimodules, respectively.
Moreover, $(P'\# G)'$ (resp.\ $(Q'\# G)'$) is f.g.\ projective $G$-invariant
$R$-$R$-bimodule (resp.\ $S$-$S$-bimodule)
because $\ovl{\ze}$ in \eqref{eq:zeta-bar} is an equivalence, and
$P'\# G$ (resp.\ $Q'\# G$) is f.g.\ projective
$G$-invariant $((R/G)\# G)$-$((R/G)\# G)$-bimodule
(resp.\ $((S/G)\# G)$-$((S/G)\# G)$-bimodule)
by Proposition \ref{prp:smsh-gr-inv}.
As a consequence, the pair $(M, N)$ induces a stable equivalence of Morita type.

(2)(\impliedby). This is proved similarly. 
\end{proof}

\subsection{Singular equivalences of Morita type}

Let $\calC$ be a $\k$-category.
Then we denote by $\calC^\rme$ the enveloping category $\calC \ox_\k \calC\op$ of $\calC$.
Then each $\calC$-$\calC$-bimodule $M$ can be seen as a
left $\calC^\rme$-module by setting $(h,g)f:= hfg$
for all $h \in \Bi(y',\calC,y),\, f \in \Bi(y,\calC,x),\, g \in \Bi(x,\calC, x')$ and $y',\,y,\,x,\,x' \in \calC_0$.

\begin{dfn}
\label{dfn:inv-graded-singular-Morita}
\begin{enumerate}
\item
A pair $({}_SM_R, {}_RN_S)$ of bimodules is said to induce a $G$-{\em invariant singular equivalence of Morita type} between $R$ and $S$ if  
the bimodules ${}_SM_R$ and ${}_RN_S$ are
$G$-invariant and finitely generated projective on each side such that 
$$
\begin{aligned}
N \ox_S M &\iso R \ds {}_RP_R \text{ and}\\
M \ox_R N &\iso S \ds {}_SQ_S
\end{aligned}
$$
 as $G$-invariant bimodules
for some f.g.\ $G$-invariant
bimodules ${}_RP_R, {}_SQ_S$ of finite projective dimension in $\Ginv(\biMod{R}{R})$ and $\Ginv(\biMod{S}{S})$, respectively.
\item
A pair $({}_BM_A, {}_AN_B)$ of bimodules is said to induce a $G$-{\em graded singular equivalence of Morita type} between $A$ and $B$ if 
the bimodules ${}_BM_A$ and ${}_AN_B$ are
$G$-graded and finitely generated projective on each side such that 
$$
\begin{aligned}
N \ox_B M &\iso A \ds {}_AP_A\text{ and} \\
M \ox_A N &\iso B \ds {}_BQ_B
\end{aligned}
$$
as $G$-graded bimodules
for some finitely generated $G$-graded bimodules
${}_AP_A, {}_BQ_B$ of finite projective dimension
over $A^\rme$ and $B^\rme$, respectively.
\end{enumerate}
\end{dfn}

\begin{rmk}
\label{rmk:one-sided-proj}
(1) In Definition \ref{dfn:inv-graded-singular-Morita}(1),
note that both ${}_RP$ and $P_R$ (resp.\ ${}_SQ$ and $Q_S$)
turn out to be projective
because they are direct summands of $N \ox_S M$
(resp.\ $M \ox_R N$).

Similarly in Definition \ref{dfn:inv-graded-singular-Morita}(2),
all of ${}_AP, P_A, {}_BQ, Q_B$ turn out to be projective.

(2) 
In Definition \ref{dfn:inv-graded-singular-Morita}(1),
$P$ has a finite projective dimension in $\Ginv(\biMod{R}{R})$
if and only if so does in $\biMod{R}{R}= \ltMod{R^\rme}$
by Proposition \ref{prp:Ginv-prj}.

Similarly
in Definition \ref{dfn:inv-graded-singular-Morita}(2),
$P$ has a finite projective dimension in $\Ggr(\biMod{A}{A})$
if and only if so does in $\biMod{A}{A}= \ltMod{A^\rme}$
by Remark \ref{rmk:Ggr-prj}.
\end{rmk}

Recall that the \emph{bar resolution} of $S$ is defined to be the
following exact sequence of $S$-$S$-bimodules:
\begin{equation}
\tag{\#}
\cdots \to S \ox_\k S \ox_\k S \ya{d_1} S\ox_\k S \ya{d_0} S \to 0,
\end{equation}
where for each $i \ge 0$, $d_{i} \colon S^{\ox(i+2)} \to S^{\ox (i+1)}$
is defined by
$$
d_{i}(x_0 \ox x_1 \ox \cdots x_{i+1}):= \sum_{j=0}^{i}(-1)^j
(x_0 \ox x_1 \ox \cdots \ox (x_jx_{j+1})\ox\cdots \ox x_{i+1}).
$$

\begin{rmk}
\label{rmk:G-inv+G-gr}
(1) Since $S = (S, (Y_a)_{a\in G})$ becomes a $G$-invariant bimodule,
$S^{\ox i} = (S^{\ox i}, (Y_a^{\ox i})_{a \in G})$ turn out to be
a $G$-invariant bimodule for all $i \ge 1$.
By these $G$-invariant structures,
$d_{i} \colon S^{\ox(i+2)} \to S^{\ox (i+1)}$ are
morphisms in $\Ginv(\biMod{S}{S})$ for all $i \ge 0$.

(2) Similar remarks are valid for the $G$-graded case.
\end{rmk}

\begin{lem}
\label{lem:prj-ox-prj}
Let $\calC,\, \calD$ be small $\k$-categories, and 
${}_\calD Q,\, P_\calC$ projective modules.
Then $Q\ox_\k P$ is a projective $\calD$-$\calC$-bimodule.
\end{lem}

\begin{proof}
The functors
$\Hom_\calD(Q,\blank) \colon
{}_{\calD}\Mod \to \kMod$
and 
$\Hom_{\calC\op}(P,\blank) \colon \biMod{\calD}{\calC} \to {}_\calD\Mod$
are exact by assumption.
Hence by noting that we have an isomorphism 
$$
\Hom_{\calD\ox_\k\calC\op}(Q\ox_\k P, \blank)
\iso \Hom_{\calD}(Q, \Hom_{\calC\op}(P, \blank))
$$
of functors $\biMod{\calD}{\calC} \to \kMod$,
the functor $\Hom_{\calD\ox_\k\calC\op}(Q\ox_\k P, \blank)$ is exact as the composite of exact functors.
\end{proof}

\begin{lem}
\label{lem:proj-bimod--onesided-proj}
Let $\calC,\, \calD$ be small $\k$-categories, and
$M$ a projective $\calD$-$\calC$-bimodule.
Assume that $\calD$ is $\k$-projective.
Then $M_\calC$ is projective.
\end{lem}

\begin{proof}
Since $M$ is a projective $\calD$-$\calC$-bimodule,
as a left $\calD \ox_\k \calC\op$-module,
$M$ is a direct summand of $(\calD \ox_\k \calC\op)^{(I)}$
for some set $I$.
By assumption,
$\Bi(y',\calD,y) \ox_\k \calC\op$ becomes a projective right $\calC$-module for all $y,\, y' \in \calD_0$, which means that
$\calD \ox_\k \calC\op$ is a projective right $\calC$-module.
Hence so is $M$.
\end{proof}
\begin{lem}
\label{lem:proj-resol-bimod}
Let $\calC,\, \calD$ be small $\k$-categories. Assume that $\calC$ is $\k$-projective.
Let $P$ be an $\calC$-$\calD$-bimodule with both ${}_{\calC}P, P_{\calD}$ projective, and let $(\#)$ be the bar resolution
of $\calC$.
Then $(\#) \ox_{\calC} P$ is a projective resolution of the $\calC$-$\calD$-bimodule $\calC \ox_\calC P \iso P$.
\end{lem}

\begin{proof}
Since ${}_{\calC}P$ is projective, the sequence $(\#) \ox_{\calC} P$ is an exact
sequence in $\biMod{\calC}{\calD}$.
Since for each $n \ge 1$,
${}_{\calC} \calC^{\ox n}$ and $P_\calD$ are projective,
we see that $\calC^{\ox n} \ox_\k P$ is projective in $\biMod{\calC}{\calD}$
by Lemma \ref{lem:prj-ox-prj}.
Hence $(\#) \ox_\calC P$ becomes a projective resolution of
$\calC \ox_\calC P \iso P$ as a $\calC$-$\calD$-bimodule.
\end{proof}

\begin{lem}
\label{lem:proj-resol-RS}
Assume that $S$ is $\k$-projective.
If $P \in \Ginv(\biMod{S}{R})$ with ${}_SP, P_R$ projective,
then $(\#) \ox_S P$ is a projective resolution of $P$ as an
$S$-$R$-bimodule that is in $\Ginv(\biMod{S}{R})$.
\end{lem}

\begin{proof}
Lemma \ref{lem:proj-resol-bimod} shows that
$(\#) \ox_S P$ is a projective resolution of $P$ as an
$S$-$R$-bimodule.
Since $S = (\Bi(S,S,S), Y)$ is in $\Ginv(\biMod{S}{S})$,
and $P \in \Ginv(\biMod{S}{R})$, we have
$S^{\ox n} \ox_S P \in \Ginv(\biMod{S}{R})$
by Proposition \ref{prp:tensor-orbit}(1).
Hence by Remark \ref{rmk:G-inv+G-gr},
it is easy to verify that
$(\#) \ox_{\calC} P$ is in $\Ginv(\biMod{S}{R})$.
\end{proof}

\begin{lem}
Assume that $B$ is $\k$-projective.
If $P \in \Ggr(\biMod{B}{A})$ with ${}_BP, P_A$ projective,
then $(\#) \ox_B P$ is a projective resolution of $P$ as
a $B$-$A$-bimodule that is in $\Ggr(\biMod{B}{A})$.
\end{lem}

\begin{proof}
This follows by Remark \ref{rmk:G-inv+G-gr}, Lemma \ref{lem:proj-resol-bimod} and Proposition \ref{prp:tensor-smash}.
\end{proof}

\begin{prp}
\label{prp:fin-prj-dim-bimod}
If $P$ is a f.g.\ $G$-invariant
$R$-$R$-bimodule of finite projective dimension
over $R^\rme$, then $P/G$ is a f.g.\ 
$G$-graded $R/G$-$R/G$-bimodule of finite projective dimension over $(R/G)^\rme$.
\end{prp}

\begin{proof}
Since $P$ is a f.g.\ $G$-invariant 
$R$-$R$-bimodule,
there exists a finite set $I$ and a family $(x_i, y_i) \in (R_0 \times R_0)^I$ such that there exists an epimorphism
\[
\textstyle f \colon \Ds_{i \in I} (\Ds_{a \in G} R_{ay_i}\ox_\k {}_{ax_i}R) \to P
\]
in the category $\Ginv(\biMod{R}{R})$, where 
$\Ds_{a \in G} R_{ay_i}\ox_\k {}_{ax_i}R$
has the canonical $G$-invariant structure for all $i \in I$.
Since $?/G$ is exact, we have an epimorphism
\[
\textstyle f/G \colon \Ds_{i \in I} (\Ds_{a \in G} R_{ay_i}\ox_\k {}_{ax_i}R)/G \to P/G.
\]
By Proposition \ref{prp:orbitr-bi}
we have
$$
\textstyle(\Ds_{a \in G} R_{ay_i}\ox_\k {}_{ax_i}R)/G \iso (R/G)_{P(y_i)}\ox_\k {}_{P(x_i)}(R/G),
$$
which is canonically f.g.\ projective
$G$-graded $R/G$-$R/G$-bimodule.
Hence $P/G$ is an f.g.\ $R/G$-$R/G$-bimodule.

Let $(\#)$ be the bar resolution of $R$.
Then $(\#)\ox_R P$ is a projective resolution of $P$ in $\Ginv(\biMod{R}{R})$ by Remark \ref{rmk:one-sided-proj} and
Lemma \ref{lem:proj-resol-RS},
where $\Im (d_i\ox_R P)$
is projective $R^\rme$-module for some $i \ge 0$
because $P$ has finite projective dimension as a left $R^\rme$-module.
Hence $d_i\ox_R P \colon R^{\ox(i+2)}\ox_RP \to \Im (d_i\ox_R P)$
is a retraction, and $\Im (d_i\ox_R P)$ is a
canonically $G$-invariant f.g.\ projective $R$-$R$-bimodule.
Again by Proposition \ref{prp:orbitr-bi},
$\Im (d_i\ox_R P)/G$ is a $G$-graded f.g.\ projective
$R/G$-$R/G$-bimodule.
Therefore, the projective resolution $((\#)\ox_R P)/G \iso (\#)/G \ox_{R/G} P/G$ of $P/G$ as left $(R/G)^\rme$-module has some
projective image $\Im (d_i\ox_R P)/G$.
Thus $P/G$ has finite projective dimension as a left $(R/G)^\rme$-module.
\end{proof}

The following is proved similarly by using
Proposition \ref{prp:smsh-gr-inv}.

\begin{prp}
If $P$ is a f.g.\ $G$-graded $A$-$A$-bimodule of finite projective dimension over $A^\rme$, then $P\# G$ is a f.g.\  $G$-invariant $(A\# G)$-$(A\# G)$-bimodule of finite projective dimension
over $(A\# G)^\rme$.
\qed
\end{prp}

\begin{thm}
\label{thm:sing-eq-M}
Assume that all of $R$, $S$, $A$ and $B$ are $\k$-projective. 
Then the following statements hold.
\begin{enumerate}
\item
A pair $({}_SM_R, {}_RN_S)$ of bimodules induces a $G$-invariant singular equivalence of Morita type between $R$ and $S$ if and only if the pair $(M/G, N/G)$ induces a $G$-graded
singular equivalence of Morita type between $R/G$ and $S/G$.
\item
A pair $({}_BM_A, {}_AN_B)$ of bimodules induces a $G$-graded singular equivalence of Morita type between $A$ and $B$ if and only if
the pair $(M\# G, N\# G)$ induces a $G$-invariant singular equivalence of Morita type
between $A\# G$ and $B\# G$.
\end{enumerate}
\end{thm}

\begin{proof}
(1)(\implies). Assume that a pair $({}_SM_R, {}_RN_S)$ of bimodules induces a $G$-invariant singular equivalence of Morita type between $R$ and $S$.
Then the bimodules ${}_SM_R, {}_RN_S$ are
$G$-invariant and finitely generated projective on each side such that 
\begin{align}
N \ox_S M &\iso R \ds {}_RP_R \text{ and}\label{eqn:NM-stq}\\
M \ox_R N &\iso S \ds {}_SQ_S  \label{eqn:MN-stq}
\end{align}
 as $G$-invariant bimodules, where ${}_RP_R, {}_SQ_S$
are  f.g.\ $G$-invariant 
bimodules of finite projective dimensions in $\Ginv(\biMod{R}{R})$ and $\Ginv(\biMod{S}{S})$, respectively.
Apply the functor $?/G$ to \eqref{eqn:NM-stq} to have $(N \ox_S M)/G \iso (R \ds P)/G$, which shows that 
$$
(N/G) \ox_{S/G}(M/G) \iso R/G \ds P/G
$$
as $G$-graded $R/G$-$R/G$-bimodules by Proposition \ref{prp:tensor-orbit}(2).
Here ${}_{S/G}M/G, M/G_{R/G}$,
${}_{R/G}N/G, N/G_{S/G}$ are canonically $G$-graded projective by Corollary \ref{cor:proj-orbit}.
And ${}_{R/G}P/G_{R/G}$ has finite projective dimension by Proposition \ref{prp:fin-prj-dim-bimod}.
Similarly from \eqref{eqn:MN-stq} we obtain the remaining isomorphism.

(2)(\implies). This is proved similarly. 

(1)(\impliedby).
This follows by (2)(\implies) and the equivalences 
$(R/G)\# G \simeq R$ and $(S/G)\# G \simeq S$ of categories.

(2)(\impliedby).
This is proved similarly.
\end{proof}

\subsection{Singular equivalences of Morita type with level}
\begin{dfn}
\label{dfn:inv-graded-singular-Morita-lev}
\begin{enumerate}
\item
A pair $({}_SM_R, {}_RN_S)$ of bimodules is said to induce a $G$-{\em invariant singular equivalence of Morita type with level $l\geq 0$} between $R$ and $S$ if 
the bimodules ${}_SM_R, {}_RN_S$ are
$G$-invariant, and finitely generated projective on each side such that 
$$
\begin{aligned}
N \ox_S M &\iso \Om^l_{R^e} (R)\ds {}_RP_R
\text{ and}\\
M \ox_R N &\iso \Om^l_{S^e} (S) \ds {}_SQ_S
\end{aligned}
$$
 as $G$-invariant bimodules
for some f.g.\ projective $G$-invariant bimodules
${}_RP_R, {}_SQ_S$, where $\Om_{R^e}$ and $\Om_{S^e}$ denote the Heller shifts in $\Ginv(\biMod{R}{R})$ and $\Ginv(\biMod{S}{S})$, respectively.

\item
A pair $({}_BM_A, {}_AN_B)$ of bimodules is said to induce a $G$-{\em graded singular equivalence of Morita type with level $l\geq 0$} between $A$ and $B$ if 
the bimodules ${}_BM_A, {}_AN_B$ are
$G$-graded bimodules and finitely generated projective on each side such that 
$$
\begin{aligned}
N \ox_S M &\iso \Om^l_{R^e} (A)\ds {}_AP_A
\text{ and}\\
M \ox_R N &\iso \Om^l_{S^e} (B) \ds {}_BQ_B
\end{aligned}
$$
as $G$-graded bimodules
for some finitely generated projective 
$G$-graded bimodules
${}_AP_A, {}_BQ_B$.
\end{enumerate}
\end{dfn}

\begin{thm}
\label{thm:sing-eq-M-l}
Assume that all of $R$, $S$, $A$ and $B$ are $\k$-projective,
and let $l$ be a non-negative integer.
Then the following statements hold.
\begin{enumerate}
\item
A pair $({}_SM_R, {}_RN_S)$ of bimodules induces a $G$-invariant singular equivalence of Morita type with level $l$ between $R$ and $S$ if and only if the pair $(M/G, N/G)$ induces a $G$-graded
singular equivalence of Morita type with level $l$ between $R/G$ and $S/G$.
\item
A pair $({}_BM_A, {}_AN_B)$ of bimodules induces a $G$-graded singular equivalence of Morita type with level $l$ between $A$ and $B$ if and only if
the pair $(M\# G, N\# G)$ induces a $G$-invariant singular equivalence of Morita type with level $l$
between $A\# G$ and $B\# G$.
\end{enumerate}
\end{thm}

\begin{proof}
(1)(\implies). Assume that a pair $({}_SM_R, {}_RN_S)$ of bimodules induces a $G$-invariant singular equivalence of Morita type with level $l$ between $R$ and $S$.
Then the bimodules ${}_SM_R, {}_RN_S$ are $G$-invariant and finitely generated projective on each side such that 
\begin{align}
N \ox_S M &\iso \Om^l_{R^e}(R) \ds {}_RP_R \text{ and}\label{eqn:NM-stql}\\
M \ox_R N &\iso \Om^l_{S^e}(S) \ds {}_SQ_S  \label{eqn:MN-stql}
\end{align}
 as $G$-invariant bimodules, where ${}_RP_R, {}_SQ_S$
are finitely generated $G$-invariant bimodules with finite projective dimensions in $\Ginv(\biMod{R}{R})$ and $\Ginv(\biMod{S}{S})$, respectively. 
Apply the functor $?/G$ to \eqref{eqn:NM-stql} to have $(N \ox_S M)/G \iso (\Om^l_{R^e}(R) \ds P)/G$, which shows that 
$$
(N/G) \ox_{S/G}(M/G) \iso \Om^l_{R^e}(R)/G \ds P/G
$$
as $G$-graded $R/G$-$R/G$-bimodules by Proposition \ref{prp:tensor-orbit}(2).
Here ${}_{S/G}M/G, M/G_{R/G}$,
${}_{R/G}N/G, N/G_{S/G}$ are $G$-graded projective by Corollary \ref{cor:proj-orbit}, and ${}_{R/G}P/G_{R/G}$ is finitely generated $G$-graded projective by Proposition \ref{prp:orbitr-bi}.

Consider the bar resolution $(\#)$ of $R$.
Then since it is a projective resolution of $R$
as a left $R^\rme$-module, we have
$\Om^i_{R^\rme}(R) \iso \Im d_i$ for all $i \ge 0$.
Since $?/G$ is exact and $(\#)/G$ becomes a projective
resolution of $R/G$ as a left $(R/G)^\rme$-module, we have
$$
\Om^l_{R^\rme}(R)/G 
\iso (\Im d_l)/G
\iso \Om^l_{(R/G)^e}(R/G).
$$
Hence we have
$$
(N/G) \ox_{S/G}(M/G) \iso \Om^l_{(R/G)^e}(R/G) \ds P/G
$$
as desired.
The remaining isomorphism follows similarly.

(2)(\implies). This is proved similarly. 

(1)(\impliedby).
This follows by (2)(\implies) and the equivalences
$(R/G)\# G \simeq R$ and $(S/G)\# G \simeq S$ of categories.

(2)(\impliedby).
This is proved similarly.
\end{proof}

\section{Examples}
\label{sec:exmples}
In this section, we give examples of
$G$-invariant $S$-$R$-bimodule $M$
and $G$-graded $B$-$A$-bimodule $M'$ such that
$M/G \iso M'$, and hence $(M'\# G)' \simeq M$.
Moreover, we further define an $A$-$B$-bimodule $N'$
such that the pair $(M', N')$ induces a $G$-graded
stable equivalence of Morita type between
$A$ and $B$.
Hence by taking $M:= M'\# G,\, N:= N'\# G$, the pair
$(M,N)$ defines a $G$-invariant stable equivalence of Morita
type between $R:= A\# G$ and $S:= B\# G$.

\subsection{Presentation of bimodules as triangular matrix algebras}

In this section, we assume that
$S, R, B, A$ are path-categories of finite  bound quivers.
Therefore, they can be regarded as finite-dimensional algebras.
To express the $S$-$R$-bimodule $M$ and the $B$-$A$-bimodule $M'$,
we use bound quiver presentations of the triangular matrix algebras
$T(M):= \bmat{R & 0\\M & S}$ and $T(M'):= \bmat{A & 0\\M' & B}$, respectively, where we identify $M$ with $\bmat{0&0\\M&0} \subseteq T(M)$
and $M'$ with $\bmat{0&0\\M'&0} \subseteq T(M')$.
We refer the reader to \cite{Asa-Kim} and \cite{Asa18}
for the computations of $R/G$ and $A\# G$, respectively.

As is easily seen, the quiver presentation of $T(M)$ is given as follows.

\begin{prp}
Let $(Q^R, I^R)$ and $(Q^S, I^S)$ be finite bound quivers of
$R$ and $S$ with $Q^R = (Q^R_0, Q^R_1, s^R, t^R)$,
$Q^S = (Q^S_0, Q^S_1, s^S, t^S)$, respectively.
Let $\Ph^R \colon \k Q^R \to R$,\ $\Ph^S \colon \k Q^S \to S$
be algebra morphisms with $\Ker \Ph^R = I_R,\, \Ker \Ph^S = I_S$.
Set $J_R,\, J_S$ to be the Jacobson radicals of $R,\, S$, respectively.
For each $x \in Q^R_0,\, y \in Q^S_0$,
let $\{m_{y,x}^{(i)} \mid 1 \le i \le d_{y,x}\} \subseteq \Bi(y,M,x)$
such that the residue classes of its elements form a basis
of $\Bi(y,(M/(J_S M + MJ_R)),x)$.
Thus $d_{y,x} = \dim \Bi(y,(M/(J_S M + MJ_R)),x)$.
Then the quiver $Q = (Q_0, Q_1, s, t)$ of $T(M)$ is defined as follows.

$Q_0:= Q^R_0 \sqcup Q^S_0$.
$Q_1:= Q^R_1 \sqcup Q^S_1 \sqcup Q^M_1$,
where $Q^M_1$ is the set of symbols $\al_{y,x}^{(i)}$
for all $x \in Q^R_0,\, y \in Q^S_0,\ 1 \le i \le d_{y,x}$.
For each $\al \in Q_1$,
$$
s(\al):=
\begin{cases}
s^R(\al) & (\al \in Q^R_1),\\
s^S(\al) & (\al \in Q^S_1),\\
x & (\al = \al_{y,x}^{(i)} \in Q^M_1),
\end{cases}\ 
t(\al):=
\begin{cases}
t^R(\al) & (\al \in Q^R_1),\\
t^S(\al) & (\al \in Q^S_1),\\
y & (\al = \al_{y,x}^{(i)} \in Q^M_1).
\end{cases}
$$
Define an algebra morphism $\Phi \colon \k Q \to T(M)$
as follows:
For each $x \in Q_0$, the trivial path $e_x$ is sent as
$$
\Ph(e_x):=
\begin{cases}
\Ph^R(e_x) & (x \in Q^R_0),\\
\Ph^S(e_x) & (x \in Q^S_0);
\end{cases}
$$
and for each $\al \in Q_1$,
$$
\Ph(\al):=
\begin{cases}
\Ph^R(\al) & (\al \in Q^R_1),\\
\Ph^S(\al) & (\al \in Q^S_1),\\
m_{y,x}^{(i)} & (\al = \al_{y,x}^{(i)} \in Q^M_1).
\end{cases}
$$
Then by setting $I:= \Ker \Ph$, $T(M)$ is presented by
the bound quiver $(Q, I)$.
\end{prp}

The following are easy to verify by definitions.

\begin{lem}
$
T(M)/G \iso \bmat{R/G & 0 \\M/G & S/G} =: T(M/G).
$
\qed
\end{lem}

\begin{lem}
$
T(M')\# G \iso \bmat{A\# G & 0 \\M'\# G & S\# G} =: T(M'\# G).
$
\qed
\end{lem}

\subsection{Example}

Throughout the rest of this section, let $G=\ang{g \mid g^2 = 1}$ be the
cyclic group of order 2.

\subsubsection{$G$-categories $R$ and $S$}

Consider the Brauer tree algebras $R$ and $S$
given by the following quivers with relations, respectively:
$$
\begin{aligned}
&Q_R=
\begin{tikzcd}
1 & 2 & 2' &  1'
\Ar{1-1}{1-2}{"\be_1", bend left=20pt}
\Ar{1-2}{1-3}{"\al_2", bend left=20pt}
\Ar{1-3}{1-4}{"\be_{2'}", bend left=20pt}
\Ar{1-2}{1-1}{"\be_2", bend left=20pt}
\Ar{1-3}{1-2}{"\al_{2'}", bend left=20pt}
\Ar{1-4}{1-3}{"\be_{1'}", bend left=20pt}
\end{tikzcd};
\\
&\be_1\be_2\be_1 = 0,\, \be_{1'}\be_{2'}\be_{1'} = 0,\,
\al_2\be_1 = 0,\, \be_{2'}\al_2 = 0,\,
\be_2\al_{2'} =0,\, \al_{2'}\be_{1'} = 0,\\
&\al_{2'}\al_2 = \be_1\be_2,\,
\al_2\al_{2'} = \be_{1'}\be_{2'}
\end{aligned}
$$
and
$$ 
\begin{aligned}
Q_S=
\begin{tikzcd}
&2\\
1 && 1'\\
& 2'
\Ar{2-1}{1-2}{"a_1"}
\Ar{1-2}{2-3}{"a_2"}
\Ar{3-2}{2-1}{"a_{2'}"}
\Ar{2-3}{3-2}{"a_{1'}"}
\end{tikzcd};
\\
\text{paths of length 5 are zero.}
\end{aligned}
$$
Define a $G$-action on $R$ (resp.\ $S$) by the
the automorphism $X_g$ (resp.\ $Y_g$) of $R$ (resp.\ $S$)
that exchanges vertices
$x$ and $x'$ for $x = 1, 2$.

\subsubsection{$G$-graded categories $A$ and $B$}
By taking $?/G$, skeletons $A$ and $B$ of $R/G$ and
$S/G$ turns out to be $G$-graded Brauer tree algebras
given by the weighted quivers $(Q_A, W_A)$ and
$(Q_B, W_B)$ with relations (see \cite[Definition 1.6]{Asa18} or \cite[Definition 6.2.9]{Asa-book})
presented by
$$
\begin{aligned}
&Q_A=
\begin{tikzcd}
\Nname{1}1 & &\Nname{2}2 
\Ar{1}{2}{"\be_1", bend left=20pt}
\Ar{2}{1}{"\be_2", bend left=20pt}
\Ar{2}{2}{"\al_2", "(2)"', loop, in=-30, out=30, distance=8mm}
\end{tikzcd},\ 
W_A(\al_2) = g,\, W_A(\be_1) = W_A(\be_2) = 1;
\\
&\be_1\be_2\be_1 = 0,\, 
\al_2\be_1 = 0,\,\be_2\al_2 = 0,\,
\al_2^2 = \be_1\be_2,\,
\end{aligned}
$$
and 
$$
\begin{aligned}
&Q_B=
\begin{tikzcd}
1 & &2 
\Ar{1-1}{1-3}{"a_1", bend left=20pt}
\Ar{1-3}{1-1}{"a_2", bend left=20pt}
\Ar{1-1}{1-3}{"(2)" description, draw=none}
\end{tikzcd},\ 
W_B(a_1) = 1,\, W_B(a_2) = g;
\\
&\text{paths of length 5 are zero},
\end{aligned}
$$
where (2) in the cycle of each quiver denotes the multiplicity 2 of
the exceptional vertex corresponding to the cycle.

\subsubsection{$G$-invariant $S$-$R$-bimodule $M$}
Define an $S$-$R$-bimodule $M$ so that $T(M)$ is presented by the quiver 
\[
\begin{tikzcd}[row sep=30pt]
1 & 2 & 2' &  1'\\
1 & 2 & 1' &  2'
\Ar{1-1}{1-2}{"\be_1", bend left=15pt}
\Ar{1-2}{1-3}{"\al_2", bend left=15pt}
\Ar{1-3}{1-4}{"\be_{2'}", bend left=15pt}
\Ar{1-2}{1-1}{"\be_2", bend left=15pt}
\Ar{1-3}{1-2}{"\al_{2'}", bend left=15pt}
\Ar{1-4}{1-3}{"\be_{1'}", bend left=15pt}
\Ar{1-1}{2-1}{"m_1" '}
\Ar{1-2}{2-2}{"m_2" '}
\Ar{1-3}{2-4}{"m_{2'}"' near start}
\Ar{1-4}{2-3}{"m_{1'}" near start}
\Ar{2-1}{2-2}{"a_1"}
\Ar{2-2}{2-3}{"a_2"}
\Ar{2-3}{2-4}{"a_{1'}"}
\Ar{2-4}{2-1}{"a_{2'}", bend left=15pt}
\end{tikzcd}
\]
with the Brauer quiver relations for $R$ and $S$ and relations
$$
\left\{
\begin{aligned}
&a_1m_1-m_2\be_1=0,\\
&m_1 \be_2 + a_{2'}a_{1'}a_2 m_2 - a_{2'} m_{2'} \al_2 = 0,
\end{aligned}
\right.\quad
\left\{
\begin{aligned}
&a_{1'}m_{1'}-m_{2'}\be_{1'}=0,\\
&m_{1'} \be_{2'} + a_{2}a_{1}a_{2'} m_{2'} - a_{2} m_{2} \al_{2'} = 0.
\end{aligned}
\right.
$$
We can define a $G$-invariant structure $\ph$ of $M$ by $\ph_g$,
which exchanges $m_i$ and $m_{i'}$ for $i = 1,\, 2$.

\subsubsection{$G$-graded $B$-$A$-bimodule $M'$}
Then the bimodule
${}_{S/G}M/G_{R/G} \simeq {}_BM'_A$ is expressed by $T(M/G)$,
which can be computed as in \cite{Asa-Kim},
where $T(M')$ is a skeleton of $T(M/G)$ defined as its full subcategory consisting of the objects without the prime sign.
Then we see that the bimodule
${}_BM'_A$
is expressed by $T(M')$
that is presented by the quiver
\[
\begin{tikzcd}[column sep=50pt, row sep=40pt]
1 & 2 \\
1 & 2 
\Ar{1-1}{1-2}{"\be_1", bend left=20pt}
\Ar{1-2}{1-1}{"\be_2", bend left=20pt}
\Ar{2-1}{2-2}{"a_1", bend left=20pt}
\Ar{2-2}{2-1}{"a_2", bend left=20pt}
\Ar{2-2}{2-1}{"(2)" description, draw=none}
\Ar{1-1}{2-1}{"m_1" '}
\Ar{1-2}{2-2}{"m_2"}
\Ar{1-2}{1-2}{"\al_2", "(2)"', loop, in=-30, out=30, distance=8mm}
\end{tikzcd}
\]
with the Brauer quiver relations for $A$ and $B$ and relations
$a_1 m_1 - m_2 \be_1=0,\,
m_1 \be_2 + a_2 a_1 a_2 m_2 - a_2 m_2 \al_2=0$,
where the $G$-grading is defined by the weight $W$ given by
$$
W(\al):=
\begin{cases}
W_A(\al) & (\al \in Q^A_1),\\
W_B(\al) & (\al \in Q^B_1),\\
1 & (\al \in \{m_1, m_2\})
\end{cases}
$$
for all arrows $\al$.

In this case, as is easily seen,
we have $M'\# G \iso M$ and $T(M')\# G \iso T(M)$.

\subsubsection{$G$-graded $A$-$B$-bimodule $N'$}
Define an $A$-$B$-bimodule $N'$ by setting
$N':= \Hom_{A\op}(M', A)$.
Recall that a finite dimensional $\k$-algebra $A$ is said to be symmetric if $A$ and its $\k$-dual $ \Hom_\k(A, \k)$ are isomorphic as $A$-$A$-bimodules.
In our case, since $A$ is a symmetric algebra, we have
$$
\Hom_A(\blank, A) \iso \Hom_A(\blank, \Hom_\k(A, \k))
\iso \Hom_\k(\blank \ox_A A, \k) \iso \Hom_\k(\blank, \k),
$$
where $\Hom_A$ denotes the Hom-functor of right $A$-modules.
In particular, we have
$\Hom_A(M', A) \iso \Hom_\k(M', \k)$
as $A$-$B$-bimodules.
Then $N'$ can be seen as a $G$-graded $A$-$B$-bimodule as follows:
Define
$$
(N')^a:= \Hom_\k((M')^{a\inv}, \k)
$$
for all $a \in G$.
Then $N' \iso \Ds_{a \in G}(N')^a$ as $\k$-modules.
Indeed, since $G$ is a {\it finite group}, we have
$$
N' \iso \Hom_{\k}(\Ds_{a\in G}(M')^a, \k) \iso \Ds_{a\in G}\Hom_{A\op}((M')^a, \k) = \Ds_{a\in G}(N')^a.
$$
We verify that this decomposition of $N'$ gives a
$G$-graded $A$-$B$-bimodule structure.
Let $a, b, c \in G$ and $u \in A^a$, $v \in B^b$,
$f \in (N')^c = \Hom_\k((M')^{c\inv},k)$.
It is enough to show that
$ufv \in (N')^{acb} = \Hom_\k((M')^{(acb)\inv},\k)$.
For each $x \in (M')^{(acb)\inv}$,
we have $(ufv)(x) = f(vxu)$,
where since $vxu \in (M')^{b(acb)\inv a} = (M')^{c\inv}$,
$f(vxu) \in \k$ is defined.
Hence $ufv \in \Hom_\k((M')^{(acb)\inv},\k)$.

\begin{prp}
\label{prp:Ggr-st-eq-M}
The pair $(M', N')$ induces a $G$-graded stable equivalence
of Morita type between $A$ and $B$.
\end{prp}

\begin{proof}
We denote by $e_i$ (resp.\ $f_i$) the primitive idempotents corresponding to
the vertices $i \in \{1, 2\}$ in $Q_A$ (resp.\ $Q_B$).
In our case, since $G$ is an abelian group, the left grading shift
and the right grading shift by an $a \in G$ are the same,
and we set $\si_a:= \ro_a = \la_a$ for all $a \in G$. 

It is not hard to check that 
the $G$-graded $B$-$A$-bimodule $M'$ as a representation of $A\op$
in the category $\prj B$ is given as follows:
$$
({}_BM')_A : \quad
\begin{tikzcd}[ampersand replacement=\&, column sep=60pt]
Bf_1 \& Bf_2 \ds \si_g(Bf_2)
\Ar{1-2}{1-1}{"\cdot\sbmat{\cdot a_1\\0}" ', shift right=2pt}
\Ar{1-1}{1-2}{"{\cdot[\cdot(-a_2a_1a_2),\, \cdot a_2]}" ', shift right=2pt}
\Ar{1-2}{1-2}{"\cdot \sbmat{0&\cdot f_2\\
\cdot (-a_1a_2a_1a_2) & \cdot a_1a_2}", loop, in=-30, out=30, distance=60pt}
\end{tikzcd},
$$
where morphisms as matrices act from the right and all entries of
matrices are given by the right multiplication of elements in $B$, which is expressed by writing ``$\cdot$'' on the left of matrices and elements of $B$.

On the other hand, the $G$-graded $B$-$A$-bimodule $M'$ as a representation of $B$
in the category $\prj A\op$ is given as follows:
$$
{}_B(M'_A) : \quad
\begin{tikzcd}[ampersand replacement=\&, column sep=60pt]
e_1A \ds \si_g(e_2A) \& e_2A \ds \si_g(e_2A)
\Ar{1-2}{1-1}{"\sbmat{\be_1\cdot& 0\\0& e_2\cdot}\cdot" ', shift right=2pt}
\Ar{1-1}{1-2}{"\sbmat{0&-\be_2\cdot\\
e_2\cdot& \al_2\cdot}\cdot" ', shift right=2pt}
\end{tikzcd},
$$
where morphisms as matrices act from the left and all entries of
matrices are given by left multiplication of elements in $A$, which is expressed by writing ``$\cdot$'' on the right of matrices and elements of $A$.
Hence, the $G$-graded $A$-$B$-bimodule $N'$ as a representation of $B$
in the category $\prj A$
is given as follows:
$$
\begin{tikzcd}[ampersand replacement=\&, column sep=60pt]
({}_AN')_B: \quad
Ae_1 \ds \si_g(Ae_2) \& Ae_2 \ds \si_g(Ae_2)
\Ar{1-2}{1-1}{"\cdot\sbmat{\cdot\be_1& 0\\0& \cdot e_2}" ', shift right=2pt}
\Ar{1-1}{1-2}{"\cdot\sbmat{0&\cdot(-\be_2)\\
\cdot e_2& \cdot\al_2}" ', shift right=2pt}
\end{tikzcd}.
$$
Therefore, the $A$-$A$-bimodule $N' \ox_B M'$ as a representation of $A\op$ in the category $\prj A$ is
given as follows (note that $g^2 = 1$):
$$
\begin{tikzcd}[ampersand replacement=\&, column sep=80pt]
Ae_1\ds \si_g(Ae_2) \& (Ae_2\ds \si_g(Ae_2)) \ds (\si_g(Ae_2)\ds Ae_2)
\Ar{1-2}{1-1}{"\cdot\sbmat{\cdot \be_1& 0\\0&\cdot e_2\\0& 0\\0&0}" ', shift right=2pt}
\Ar{1-1}{1-2}{"{\cdot\sbmat{\cdot\be_2 & 0 &0 & \cdot(-\be_2)\\\cdot(-\al_2)&0&\cdot e_2& \cdot\al_2}}" ', shift right=2pt}
\Ar{1-2}{1-2}{"\cdot \sbmat{0&0 &\cdot e_2&0 \\0& 0&0&\cdot e_2\\
\cdot\be_1\be_2 &0&0& \cdot(-\be_1\be_2)\\\cdot(-\al_2)&0&\cdot e_2&\cdot\al_2}", loop, in=-15, out=15, distance=75pt}
\end{tikzcd}.
$$
It is not hard to check that $N' \ox_B M' \iso A \ds P$
as $G$-graded $A$-$A$-bimodules
(note that $\al_2$ and $a_2$ have degree $g$, and the others degree $1$),
where $P:= Ae_2 \ox_\k \si_g(e_2 A)$ is a $G$-graded projective $A$-$A$-bimodule.

Similarly,
the $B$-$B$-bimodule $M' \ox_A N'$ as a representation
of $B$ in the category $\prj B$ is given as follows:
$$
\begin{tikzcd}[ampersand replacement=\&, column sep=120pt]
Bf_1\ds (\si_g(Bf_2)\ds Bf_2) \& (Bf_2\ds \si_g(Bf_2)) \ds (\si_g(Bf_2)\ds Bf_2)
\Ar{1-2}{1-1}{"\cdot\sbmat{\cdot a_1& 0&0\\0&0&0\\0& \cdot f_2&0\\0&0&\cdot f_2}" ', shift right=2pt}
\Ar{1-1}{1-2}{"{\cdot\sbmat{0&0&\cdot a_2a_1a_2 & \cdot(-a_2)\\\cdot f_2&0&0&\cdot f_2\\0&\cdot f_2& \cdot(-a_1a_2a_1a_2)&\cdot a_1a_2}}" ', shift right=2pt}
\end{tikzcd}.
$$
It is not hard to check that $M' \ox_A N' \iso B \ds Q$
as  $G$-graded $B$-$B$-bimodules,
where $Q:= Bf_2 \ox_\k \si_g(f_2 B)$ is a $G$-graded projective $B$-$B$-bimodule.
\end{proof}


\end{document}